\documentstyle[amscd,amssymb,epsf]{amsart}
\newtheorem{thm}{Theorem}[section]
\newtheorem{lem}[thm]{Lemma}
\newtheorem{conj}[thm]{Conjecture}
\newtheorem{prop}[thm]{Proposition}
\newtheorem{cor}[thm]{Corollary}

\theoremstyle{definition}
\newtheorem{defn}[thm]{Definition}

\newtheorem{rem}[thm]{Remark}
\newtheorem{nota}[thm]{Notation}
\newtheorem{conve}[thm]{Convention}

\makeatletter
\def\theequation{\@arabic\c@equation}
\makeatother
\numberwithin{equation}{thm}

\def\intt{\int\limits}

\def\U{{\cal U}}
\def\spec{{\hbox{\rom{Spec}\,}}}
\def\ppt{{\Bbb P}_2^*}
\def\Pt{{\Bbb P}_2}
\def\M{{\cal M}}
\def\mmm{{m}}
\def\pt{\Pt}
\def\SS{{\cal S}}
\def\eps{\varepsilon}
\def\R{{\Bbb R}}
\def\pp{{\cal P}}
\def\A{{\cal A}}
\def\G{{\cal G}}
\def\V{{\cal V}}
\def\Q{{\Bbb Q}}
\def\Z{{\Bbb Z}}
\def\N{{\Bbb N}}
\def\P{{\Bbb P}}
\def\C{{\Bbb C}}
\def\I{{\cal I}}
\def\L{{\tilde L}}
\def\X#1#2{X(#1,#2)}
\def\id{{\hbox{\rom{id}}}}
\def\Coh{{\cal Coh}}
\def\II{{\cal I}}
\def\F{{\cal F}}
\def\E{{\cal E}}
\def\OO{{\cal O}}
\def\oo{\OO}
\def\PP{{\cal P}}
\def\eps{\varepsilon}
\def\Ext{{\hbox{\rom{Ext}}}}
\def\Hom{{\hbox{\rom{Hom}}}}
\def\EXT{{\cal Ext}}
\def\HOM{{\cal Hom}}
\def\hxn{{X^{[n]}}}
\def\hsn{{S^{[n]}}}
\def\xn{{X(n)}}
\def\xtn{{\widetilde \xn}}
\def\uberschrift#1{\bigskip\bigskip {\bf #1}\medskip}
\def\Theor#1#2{\medskip \noindent{{\bf Theorem~#1.~~}}{\it #2}\medskip}
\def\Thm#1#2{\medskip \noindent{{\bf Theorem~#1.~~}}{\it #2}\medskip}
\def\Prop#1#2{\medskip \noindent{{\bf Proposition~#1.~~}}{\it #2}\medskip}
\def\Cor#1#2{\medskip \noindent{{\bf Corollary~#1.~~}}{\it #2}\medskip}
\def\Lemma#1#2{\medskip \noindent{{\bf Lemma~#1.~~}}{\it #2}\medskip}
\def\Def#1{\medskip \noindent{{\bf Definition~#1.~~}}}
\def\Not#1{\medskip \noindent{{\bf Notation~#1.~~}}}
\def\Rmk#1{\medskip \noindent{{\bf Remark~#1.~~}}}
\def\Rem#1{\Rmk{#1}}
\def\Ex#1{\medskip \noindent{{\bf Example~#1.~~}}}
\def\Parag#1#2{\bigskip \noindent {\S}{\bf #1.~#2.} \medskip}
\def\Pf{\demo{Proof}}
\def\acapo{\vskip 1pt \noindent}
\def\whsq{\vbox to 5.8pt
{\offinterlineskip\hrule
\hbox to 5.8pt{\vrule height
5.1pt\hss\vrule height 5.1pt}\hrule}}
\def\blsq{\vrule height 8pt width 6pt}
\def\Qed{{\hfill {\whsq}}\enddemo}
\def\Cvd{\Qed}
\def\ppp{{\Bbb P}}
\def\maps{\longrightarrow}
\def\oo{{\cal O}}
\def\ii{{\cal I}}
\def\diagramm#1{
\def\normalbaselines{\baselineskip18pt\lineskip7pt\lineskiplimit7pt}
 \matrix{#1}}

\def\mapr#1{\smash{
 \mathop{\longrightarrow}\limits^{#1}}}

\def\mapsr#1{{
 \mathop{\longrightarrow}\limits^{#1}}}

\def\mapinto#1{\smash{
 \mathop{\longrighthook}\limits^{#1}}}

\def\mapl#1{\smash{
 \mathop{\longleftarrow}\limits^{#1}}}

\def\mapd#1{\Big\downarrow
\rlap{$\vcenter{\hbox{$\scriptstyle#1$}}$}}
\def\mapu#1{\Big\uparrow
\rlap{$\vcenter{\hbox{$\scriptstyle#1$}}$}}

\def\mapsec#1#2{\mapdownl#1\mapup#2}

\def\mapne#1{\nearrow
\rlap{$\vcenter{\hbox{$\scriptstyle#1$}}$}}
\def\mapse#1{\searrow
\rlap{$\vcenter{\hbox{$\scriptstyle#1$}}$}}
\def\mapnw#1{\nwarrow
\rlap{$\vcenter{\hbox{$\scriptstyle#1$}}$}}
\def\mapsw#1{\swarrow
\rlap{$\vcenter{\hbox{$\scriptstyle#1$}}$}}

\def\relgl{\mathrel{\smash=}}
\def\longleich{\relgl\joinrel\relgl\joinrel\relgl}

\def\sn{\Hilb^n(S)}
\def\snn{{S^{[n-1,n]}}}
\def\sym{S^{(n)}}
\def\H{{\cal H}}
\def\C{{\Bbb C}}
\def\R{{\Bbb R}}
\def\P{{\Bbb P}}
\def\N{{\Bbb N}}
\def\PP{{\cal P}}
\def\Z{{\Bbb Z}}
\def\Q{{\Bbb Q}}
\def\T{{\cal T}}
\def\eps{\varepsilon}
\def\Oplus{\textstyle\bigoplus\limits}
\def\phi{\varphi}
\def\cinf{C^{\infty}}
\def\Symg{\frak S}

\def\({\left(}
\def\){\right)}
\def\pB{{\cal B}}
\def\sB{{\cal B}^\prime}
\def\Supp{\text{Supp}\,}

\def\part{P(n)}
\def\bb{{\text{\rom{\bf b}}}}
\def\hsb{S^{(\bb)}}

\def\EE{{\text{\rom{\bf E}}}}

\def\al{\alpha}
\def\Ga{\Gamma}
\def\bGa{\bar \Gamma}
\def\ga{\gamma}
\def\peta{\bar \eta}

\def\ZZ{{\cal Z}}

\def\pt{{pt}}
\def\al{\alpha}

\def\Sn{{S^{(n)}}}
\def\Sm{{S^{(m)}}}
\def\<{\langle}
\def\>{\rangle}
\def\Hilb{{\text{\rom{Hilb}}}}
\def\W{{\cal W}}
\def\Tor{{\hbox{\rom{Tor}}}}
\def\Ext{{\hbox{\rom{Ext}}}}

\def\stil#1{\tilde S^#1}
\def\pigr{\pi}

\def\gi{g}
\def\Fi{\phi}
\def\ti{\tilde}
\def\alb{{\bar\al}}

\begin{document}

\title[Variation of Donaldson invariants]{Variation of moduli spaces
and Donaldson invariants under change of
polarization}
\author{Geir Ellingsrud}
\address{Mathematical Institute\\University of Oslo\\P.~O.~Box~1053\\
         N--0316 Oslo, Norway}
\email{ellingsr@@math.uio.no}
\author{Lothar G\"ottsche}
\address{Max--Planck--Institut f\"ur Mathematik\\Gottfried--Claren--Stra\ss e
26\\
D-53225 Bonn, Germany}
\email{lothar@@mpim-bonn.mpg.de}

\maketitle\

%\begin{abstract}
%\end{abstract}

\section{Introduction}

Let $S$ be a smooth projective surface over the complex numbers and let
 $c_1\in H^2(S,\Z)$ and
$c_2\in H^4(S,\Z)$ be two classes.
For an ample divisor $H$   on $S$, one can study the
moduli space $M_H(c_1,c_2)$ of $H$-semistable
torsion-free sheaves $E$ on $S$ of rank $2$
with $c_1(E)=c_1$
and $c_2(E)=c_2$. We want to study the change of $M_H(c_1,c_2)$
under variation of  $H$. It is known that the ample cone of $S$ has a
chamber structure, and that $M_H(c_1,c_2)$ depends only on the
 chamber containing  $H$.
In this article we will  try to understand how $M_H(c_1,c_2)$ changes,
 when $H$ passes
through a wall separating two chambers. The set-theoretic changes of
the subspace
consisting  of locally free sheaves and of
$M_H(c_1,c_2)$ have been treated in \cite{Q1} and \cite{Go1} respectively.
We show that the change of the moduli space when $H$ passes  through a wall,
can be expressed as a sequence of operations similar to a flip.
In fact the  moduli spaces at each step can be identified as moduli spaces
 of torsion-free sheaves with a
suitable parabolic structure of length $1$.
We  assume that either the geometric genus $p_g(S)$ is $0$ or that $K_S$ is
trivial. We
shall also  make suitable hypotheses on the wall, and walls fulfilling this
condition
we call good. This assumption is reasonably weak if the Kodaira
dimension of $S$ is at most $0$, but gets stronger if
e.g., $S$ is of general type.
When the polarization passes through a good wall, each of the steps above
is realized by  a smooth blow-up
along a projective bundle
over a product of Hilbert schemes  of points on $S$, followed by a smooth
 blow-down of the exceptional divisor in another direction.

The change of moduli spaces can  be
viewed as a change of GIT quotients, treated in \cite{Th2} and \cite{D-H}.
These results  could  in principle be applicable, although it
would still take quite some work to do so.
We have however chosen a more direct approach via elementary
transforms of universal families, which is more in the spirit of \cite{Th1},
and which also immediately gives the change of the universal sheaves
 needed for the computation of Donaldson invariants.

In the case that $K_S$ is trivial, i.e., $S$ is an abelian or a $K3$ surface,
we see that the change of $M_H(c_1,c_2)$, when $H$ passes through a wall,
is given by elementary transformations of symplectic varieties.

In the case that $p_g(S)=q(S)=0$ we   use these results in
 order to compute the change of the Donaldson
polynomials under change of polarisation.
The Donaldson polynomials of a $C^\infty$-manifold $M$ of dimension $4$
are defined using a Riemannian metric on $S$,
but in case $b^+(M)>1$ they are known to be independent of the metric, as
long as it is generic.
In case $b_+(M)=1$, (which for an algebraic surface $S$ corresponds to
$p_g(S)=0$), the invariants have been introduced and studied by
Kotschick in \cite{Ko}.
In \cite{K-M} Kotschick and Morgan  show that the invariants only
depend on the chamber
of the period point of the metric in the positive cone of $H^2(M,\R)$.
They also compute  the lowest order term of the change and  conjecture
the shape of a formula for the change.

The case we are studying corresponds to $M$ being an
algebraic surface $S$ with $p_g(S)=q(S)=0$ and a wall
lying in the ample cone, in addition we assume
that the wall is good.

In a first step we compute the change of the Donaldson invariants
 in terms of
natural cohomology classes on Hilbert schemes of points on $S$ and then we
 use  some computations
in the cohomology rings of these Hilbert schemes  to determine
 the six lowest order  terms of the change of the Donaldson invariants
explicitly. The results are compatible with the conjecture of \cite{K-M}
(which in particular predicts that three of the terms above are zero.

Parallelly and independently similar results to ours have been
obtained by other authors.
Matsuki and Wentworth show in \cite{M-W}  that the change of moduli spaces
of torsion-free sheaves of arbitrary rank on a projective variety
under change of polarisation can be described as a sequence of flips.
In \cite{F-Q} Friedman and Qin  obtain  very similar results to ours.

\section{Background material}

In this paper let  $S$ be a projective surface over $\C$. By  the Neron-Severi
group $NS(S)$ of $S$ we mean the group of divisors modulo homological
equivalence, i.e., the image of   $Div(S)$ in $H^2(S,\Z)$ under the map sending
the class of a divisor $D$ to its fundamental cycle $[D]$. Let
$Div^0(S)$ be its kernel. Let
$c_1\in H^2(S,\Z)$ and $c_2\in H^4(S,\Z)=\Z$ be elements which will be fixed
throughout the paper.

Let $H$ be a polarization of $S$.
As we mostly shall consider stability and semistability in the sense of
Gieseker
and Maruyama we  shall write
$H$-stable (resp. $H$-semistable) instead of Gieseker stable  (resp.
semistable)
with respect to $H$ and $H$-slope stable (resp. $H$-slope  semistable)
instead of
stable (resp. semistable) with respect to $H$ in the sense of Mumford-Takemoto.
Denote by   $M_H(c_1,c_2)$
 the moduli space of H-semistable torsion-free sheaves $E$ on $S$ of
rank $2$ with $c_1(E)=c_1$ and $c_2(E)=c_2$ and
$M_H^s(c_1,c_2)$ the open subscheme of $M_H^s(c_1,c_2)$ of stable sheaves.
Let $Spl(c_1,c_2)$ be the moduli space of simple
torsion-free sheaves with $c_1(E)=c_1$ and $c_2(E)=c_2$ (see \cite{A-K}).

\begin{nota}
For a sheaf $\F$ on a scheme $X$ and a divisor $D$ let $\F(D):=\F\otimes
\oo_X(D)$ .

Many of our arguments will take place over products $S\times X$, where $X$ is a
scheme.  We shall denote by $p:S\times X\maps S$ and $q_X:S\times X\maps X$ the
two projections and if there is no danger of confusion, we will  drop the index
$X$. For  a divisor
$D$  on $X$ we denote $D_S:=q_X^*(D)$. For a sheaf $\F$ on $S\times X$
and a divisor or divisor
class, $D$  on $S$ we denote by $\F(D)$  the sheaf
$\F(p^*(D))$.

If $X$ is a smooth variety of dimension $n$, we denote the cup product  of two
elements $\alpha$ and $\beta$ in $H^*(X,\Z)$ by $\alpha\cdot\beta$ and the
degree of a class
$\alpha\in H^{2n}(X,\Z)$ by $\int_X\alpha$. For $\alpha,\beta\in H^2(S,\Z)$ let
$\<\alpha\cdot\beta\>:=\int_S\alpha\cdot\beta$. We write
$\alpha^2$ for $\<\alpha\cdot\alpha\>$ and, for $\gamma\in H^2(S,\Z)$, we put
$\<\alpha,\gamma\>:=\<\alpha\cdot \check\gamma\>$, where
$\check\gamma$ is the Poincar\'e dual of $\gamma$.
\end{nota}

\begin{conve} \label{convent}
If $Y,X$ are schemes and there is a "canonical" map $f:X\maps Y$,
then for a cohomology class $\alpha\in H^*(Y,\Z)$ (resp. for a vector bundle
$E$ on $Y$) we will very often also denote the pull-back via $f$
by $\alpha$ (resp. $E$).
\end{conve}

\begin{defn} \cite{OG2}\label{defpseudofam} Let $B$ be a scheme.
A family of sheaves, $\F$, on $S$ parametrized by $B$ is a $B$-flat sheaf on
$S\times B$. Two families of sheaves $\F$ and $\G$ on $S$ parametrized by
$B$
are called equivalent if there exists an isomorphism
$\F\simeq \G\otimes q_B^*M$, for some line bundle $M$  on $B$.
 Let
$(B_j)_{j\in J}$ be an \'etale cover of $B$ by  schemes.
Assume that for each $j\in J$ there is a family $\F_j$ of sheaves on
$X$ parametrized by $B_j$, and that for each pair $k,l\in J$
the pullbacks of $\F_k $ and $\F_l$ to $B_k\times_B B_l$ are equivalent.
Then we will say that the above data defines a pseudo-family of sheaves on $S$
parametrized by $B$. We will denote it by $\F$.

It is clear what is meant by a map of pseudo-families and by
two pseudo-families being equivalent.
\end{defn}

The main reason to introduce pseudo-families is that the moduli space
$M^s_H(c_1,c_2)$ does not always carry a universal
family of sheaves, but  there will always be a universal pseudo-family.

By the universal property of $M^s_H(c_1,c_2)$ a pseudo family of $H$-stable
torsion-free sheaves $E$ on $S$ with $c_1(E)=c_1$,
$c_2(E)=c_2$ parametrized by $B$ gives rise to a morphism $B\maps
M^s_H(c_1,c_2)$.

\bigskip

{\bf Walls and chambers for torsion-free sheaves}

We now recall some results about walls and chambers from \cite{Q1},
\cite{Q2} and \cite{Go1}.

\begin{defn}\label{defwall}(for the first part see \cite{Q1} Def I.2.1.5)
 Let $C_S$ be the ample cone in $NS(S)\otimes \R$.
For $\xi\in NS(S)$ let
$$W^\xi:=C_S\cap\big \{ x\in NS(S)\otimes\R \bigm| \<x\cdot\xi\>=0\big\}.$$
We shall call $W^\xi$ a wall of type $(c_1,c_2)$, and say that it is
defined by  $\xi$ if the following conditions are satisfied:
\begin{enumerate}
\item  $\xi+c_1$ is divisible by $2$ in
$NS(S)$,
\item $c_1^2-4c_2\le \xi^2<0$,
\item there is a polarisation $H$ with $\< H\cdot \xi\>=0$.
\end{enumerate}

In particular  $d_\xi:= (4c_2-c_1^2+\xi^2)/4$ is  a nonnegative integer.
An ample divisor $H$ is said to lie in the wall $W$ if $[H]\in W$.
If $D$ is a divisor with $[D]=\xi$, we will also say that
$D$ defines the wall $W$.

A {\it chamber} of type $(c_1,c_2)$ or simply a chamber, is a connected
component of the complement of the union of all the walls of
type
$(c_1,c_2)$. Two different chambers will be said to be
{\it neighbouring chambers} if the intersection of
 their closures contains a nonempty open subset
of a wall.

We will call a wall $W$ {\it good}, if $D+K_S$ is not effective
for any divisor $D$ defining  the wall $W$.
\end{defn}

If $D$ defines a wall, then neither  $D$ nor $-D$ can be effective because
$D$ is orthogonal to an ample divisor.
In particular  every  wall will be  good  if $-K_S$ is
effective or
if $[K_S]$ is a torsion class. More generally, a wall $W$ will be  good
if there exists an
ample divisor
$H$ in $W$ with $\<K_S \cdot H\>\le 0$.

\begin{defn}\label{defenm}
Let $\Hilb^l(S)$ be the Hilbert scheme of subschemes of length $l$ on $S$.
For  $\alpha\in NS(S)$ and  $l\in \Z$, let
$M(1,\alpha,l)$ be the moduli space of rank $1$ torsion-free sheaves
$\I_{Z}(F)$ on $S$
with $c_1(\I_{Z}(F))=[F]=\alpha$, $c_2(F)=length(Z)=l$.
Let
$$T^{n,m}_\xi:=\coprod_{2\alpha=c_1+\xi} M(1,\alpha,n)\times
M(1,c_1-\alpha,m).$$
Let $N_2 \subset NS(S)$ be the subgroup of  $2$-torsion elements.
There is a (noncanonical) isomorphism $$T^{n,m}_\xi\simeq N_2\times\Hilb^n
S\times Div^0(S)\times \Hilb^mS\times Div^0(S),$$
which depends on the choice of an $\alpha\in NS(S)$ with $2\alpha=c_1+\xi$
and on a representative $F$ in $Div(S)$ for $\alpha$.

For any extension
\begin{equation*}\tag{$\epsilon$}0\maps A_1\maps E\maps  A_2\maps
0\end{equation*} where $A_1$ and $A_2$ are  torsion-free rank
one sheaves, we define $\Delta(\epsilon):=\chi(A_1)-\chi(A_2)$. Then if
$\alpha=c_1(A_1)-c_1(A_2)$, the Riemann-Roch theorem gives
$\Delta(\epsilon)=1/2\<(c_1(E)-K_S)\cdot\alpha\>+c_2(A_2)-c_2(A_1)$.
Furthermore for any divisor $D$ we have
$\Delta(\epsilon(D))=\Delta(\epsilon)+\<\alpha\cdot
D\>$, where $ \epsilon(D)$ denotes the extension
$\epsilon$ twisted by the line bundle $\oo(D)$. This follows immediately
from the fact that $c_1(E(D))=c_1(E)+2[D]$.

 Assume that $\xi$ defines a wall of type $(c_1,c_2)$,
and that $n$ and $m$ are nonnegative integers with
$n+m=d_\xi=c_2-(c_1^2-\xi^2)/4$.
Let $\EE^{n,m}_\xi$ be  the set of sheaves lying in nontrivial extensions
\begin{eqnarray}\label{splitting}
&&0\maps \I_{Z_1}(F_1)\maps E\maps  \I_{Z_2}(F_2)\maps 0\end{eqnarray}
 where $(\I_{Z_1}(F_1),\I_{Z_2}(F_2))$ runs through $T^{n,m}_\xi$.
It is easy to see that every sheaf $E\in {\EE}^{n,m}_\xi$ is simple
(\cite{Go1}, lemma 2.3).
Let $$V^{n,m}_\xi:=\EE^{n,m}_\xi\setminus\Big(\bigcup_{l,s}
\EE_{-\xi}^{l,s}\Big).$$
\end{defn}

\begin{nota}\label{defaplus}
Assume  that  $H_+$ and $H_-$ are ample divisors  lying in  neighbouring
chambers separated
by the wall $W$. Then we define
$$A^+(W):=\Big\{\xi\in NS(S)\Bigm| \xi \text{ defines } W \text{ and } \
\<\xi\cdot H_+\> >0\Big\}$$
and $A^-(W):=-A^+(W)$.\end{nota}

The following proposition mostly comprizes some of the results of
\cite{Go1}, that are generalizations of the corresponding results of \cite{Q1},
\cite{Q2} and
 will be important for the rest of the paper.
Note that unlike  \cite{Go1} we assume walls to be defined by
classes in $NS(S)$ and not by numerical equivalence classes, and that we
look at moduli spaces with fixed first Chern class and not with fixed
determinant. The proofs in  \cite{Go1} stay however valid with very few
changes.

\begin{prop} \label{wall}
\begin{enumerate}
\item For $H$ not lying on a wall, $M_H(c_1,c_2)\setminus
M_H^s(c_1,c_2)$ is independent of $H$ and
$M_H(c_1,c_2)$ depends only on the chamber of $H$.

\noindent For the rest of the proposition we assume that we are in the
situation of \ref{defaplus} and  that $\xi\in A^+(W)$.
\item Every $E\in \EE^{n,m}_\xi$ is $H_+$ slope-unstable and the sequence
(\ref{splitting})
is its Harder-Narasimhan filtration with respect to $H_+$.
\item
$\Hom(\I_{Z_1}(F),E)=\C$. Thus,
for $E\in \EE^{n,m}_\xi$, the sequence \ref{splitting} is the
unique extension
$$0\maps \I_{W_1}(F_1)\maps E\maps \I_{W_2}(G_2) \maps 0$$
with $\<(2F_1-c_1)\cdot H_+\> >0$.
\item In particular we see that,  for $\xi,\eta\in A^+(W)$,
the subsets $\EE_\xi^{n,m}$, $\EE_{\eta}^{k,l}$ of
$Spl(c_1,c_2)$ do not intersect, unless
$\xi=\eta$ and $(n,m)=(k,l)$.
\item If $E\in \EE_\xi^{n,m}$ then $E$ is $H_-$-slope stable
if and only if $E\in V_\xi^{n,m}$ and $H_-$-slope unstable otherwise.
\item
On the other hand let $E$ be a torsion-free sheaf with $c_1(E)=c_1$
and $c_2(E)=c_2$, which is $H_-$-semistable and $H_+$-unstable.
Then $E$ is $H_-$-slope stable and $E\in {\EE}^{n,m}_\xi$ for suitable
numbers $n$ and $m$ and
 $\xi\in A^+(W)$.
\end{enumerate}
\end{prop}
\begin{pf}
(1) is (\cite{Go1}, theorem  2.9(1)).  (2) is easy. (3) follows from (2)
and  (\cite{Go1},
lemma  2.3). (4) follows from (2). (5) is (\cite{Go1},
prop 2.5). (6) is  (\cite{Go1},
lemma  2.2).
\end{pf}

\section{Parabolic structures and the passage through a wall}

As mentioned in the previous section, $M_H(c_1,c_2)$ depends only on the
chamber to which $H$ belongs. If $H'$ lies in a
neighbouring chamber to $H$ the moduli space $M_{H'}(c_1,c_2)$  will in
most cases be birational to $M_H(c_1,c_2)$, although  new components do
 occur in some cases. If  the wall
separating the two chambers is good, we will describe
the birational transformation in detail  by giving an explicit sequence
  of blow-ups and blow-downs with smooth
centers which are known.

If the wall is good, but the transformation  is not birational,  our
arguments give a description of the components which are
 added to or deleted from the moduli space.

For the rest of the paper we will assume
that  $H_+$ and $H_-$ are ample divisors  lying in  neighbouring chambers
separated
by the wall $W$, and that $H$ is an ample divisor in the wall
$W$ which lies in the closure of the
chambers containing $H_-$ and $H_+$ respectively
and which does not lie in any other wall. Furthermore
we shall assume that $M=H_+-H_-$ is effective.
By
replacing $H_+$ by a high multiple if necessary,
we can always achieve this.

Our aim is to divide the passage  through a wall into a number of smaller
steps. To this purpose we will   introduce a finer notion of stability.
The starting point is the observation that unlike slope stability,
Gieseker stability
is not invariant under tensorization by a line bundle.

\begin{lem}\label{tensor}
There is a positive integer $n_0$ such that for all $l\ge n_0$ and all
torsion-free rank $2$ sheaves $E$
on $S$ with $c_1(E)=c_1$, $c_2(E)=c_2$
\begin{enumerate}
\item $E$ is $H_-$-stable (resp. semistable) if and only if $E(-lM)$ is
$H$-stable (resp. semistable).
\item $E$ is $H_+$-stable (resp. semistable) if and only if $E(lM)$ is
$H$-stable (resp. semistable).
\end{enumerate}\end{lem}

\begin{pf}
It will be enough to show (1).
As $H_-$ does not lie on a wall, it is easy to see that  $E$ is
$H_-$-(semi)stable if and only if $E(M)$ is.
Also there are only finitely
many $\xi\in NS(S)$ defining the wall $W$.
Therefore lemma \ref{tensor} follows immediately from
lemma \ref{tensor1} and lemma \ref{tensor2} below.\end{pf}

\begin{lem}\label{tensor1}
 \begin{enumerate}
\item Assume $E$ is $H_-$-semistable but  $H$-unstable.
Then $E\in {\EE}^{n,m}_\xi$ for suitable
$n,m$ and $\xi\in A^+(W)$.
\item Assume $E$ is $H_-$-unstable but  $H$-semistable.
Then $E\in {\EE}^{n,m}_{-\xi}$ for suitable
$n,m$ and $\xi\in A^+(W)$.
\end{enumerate}\end{lem}

\begin{pf}
We just prove (1), the proof of (2) being analoguous. By assumption there
is an extension
\begin{equation*}\tag{$\epsilon$}0\maps \I_{Z_1}(F_1)\maps E\maps
\I_{Z_2}(F_2)\maps 0\end{equation*}
with
$\Delta(\epsilon(lH_-))\le 0$ and
$\Delta(\epsilon(lH))>0$ for $l>>0$.
In particular we have
$\<\eta \cdot H_-\>\le 0\le\< \eta \cdot H\>$ where $\eta:=2[F_1]-c_1$.
If  $0<\< \eta\cdot H\>$, there would be  a wall separating $H_-$ and $H$.
So $\< \eta\cdot H\>=0$, and  unless $\eta$ is a torsion
class, it defines a wall in which $H$ lies. As $H$ lies in a unique wall it
must be $W$.  Hence $\eta\in A^+(W)$, and $E\in
{\EE}^{n,m}_{\eta}$.

Assume that $\eta$ is a torsion class. Then $F_1$ and $F_2$ are numerically
equivalent, and it is easily verified that
\begin{equation*}
\Delta(\epsilon(lH_-))=\Delta(\epsilon(lH))
\end{equation*}
which is a contradicition.\end{pf}

\begin{lem}\label{tensor2}
Given $n,m,\xi$. Then there exists an integer $k_0$ such that for all
$k>k_0$ and all $E\in {\EE}^{n,m}_\xi$, the sheaf  $E(-k M)$ is $H$-stable
if and only $E$ is $H_-$-slope stable.
Otherwise  $E$ is both  $H_-$-slope unstable and $H$-unstable.
\end{lem}

\begin{pf}
Let $E\in E^{n,m}_\xi$. Then there is an extension
\begin{equation*}\label{split2}\tag{$\epsilon_1$}0\maps \I_{Z_1}(F_1)\maps
E\maps  \I_{Z_2}(F_2)\maps 0\end{equation*}
 with
$\xi=2[F_1]-c_1$. Assume first that $E$ is $H_-$-slope-stable
and let
\begin{equation*}\tag{$\epsilon_2$}0\maps \I_{Y_1}(G_1)\maps E\maps
\I_{Y_2}(G_2)\maps 0\end{equation*}
be another extension. Put $\eta:=2[G_1]-c_1$. As $E$ is
  $H_-$-slope stable
 we have $\<\eta\cdot H_-\>< 0$, and because there is no wall
 between $H_-$
and $H$, we know that $\<\eta\cdot H\>\le 0$.
For any integers $k$ and $l$
$$\Delta(\epsilon_2(-kM+lH))=\Delta(\epsilon_2)-k\<\eta\cdot
M\>+l\<\eta\cdot H\>.$$
Hence if
$\<\eta\cdot H\><0$ for all extensions $\epsilon_2$ above, then $E(-k M)$
will be  $H$-stable for any $k$.
Assume that  $\<\eta\cdot H\>=0$. By assumption $H$ is contained in a
single wall $W$, so necessarily $\eta\in A^+(W)$.
Hence  by proposition  \ref{wall}(3), we get
$\I_{Z_1}(F_1)=\I_{Y_1}(G_1)$. Therefore it suffices
to see that for $k>>0$
and any $l$ we have the inequality
$$\Delta(\epsilon_1(-kM+lH))<0.$$
Now
$$\Delta(\epsilon_1(-kM+lH))=\Delta(\epsilon_1)-k\<\eta\cdot M\>,$$
which is negative for $k>>0$ as $\<\eta\cdot M\>>0$.

To prove the converse, assume that $E$ is not $H_-$-slope stable.
Then by proposition  \ref{wall}(5) there is  an extension
\begin{equation*}\tag{$\epsilon_3$}0\maps \I_{Y_1}(F_2)\maps E\maps
\I_{Y_2}(F_1)\maps 0.\end{equation*}
 Because
$2[F_2]-c_1=-\xi$ we have
$$\Delta(\epsilon_3(-kM+lH))=\Delta(\epsilon_3)+\<-\xi\cdot (-kM+lH)\>=
\Delta(\epsilon_3)+k\<\xi\cdot M\>>0$$
for $k>>0$\end{pf}

{}From now on until the end of this section we fix $n_0$ as in lemma
\ref{tensor},
and we put $C:= (n_0+1)M$.

\begin{defn}\label{defalstable}
Let  $a$ be a real number between $0$ and $1$. For any torsion-free sheaf
$E$ we define
$$P_a(E)=((1-a)\chi(E(-C))+a\chi(E(C)))/rk(E).$$
A torsion-free sheaf $E$ on $S$
 is called $a$-semistable if and only if  every subsheaf
$E'\subset E$
satisfies
$P_a(E'(lH))\le P_a(E(lH))$ for all $l>>0,$
and it is called $a$-stable if  strict inequality holds.
\end{defn}

In particular, by lemma \ref{tensor}, $E$ is $0$-semistable if
and only if it is
$H_-$-semistable, and it is $1$-semistable if and only if it is
$H_+$-semistable.

For any extension
\begin{equation*}\tag{$\epsilon$}0\maps A_1\maps E\maps  A_2\maps
0\end{equation*}
we define $\Delta_a(\epsilon):=P_a(A_1)-P_a(A_2)$. Then
$\Delta_a(\epsilon)=\Delta(\epsilon)+(2a-1)\<C\cdot\alpha\>$ where
$\alpha=c_1(A_1)-c_1(A_2)$.
Clearly $\Delta_a(\epsilon(D))=\Delta_a(\epsilon)+\<D\cdot\alpha\>$
for any divisor $D$.
A sheaf $E$ is $a$-stable (resp. $a$-semistable) if
$\Delta_a(\epsilon(lH))<0$ (resp. $\le0$) for all $l>>0$
and for any extension $\epsilon$ whose middle term is $E$.

\begin{rem} It is easy to see that $P_a(E(lH))$ is the parabolic Hilbert
polynomial of the parabolic bundle $(E(C),E(-C),a)$,
(i.e.  with a  filtration of length $1$). Therefore $E$ is $a$-semistable
if and only if  $(E(C),E(-C),a)$ is semistable. In
\cite{Ma-Yo} a coarse quasiprojective  moduli space of stable
parabolic sheaves
 with fixed Hilbert polynomial is constructed, and by \cite{Yo} there
exists a projective coarse  moduli space for
$S$-equivalence classes of semistable parabolic sheaves. In particular
there exists a coarse moduli space $M_{a}(c_1,c_2)$ for
$a$-semistable sheaves $E$ on $S$ with $c_1(E)=c_1$ and
$c_2E=c_2$. We denote by $M_{a}^s(c_1,c_2)$ its open subscheme of stable
sheaves.
\end{rem}

\begin{rem} We see that $M_{H_-}(c_1,c_2)$ and
$M_0(c_1,c_2)$ respectively $M_{H_+}(c_1,c_2)$  and $M_1(c_1,c_2)$
are coarse moduli schemes for the same functor and therefore they are
isomorphic.
\end{rem}

\begin{rem}
The same proof as in the case of $H$-stable sheaves shows that
$M_{a}^s(c_1,c_2)$
carries a universal pseudofamily.
One checks easily that every $E\in M^s_{a}(c_1,c_2)$
is simple. As $M^s_{a}(c_1,c_2)$ and $Spl(c_1,c_2)$
both carry universal
pseudofamilies, $\V$ and  $\W$ respectively, there exists a morphism
$f:M^s_{a}(c_1,c_2)\to Spl(c_1,c_2)$ such that $(\id_S\times f)^*(\W)=\V$.
Let $M$ be its image. By the same argument there exists a map
$g:M\to M^s_{a}(c_1,c_2)$, with
$(\id_S\times f)^*(\id_S\times g)^*(\V)=\V$. Hence $f$ is an open
embedding.
In particular and what is the most important thing for us,
 the tangent space to $M^s_{a}(c_1,c_2)$
at a point $E$ is $\Ext^1(E,E)$.
\end{rem}

\begin{defn}\label{defminiwall} For all $a\in
 [0,1]$ let $A^+(a)$ be the set of $(\xi,n,m)\in  A^+(W)\times \Z_{\ge
0}^2$  satisfying
\begin{eqnarray}\label{minicond} n+m&=& c_2-(c_1^2-\xi^2)/4,\\ n-m&=&
\<\xi\cdot (c_1-K_S)\>/2
+(2a-1)\<\xi\cdot[C]\>.\end{eqnarray} A number $a$ is called a {\it
miniwall} if $A^+(a)\ne \emptyset$. A {\it minichamber} is a
 connected component of the complement of the set of all miniwalls in
$[0,1]$. It is clear that there are finitely many
minichambers. Two minichambers are called neighbouring minichambers if
their closures intersect.
\end{defn}

\begin{rem}Note that $A^+(a)$ is the set of all
$\xi,n,m$ with $\xi\in A^+(W)$ for which there exists a (possibly split)
extension
\begin{equation*}\tag{$\epsilon$}0\maps A_1\maps E\maps  A_2\maps
0\end{equation*}
with
$\xi=c_1(A_1)-c_1(A_2)$, $n=c_2(A_1)$, $m=c_2(A_2)$ and
$\Delta_a(\epsilon)=0$.
\end{rem}
\begin{lem}
Let $0\le a_-< a_+\le 1$ and assume that neither $a_-$ nor $a_+$ is a
miniwall. Let $E$ be $a_-$-semistable and
$a_+$-unstable.
 Then there exists a miniwall $a$ between $a_-$ and
$a_+$ and an element $(\xi,n,m)\in A^+(a)$, such that
$E\in {\EE}_{\xi}^{n,m}$.
\end{lem}

\begin{pf} By assumption $E$ is $a_+$-unstable. Hence there is  an extension
\begin{equation*}\tag{$\epsilon$}0\maps A_1\maps E\maps  A_2\maps
0\end{equation*}
 such that for
all $l>>0$ we
have $\Delta_{a_+}(\epsilon(lH))>0$. Putting $\xi:=c_1(A_1)-c_1(A_2)$ and
 using that $E$ is $a_-$-semistable,
we obtain the following inequalities valid for all $l>>0$
$$\Delta_{a_-}(\epsilon(lH))=\Delta_{a_-}(\epsilon)+l\<H\cdot\xi\>\le0<
\Delta_{a_+}(\epsilon(lH))=\Delta_{a_+}(\epsilon)+l\<H\cdot\xi\>.$$
In particular $\<H\cdot\xi\>=0$  and
$\Delta_{a_-}(\epsilon)<0<\Delta_{a_+}(\epsilon)$.

Furthermore $\xi$ is not a torsion class
and  $\xi$ defines a wall on which $H$ is lying,
which therefore must be $W$.
There clearly is an $a$ such that $\Delta_a(\epsilon)=0$.
\end{pf}

\begin{lem} \label{stab} Let $a_-<a_+$ be in neighbouring minichambers
separated by the miniwall $a$.
Let $(\xi,n,m)\in A^+(a)$.
\begin{enumerate}
\item Any $E\in {\EE}_\xi^{n,m}$ is $a_-$-stable, strictly
$a$-semistable
and $b$-unstable for all $b>a$.
\item Any $E\in {\EE}_{-\xi}^{m,n}$ is $a_+$-stable, strictly $a$-semistable
and $b$-unstable for all $b<a$.\end{enumerate}
\end{lem}

\begin{pf}
By symmetry it is enough to show (1). Let $E\in {\EE}_\xi^{n,m}$. Then
$E$ is given by an extension
\begin{equation*}\tag{$\epsilon$}0\maps \I_{Z_1}(F_1)\maps E\maps
\I_{Z_2}(F_2)\maps 0\end{equation*}
 with $\xi=2[F_1]-c_1$ and $length(Z_1)=n$, $length(Z_2)=m$.
Now if $b>a$ we have
$\Delta_b(\epsilon(lH))=\Delta_a(\epsilon(lH))+2(b-a)\<C\cdot\xi\>=2(b-a)\<C
\cdot\xi\>>0$
since $\Delta_a(\epsilon)=0$ and $\<C\cdot\xi\>>0$. Thus $E$ is $b$-unstable.
Assume that $E$ is not $a_-$-stable. Then it lies in an extension

\begin{equation*}\tag{$\epsilon_1$}0\maps \II_{Y_1}(G_1)\maps E\maps
\II_{Y_2}(G_2)\maps 0,\end{equation*}
 for which
 $\Delta_{a_-}(\epsilon_1(lH))\ge 0>\Delta_{a}(\epsilon_1(lH))$ for $l>>0$.

Hence we obtain $\<(2G_1-c_1)\cdot H\>\ge \<(2F_1-c_1)\cdot H\>$ and
$P_{a_-}(\I_{Z_1}(F_1+lH))<P_{a_-}(\I_{Y_1}(G_1+lH))$ and thus
$\chi(\I_{Z_1}(F_1+lH-C))<\chi(\I_{Y_1}(G_1+lH-C))$ or
$\chi(\I_{Z_1}(F_1+lH+C))<\chi(\I_{Y_1}(G_1+lH+C))$.
Consequently
$\Hom(\I_{Y_1}(G_1),\I_{Z_1}(F_1))=0$ and the obvious map
$\I_{Y_1}(G_1)\maps \II_{Z_2}(F_2)$ is an injection. Hence
$F_2-G_1$ is effective.
If $F_2\ne G_1$, we would have $\<(G_1-F_2)\cdot H\><0$, and, by
$\<\xi \cdot H\>=0$, we would get the contradiction $\<(2G_1-c_1)\cdot H\><0$.
So $G_1=F_2$. By the injectivity of $\II_{Y_1}(G_1)\maps \II_{Z_2}(F_2)$
and the fact that \ref{splitting} is not split, we get
$length(Z_2)<length(Y_1)$ which shoes that  $E$ is $a_-$-stable.
A similar argument shows that $E$ is strictly $a$-semistable.
\end{pf}

\begin{rem}\label{stab1} We can also easily see from the above arguments that
in the situation of \ref{stab} any sheaf $E\in M_{a_-}(c_1,c_2)$,
which does not lie in any $\EE_{\xi}^{n,m}$ for $(\xi,n,m)\in A^+(a)$
is $a_-$-stable (resp. semistable) if and only if it is $a$-stable  (resp.
semistable).
\end{rem}

\begin{rem}\label{fine} \begin{enumerate}
\item Looking at the proof of \cite{Ma2}
for the sufficient
criterion for the existence of a universal family on $M_H(c_1,c_2)$,
we see that the same proof also works for $M_{a}(c_1,c_2)$ and we
get the same criterion, i.e. if $c_1$ is not divisible by $2$ in
$NS(S)$ or otherwise $4c_2-c_2^2$ is not divisible by $8$ and
$M_a(c_1,c_2)=
M_a^s(c_1,c_2)$, then $M^s_a(c_1,c_2)$ carries a universal family.
\item From the results obtained so far it follows easily that,
under the above conditions
for the Chern classes, $M_a(c_1,c_2)=M^s_a(c_1,c_2)$ if and only if
$a$ is not a miniwall.
\end{enumerate}

\end{rem}

\begin{prop}\label{flip}
\begin{enumerate}
\item $M_0(c_1,c_2)=M_{H_-}(c_1,c_2)$ and $M_1(c_1,c_2)=M_{H_+}(c_1,c_2)$.
\item If $b\in [0,1]$ is not on a miniwall, the moduli space
$M_b(c_1,c_2)$ depends only on the minichamber in which $b$ is
lying, and
$M_b(c_1,c_2)\setminus M^s_b(c_1,c_2)$ is independent of $b$.
\item Let $a_-<a_+$ be in neighbouring minichambers separated
by the  miniwall
$a$. Then we have a set-theoretical decomposition
$$M_{a_+}(c_1,c_2)=\left(M_{a_-}(c_1,c_2)\setminus
\coprod_{(\xi,n,m)\in A^+(a)}
\EE^{n,m}_{\xi} \right){\sqcup}
\left(\coprod_{(\xi,n,m)\in A^+(a)}
\EE^{m,n}_{-\xi}\right),$$
 and there are morphisms
$$\matrix M_{a_-}(c_1,c_2)&&&&M_{a_+}(c_1,c_2)\cr
&\mapse{\psi_-}&&\mapsw{\psi_+}\cr
&&M_{a}(c_1,c_2)\cr\endmatrix$$
which are open embeddings over
$$M_{a_-}(c_1,c_2)\setminus
\coprod_{(\xi,n,m)\in A^+(a)}
\EE^{n,m}_{\xi} \text{ and }\ \ M_{a_+}(c_1,c_2)\setminus
\coprod_{(\xi,n,m)\in A^+(a)}
\EE^{m,n}_{-\xi}.$$
\end{enumerate}\end{prop}

\begin{pf} (1), (2), (3)  follow by putting together the results of this
section. By lemma \ref{stab} all the points of
$M_{a_-}(c_1,c_2)$ and $M_{a_+}(c_1,c_2)$ are
$a$-semistable and hence we get the morphisms
$\psi_-$ and $\psi_+$. The statement that they be open embeddings over the
indicated open subsets, follows  from remark
\ref{stab1}.\end{pf}

\section{The normal bundles of the exceptional sets}

Our aim in this and the next  chapter is to describe the
passage through a miniwall
which corresponds to a good wall.
We keep the assumptions from the beginning of the previous section.
In addition to those we assume that  either $p_g(S)=0$ or
$K_S$  is trivial, and that the wall
$W$ is good.

Let $a$ define a miniwall and let $(\xi,n,m)\in A^+(a)$.
Let $a_-<a_+$ lie in neighbouring minichambers separated
 by $a$.
For simplicity of notation we shall assume that $A^+(a)=\{(\xi,n,m)\}$.
Because,
for $(\xi,n_1,m_1)$, $ (\xi_2,n_2,m_2)$ distinct elements of $A^+(a)$,
the sets $\EE_{\xi_1}^{n_1,m_1}$ and $\EE_{\xi_2}^{n_2,m_2}$ are disjoint
by proposition \ref{wall} and  our arguments are local in a
neighbourhood of
each $\EE_{\eta}^{l,s}$, this assumption can be made without
loss of generality.
Furthermore we assume for simplicity of notation that $NS(S)$ has no
$2$-torsion. Then the classes $(c_1+\xi)/2$,
$(c_1-\xi)/2\in NS(S)$ are well-defined and
$T_\xi^{n,m}=M(1,(c_1+\xi)/2,n)\times M(1,(c_1-\xi)/2,m)$.
 Again this assumption is not important, as otherwise
the components of  $E_{\xi}^{n,m}$ and $E_{-\xi}^{m,n}$ are disjoint.

\begin{nota} \label{notsec4} We shall write
$M_-:=M_{a_-}(c_1,c_2)$, $M_+:=M_{a_+}(c_1,c_2)$,
$M^s_-:=M^s_{a_-}(c_1,c_2)$, $M^s_+:=M^s_{a_+}(c_1,c_2)$ and
put $\EE_-:=\EE_\xi^{n,m}$ and
$\EE_+:=\EE^{n,m}_{-\xi}$.
\end{nota}

\begin{defn}
Let $ \F_1'$ (resp. $ \F_2'$) be the pull-back of a universal sheaf over
$S\times  M(1,(c_1+\xi)/2,n)$ (resp. $S\times  M(1,(c_1-\xi)/2,m)$) to
$S\times T$, where
$T:= M(1,(c_1+\xi)/2,n)\times M(1,(c_1-\xi)/2,m)$. Let
$q=q_T:S\times T\to T$ be the projection.
Let $\A_{-}':=\Ext_{q}^1(\F_2',\F_1')$ and
$\A_{+}':=\Ext_{q}^1(\F_1',\F_2')$
and $\P_{-}:=\P(\A_{-}'),$ $\P_{+}:=\P(\A_{+}').$ Let
$\pi_-$ (resp. $\pi_+$) be the projections of $\P_-$
(resp. $\P_+$) to $T$ and $\tau_-$ (resp. $\tau_+$) the
tautological sublinebundles of $\A_-:=\pi_-^*(\A_-')$
(resp. $\A_+:=\pi_+^*(\A_+')$).
Let $\F_1:=(\id_S\times\pi_-)^*\F_1'$ and
$\F_2:=(\id_S\times \pi_-)^*\F_2'$.

\end{defn}

\begin{lem}\label{enm}
\begin{enumerate}
\item $ \A_-'$ is locally free of rank $-\xi(\xi-K_S)/2 +n+m-\chi(\oo_S)$
and its formation commutes with arbitrary base change.
\item There is an isomorphism $\Ext^1(\F_2',\F_1')\maps H^0(T,\A_-')$,
hence over $S\times \P_-$ there is a tautological extension
\begin{eqnarray}\label{globext}&&
0\maps \F_1\maps {\E}\to \F_2(\tau_-)\maps 0.\end{eqnarray}
There is a morphism $i_-:\P_-\maps M_-$ with image $\EE_-$.
\end{enumerate}
\end{lem}
\begin{pf}As $\xi$ defines a wall, $\Hom_{q}(\F_2,\F_1)$ is fibrewise $0$,
and, as the wall is good, $F_1-F_2+K_S$ is not effective for
$(F_1,F_2)\in T$,
therefore
by Serre duality for the extension groups \cite{Mu2} also
$\Ext^2_{q}(\F_2,\F_1)=0$.
So (1) follows by Riemann-Roch for the extension groups \cite{Mu2}.
Now we apply   \cite{La}.
\end{pf}

\begin{prop}\label{normal}
\begin{enumerate}
\item If $p_g(S)=0$ or if $K_S$ is trivial, then
$i_-:\P_-\maps M_-$ is a closed embedding  and
 $M_-$ is smooth along $\EE_-$. The irreducible component of
$M_-$ containing
$\EE_-$ has the expected dimension.
\item If $p_g(S)=0$, then the normal bundle $N_{\EE_-/M_-}$ of $\EE_-$ in
$M_-$ is equal to  $\A_+(\tau_-)$.
\item If $K_S$ is trivial, then $N_{\EE_-/M_-}=Q^\vee(\tau_-)$,
where $Q$ is the universal quotient bundle on $\P_-=\P(\A_-)$.
\end{enumerate}
\end{prop}

\begin{pf} By proposition \ref{wall}(3) and lemma \ref{enm} the map
$\P_-\maps M_-$ is injective with image $\EE_-$.  We also
see by proposition \ref{wall} that   $\EE_-\subset M_-^s$.
In case $K_S$ is trivial, $Spl(c_1,c_2)$  and thus also
the open subscheme
$M_-^s$ are smooth by \cite{Mu1}.

In order to see that $M_-$ is smooth along $\EE_-$ in the case $p_g(S)=0$,
we have to show that
$\Ext^2(E,E)=0$ for any $E\in {\EE}_-$.
So let $E\in {\EE}_-$ be given by a nontrivial extension (\ref{splitting})
\begin{equation*} \tag{$\epsilon$}
0\maps \I_{Z_1}(F_1)\maps E\maps \I_{Z_2}(F_2)\maps 0.\end{equation*}
As the wall $W$ is good,  we obtain
by Serre duality and the fact that $p_g(S)=0$ that
$\Ext^2(\I_{Z_i}(F_i),\I_{Z_j}(F_j))=0$ for $i=1,2$ and $j=1,2$.
Hence applying $\Ext^2(\I_{Z_i}(F_i),\cdot)$ to $(\epsilon)$
we get $\Ext^2(\I_{Z_i}(F_i),E)=0$ for $i=1,2$
and this in turn shows that  $\Ext^2(E,E)=0$.

We now want to compute the normal bundle to $\EE_-$.

{\it First Case: }$p_g(S)=0$.
Applying $\Hom_{q}(\cdot,\cdot)$ on both sides of
the sequence (\ref{globext})
and denoting by $\pi_i$  the composition of $\pi_-$
with the projection to the $i^{th}$ factor
we get the following exact diagram of locally free sheaves on $\P_-$
\begin{eqnarray}\label{globdiag}&&\matrix
&&0&&0&&0\cr
&&\mapd{}&&\mapd{}&&\mapd{}\cr
0&\mapr{}&T_{\P_-/T}&\mapr{}&\Ext^1_{q}(\F_2(\tau_-),\E)&\mapr{}&
\pi_2^*T_{M(1,(c_1-\xi)/2,m)}&\maps& 0\cr
&&\mapd{}&&\mapd{\psi}&&\mapd{}\cr
0&\mapr{}&\Ext^1_{q}(\E,\F_1)&\mapr{\phi}&i_-^*(T_{M_-})&\mapr{}&
\Ext^1_{q}(\E,\F_2)(\tau_-)&\maps& 0\cr
&&\mapd{}&&\mapd{}&&\mapd{}\cr
0&\mapr{}&\pi_1^*T_{M(1,(c_1+\xi)/2,n)}&\mapr{}
&\Ext^1_{q}(\F_1,\E)&\mapr{}&
\A_+(\tau_-)&\maps& 0\cr
&&\mapd{}&&\mapd{}&&\mapd{}\cr
&&0&&0&&0\cr
\endmatrix\end{eqnarray}
To identify the entries in this diagram we have used the following facts.
\begin{enumerate}
\item $\Hom_q(\F_1,\F_1)=\Hom_q(\F_2,\F_2)=\oo_{\P_-}$.
\item If $Q$ is the universal quotient on
$\P(\A_-)$, then the relative tangent bundle is
$T_{\P_-/T}=Q(-\tau_-)$, i.e. the cokernel of the natural map
$\oo_{\P_-}=\Hom_q(\F_1,\F_1)\maps \Ext^1_q(\F_2(\tau_-),\F_1)$.
\item $\pi_2^\ast(T_{M(1,(c_1-\xi)/2,m)})=\Ext^1_{q}(\F_2,\F_2)$ and
$\pi_2^\ast T_{M(1,(c_1+\xi)/2,n)}=\Ext^1_{q}(\F_1,\F_1)$.
\item By Mukai's sheafified Kodaira-Spencer map \cite{Mu1} we have
$i_-^*T_{M_-}=\Ext^1_{q}(\E,\E).$ Mukai shows the result only if
$S$ is an abelian or K3-surface, but in his  proof he only uses that
$Spl(c_1,c_2)$ is smooth in a neighbourhood of $\EE_-$,
(which we have just seen)
and  $\Ext^1_{q}(\E,\E)$ is locally free and
compatible with base change.
\end{enumerate}
To show that the sequences in the diagram are exact
we just use standard techniques.
It is  enough to check the exactness fibrewise.
One has repeatedly to make use of the fact that $\xi$
defines a good wall,
i.e. if $E\in {\EE}_-$ is given by (\ref{splitting}),
then $F_1-F_2$, $F_2-F_1$, $F_1-F_2+K_S$, $F_2-F_1+K_S$ are not effective,
which implies that
$\Hom_q(\F_1,\F_2(\tau_-))= \Hom_q(\F_2(\tau_-),\F_1)=
\Ext^2_q(\F_1,\F_2(\tau_-))= \Ext^2_q(\F_2(\tau_-),\F_1)=0$.
In addition we use that all $E\in {\EE}_-$ are simple and that
$\Ext^2_{q}(\F_2,\F_2)=\Ext^2_{q}(\F_1,\F_1)=0$.
We also use the vanishings from the proof of the smoothness of
$M_-$ along $\EE_-$.

\noindent{\it Second Case: } $K_S$ is trivial.

We apply essentially the same arguments as in the first case.
Now however we have
$\Ext^2_{q}(\E,\F_1)=\Ext^2_{q}(\F_2(\tau_-),\E)=\Ext^2_{q}(\E,\E)=
\Ext^2_{q}(\F_1,\F_1)=\Ext^2_{q}(\F_2,\F_2)=\oo_{\P_-}$,
which  follows easily from Mukai's results \cite{Mu1}.
We also notice that by Serre-duality $\A_+$ is canonically dual to
$\A_-$.
Using all this we again get the diagram (\ref{globdiag}) with the entry
$\A_+(\tau_-)$ in the lower right corner replaced by the
kernel of the natural map
$\A_+(\tau_-)\to \oo_{\P_-}$, i.e. $Q^\vee(\tau_-)$.

\noindent{\it Claim: } The image of the Kodaira-Spencer map
$\kappa:T_{\P_-}\to \Ext^1_{q}(\E,\E)$
 is $Im(\phi)$+$Im(\psi)$ (see (\ref{globdiag})).

Note that, by what we have shown so far, the claim implies the theorem.

\noindent{\it Proof of the Claim.}
For dimension reasons it is enough to show that $Im(\phi)$ and $Im(\psi)$
both are contained in the image of $\kappa$.
We show it for $Im(\phi)$. It is enough to show this fibrewise.
Let $F_1\in M(1,(c_1+\xi)/2,n)$ and let $(\P_-)_{F_1}$  be the
fibre of the projection
$\pi_1:\P_- \maps  M(1,(c_1+\xi)/2,n)$ over $F_1$.
Then
$(\P_-)_{F_1}$ is
the space of extensions
$$0\to F_1\to E\to G\to 0$$ with
$G$ running through $M(1,(c_1-\xi)/2,m)$.
 Let $x\in ({\P}_-)_{F_1}$ be given by an
extension
\begin{equation}\tag{$\lambda_x$} 0\to F_1\to E\to G_1\to 0.\end{equation}
We will want to show that
$\kappa(T_{(\P_-)_{F_1}}(x))=\phi(\Ext^1(G_1,E))$.

The tangent space to $(\P_-)_{F_1}$ at $x$ is the
space of first order deformations of $E$ together with an injection
$F_1\to E$. For $t\in T_{(\P_-)_{F_1}}(x)$ we get therefore the diagram
$$\matrix
&&0&&0&&0\cr
&&\mapd{}&&\mapd{}&&\mapd{}\cr
0&\mapr{}&F_1&\mapr{}&E&\mapr{}&G_1&\mapr{}&0\cr
&&\mapd{}&&\mapd{}&&\mapd{}\cr
0&\mapr{}&F_1\oplus F_1&\mapr{}&\widetilde
 E&\mapr{}&\widetilde G&\mapr{}&0\cr
&&\mapd{}&\mapne{\gamma}&\mapd{}&&\mapd{}\cr
0&\mapr{}&F_1&\mapr{}&E&\mapr{}&G_1&\mapr{}&0\cr
&&\mapd{}&&\mapd{}&&\mapd{}\cr
&&0&&0&&0\cr
\endmatrix
\eqno (*)$$
and we see that $T_{(\P_-)_{F_1}}(x)$ can be identified
with the space of  diagrams $(*)$.
Furthermore $\kappa(t)$ is the extension class of the middle column
of $(*)$.
{}From $(*)$ we also get a sequence
$0\maps  E\maps  \widetilde E/\gamma(F_1)\maps G_1\to 0$
such that $\widetilde E$ is defined by pull-back
$$\matrix
\widetilde E/\gamma(F_1)&\mapr{}&G_1\cr
\mapu{}&&\mapu{}\cr
\widetilde E.&\mapr{}&E\cr
\endmatrix$$
This gives a map $\theta:T_{(\P_-)_{F_1}}(x)\maps Ext^1(G_1,E)$, such that the
restriction of
$\kappa$ to $T_{(\P_-)_{F_1}}(x)$ is $\phi\circ\theta$.
To finish the proof we  have to see that
$\theta$ is an isomorphism. We give an inverse.
Let
$$0\maps E\maps W\maps G_1\maps 0$$ be an extension.
We define $\widetilde E$ as the fibre product
$$\matrix
W&\mapr{}&G_1\cr
\mapu{}&&\mapu{}\cr
\widetilde E,&\mapr{}&E\cr
\endmatrix$$
and we see that it lies in a diagram $(*)$.
\end{pf}

\begin{rem}\label{newcomp} Assume $p_g(S)=0$.  From
lemma \ref{normal} it follows that the dimension of $\EE_-$ is at
most the expected dimension
$N=(4c_2-c_1^2)-3\chi(\oo_S)+q(S)$. We have to distinguish two cases.
\begin{enumerate}
\item  $dim(\EE_-)<N$ and   $dim(\EE_+)<N$. Then   the change
from $M_-$ to $M_+$ is a birational transformation.
\item  $dim(\EE_-)=N$ or $dim(\EE_+)=N$. We can assume that $dim(\EE_-)=N$.
 Then by  lemma \ref{normal}
$\EE_-$ is a smooth connected  component of $M_-$, which is isomorphic to
$\P_-$. And,  $\A_+(\tau_-)$ being the normal bundle to $\EE_-$,
we have $\A_+=0$ and therefore $\EE_+=\emptyset$.
This happens if and only if $\<\xi\cdot(\xi-K_S)\>/2+d_\xi=\chi(\oo_S)$.
If we allow $NS(S)$ to contain $2$-torsion, we see that all
the connected components of $\EE_-$ are connected components of $M_-$.
\end{enumerate}
\end{rem}

Assume for the following definition and corollary that we are in case (1) of
\ref{newcomp}, i.e. that the change from $M_-$ to $M_+$ is
birational.

\begin{defn} Let $\widetilde M_-$ be the blow-up of
$M_-$ along $\EE_-$ and $D$ the exceptional divisor.
Similarly let $\widetilde M_+$ be the blow up of $M_+$ along $\EE_+$.
Let $\pi_D$, $\pi_{D-}$, $\pi_{D+}$ be the projections
from $D$ to $T$, $\P_-$, $\P_+$ respectively.
\end{defn}

\begin{cor}\begin{enumerate}
\item
If $p_g(S)=0$ then $D$ is isomorphic to
$\P_-\times_T \P_+$ and with this identification
$\oo(D)|_D=\oo(\tau_-+\tau_+)$.
\item If $K_S$ is trivial, then $\A_-$ and $\A_+$ are canonically dual and
 $D$ is the incidence correspondence
$\{(l,H)\in \P(\A_-)\times_T\P^{\vee}(\A_-) \,|\, l\subset H\}$ and
$\oo(D)|_D$ is the restriction of $\oo(\tau_-+\tau_+)$.
\end{enumerate}
\end{cor}

\section{Blow-up construction}

We keep the assumptions and notations of the last section.
In addition we assume that we  in case (1) of
\ref{newcomp}, i.e. the map $\widetilde M_-\maps M_-$ is birational.
In this section we want to show that $\widetilde M_-$ and $\widetilde M_+$
are isomorphic.
We shall construct a morphism $\phi_+:\widetilde M_-
\to M_+$,  which we shall show is the blow-up of $M_+$
 along $\EE_+$.
Let $\phi_-:\widetilde M_-\to M_-$ be the blow-up  map and
$j:D\to \widetilde M_-$ be the embedding.
We denote $\widetilde M_-^s:=\phi_-^{-1}M_-^s$.
Let $\U_-$ be a universal pseudo-family on $S\times M^s_-$ and
$\V_-:=(\id_S\times \phi_-)^*\U_-$.
We want to make an elementary transform of $\V_-$ along $D_S:=S\times D$
to obtain a pseudo-family $\V_+$ of $a_+$-stable sheaves on
$\widetilde M_-^s$ and
thus the
desired map $\phi_+$. If $\U_-$ is a universal family, then also
$\V_+$ will be one.

\begin{nota}
For a sheaf $\H$ on $S\times\P_-$  (resp. $S\times\P_+$) we will write
$\H_D$ for $(\id_S\times \pi_{D-})^*\H$
(resp. $(\id_S\times \pi_{D+})^*\H$). We also write
$\F_{1D}$ and $\F_{2D}$ instead of $(\F_{1})_{D}$ and $(\F_{2})_{D}$.
\end{nota}

\begin{defn}
By the universal property of $M_-$ and lemma \ref{enm} there is
a line bundle
$\lambda$ on $D$ such that there is  an exact sequence
\begin{eqnarray}
\label{restrsec}&&0\to \F_{1D}(\lambda)\maps \V_-|_{D_S}\maps
\F_{2D}(\tau_-+\lambda)\maps 0,
\end{eqnarray}
indeed there is already a sequence like this on $\EE_-$.
Let $\gamma$ be the composition $\V_-\maps \V_-|_{D_S}\maps
\F_{2D}(\tau_-+\lambda).$ Then we put $\V_+:=ker\gamma$.
Because $\V_-$ is flat on $S\times \widetilde M^s_-$, and
 $\F_{2D}(\lambda+\tau_-)$
is flat on the Cartier divisor $S\times D$, $\V_+$ is flat over
 $S\times \widetilde M^s_-$.

The restrictions of $\V_+$ and $\V_-$ to
$S\times \widetilde M^s_-\setminus D$ are naturally isomorphic.
There are diagrams
of sheaves on $S\times \widetilde M_-^s$
\begin{eqnarray}\label{dia1}
&&{\matrix
&&&&0&&0\cr
&&&&\mapd{}&&\mapd{}\cr
0&\mapr{}&\V_-(-D_S)&\mapr{}&\V_+&\mapr{}&\F_{1D}(\lambda)&\mapr{}&0\cr
&&\Big|\Big|&&\mapd{}&&\mapd{}\cr
0&\mapr{}&\V_-(-D_S)&\mapr{}&\V_-&\mapr{}&\V_-|_{D_S}&\mapr{}&0\cr
&&&&\mapd{}&&\mapd{}\cr
&&&&\F_{2D}(\tau_-+\lambda)&\longleich&\F_{2D}(\tau_-+\lambda)
&\mapr{}&0\cr
&&&&\mapd{}&&\mapd{}\cr
&&&&0&&0\cr\endmatrix}\end{eqnarray}
\begin{eqnarray}\label{dia2}&&{\matrix
&&&&0&&0\cr
&&&&\mapd{}&&\mapd{}\cr
0&\mapr{}&\V_+(-D_S)&\mapr{}&\V_-(-D_S)&\mapr{}
&\F_{2D}(\lambda-\tau_+)&\mapr{}&0\cr
&&\Big|\Big|&&\mapd{}&&\mapd{}\cr
0&\mapr{}&\V_+(-D_S)&\mapr{}&\V_+&\mapr{}&\V_+|_{D_S}&\mapr{}&0\cr
&&&&\mapd{}&&\mapd{}\cr
&&&&\F_{1D}(\lambda)&\longleich&\F_{1D}(\lambda)&\mapr{}&0\cr
&&&&\mapd{}&&\mapd{}\cr
&&&&0&&0\cr\endmatrix}
\end{eqnarray}
By the rightmost column of (\ref{dia2}), $(\V_+)_x\in \EE_+$
for all $x\in D$.
Therefore by proposition  \ref{flip}
$\V_+$ is a  pseudo-family  of $a_+$-stable sheaves
 over $\widetilde M^s_-$ and defines
a morphism $\phi_+:\widetilde M^s_- \to M^s_+$.
We see from the definitions that the restriction of $\phi_+$ to
$\widetilde M^s_-\setminus D$ is an isomorphism to $M^s_+\setminus \EE_+$,
which coincides with the natural identification
$\widetilde M^s_-\setminus D\simeq M^s_-\setminus
\EE_-\simeq M^s_+\setminus \EE_+$.
As $\EE_-\subset M_-^s$ and $\EE_+\subset M_+^s$,
 we see that $\phi_+$ extends to a morphism
$\widetilde M_-\maps M_+$, which we still denote by
$\phi_+$.
\end{defn}

\begin{thm}
$\phi_+:\widetilde M_-\maps M_+$ is the blow up of $M_+$ along
$\EE_+$.
\end{thm}

\begin{pf} By the above $\phi_+(D)\subset \EE_+$.
We want to show that $\phi_+|_D$ is the projection $\pi_{D+}:D\maps  \EE_+$.
For this we have to show that the extension
$$0\maps\F_{2D}(\lambda-\tau_+)\maps \V_+|_{S\times D}
\maps  \F_{1D}(\lambda)\mapr{} 0$$ from the rightmost column of (\ref{dia2})
is the pull-back via $\pi_{D+}$ of the tautological extension on $\P_+$
(defined analogously to \ref{globext})
tensorized with $\oo_D(\lambda-\tau_+)$.
It is enough to show this fibrewise.

Let $x=(x_-,x_+)\in D\subset \P_-\times_T \P_+$ and
let $V_-:=(\V_-)_x$ and $V_+:=(\V_+)_x$ be given by extensions
\begin{eqnarray}\label{seq1}
&&0\maps F_1\maps V_-\maps F_2\maps 0,\\
\label{seq2}
&&0\maps F_2\maps V_+\maps F_1\maps 0.
\end{eqnarray}
Then $\pi_D(x)$ is the point $(F_1,F_2)\in T$ and $x_-\in (\P_-)_{(F_1,F_2)}=
\P(\Ext^1(F_2,F_1))$
 is the extension class of (\ref{seq1}).
Then we have to show that $x_+\in (\P_+)_{(F_1,F_2)}=\P(\Ext^1(F_1,F_2))$
 is the extension class of (\ref{seq2}).

Let $R:=\spec\C[\epsilon]/(\epsilon^2)$ and
let $t:R\to \widetilde M_-$ be a tangent vector to $\widetilde M_-$ at  $x$,
which is not tangent to $D$. Then $t$ factors through $2D$ (i.e. the subscheme
defined by $\I_D^2$). If we restrict the diagrams (\ref{dia1}), (\ref{dia2})
to $2D_S$, we see that the image of the map
$\V_-(-D_S)|_{2D_S}\maps \V_+|_{2D_S}$ is $\I_{D_S}\V_-/\I^2_{D_S}\V_-$
and the image of  the composition
$\V_+(-D_S)|_{2D_S}\maps\V_-(-D_S)|_{2D_S}\maps V_+|_{2D_S}$ is
$\I_{D_S}\F_{1D}(\lambda)\cdot/\I_{D_S}^2\F_{1D}(\lambda)$. Therefore,
by pulling back the diagrams
(\ref{dia1}), (\ref{dia2}) to $S\times R$ via $(\id_S\times t)$
and pushing down with the projection
$p:S\times R\to S$, we get the diagrams
\begin{eqnarray}\label{dia3}\matrix
&&&&0&&0\cr
&&&&\mapd{}&&\mapd{}\cr
0&\mapr{}&V_-&\mapr{}&\widetilde V_+&\mapr{}&F_1&\mapr{}&0\cr
&&\Big|\Big|&&\mapd{}&&\mapd{}\cr
0&\mapr{}&V_-&\mapr{}&\overline V_-&\mapr{}&V_-&\mapr{}&0\cr
&&&&\mapd{}&&\mapd{}\cr
&&&&F_2&\longleich&F_2&\mapr{}&0\cr
&&&&\mapd{}&&\mapd{}\cr
&&&&0&&0\cr\endmatrix\end{eqnarray}
\begin{eqnarray}
\label{dia4}
\matrix
&&&&0&&0\cr
&&&&\mapd{}&&\mapd{}\cr
0&\mapr{}&F_1&\mapr{}& V_-&\mapr{}&F_2&\mapr{}&0\cr
&&\Big|\Big|&&\mapd{}&&\mapd{}\cr
0&\mapr{}&F_1&\mapr{}&\widetilde V_+&\mapr{}&V_+&\mapr{}&0\cr
&&&&\mapd{}&&\mapd{}\cr
&&&&F_1&\longleich&F_1&\mapr{}&0\cr
&&&&\mapd{}&&\mapd{}\cr
&&&&0&&0\cr\endmatrix
\end{eqnarray}
The extension class $\delta\in \P(\Ext^1(V_-,V_-))$
of the middle row of (\ref{dia3}) is the class of the image of $t$ under
$d\phi_-:T_{\widetilde M_-}(x)\maps T_{ M_-}(\phi_-(x))=\Ext^1(V_-,V_-)$.
The image of the composition
$$T_{ M_-}(x)\maps
 \phi_-^*(T_{M_-}(x))\,\mapr{\rho} \,\phi_-^*(N_{\EE_-/M_-}(x))$$ is
the tautological subline-bundle of
$\phi_-^*(N_{\EE_-/M_-}(x))=\Ext^1(F_1,F_2)$ and the kernel is $T_{D}(x)$.
Therefore the image of $\rho(\delta)$ in
$(\P_+)_{(F_1,F_2)}=\P(\Ext^1(F_1,F_2))$ is  $x_+$.
By (\ref{globdiag}) the map $\rho$ is the composition
$$ \Ext^1(V_-,V_-)\,\mapr{\rho_1}\,\Ext^1(F_1,V_-)\,\mapr{\rho_2}\,
\Ext^1(F_1,F_2)$$ given by applying $\Hom(\cdot,\cdot)$ on both
sides to the
sequence
$0\to F_1\to V_-\to F_2\to 0$. By
(\ref{dia3}) $\rho_1(\delta)$ is the extension class of the
first row of (\ref{dia3}) giving  $\widetilde V_+$, and then,  by
(\ref{dia4}), $\rho(\delta)$ is the extension class of (\ref{seq2}).
So we see that
$\phi_+|_D$ is the projection to $\EE_+$.

If for the moment we call $\overline \phi_+:\widetilde M_+\to M_+$
the blow-up of $M_+$ along $\EE_+$ and $\overline D$ the
exceptional divisor,
 we get analogously that $\overline D\simeq \P_-\times_T \P_+$,
(or the incidence correspondence in $\P_-\times_T \P_+$
in case $K_S$ is trivial).
In the same way as above we can construct a morphism
$\overline \phi_-:\widetilde M_+\maps M_-$ such that
$\overline \phi_-|_{\overline D}$
is the projection to $\EE_-$ and
$\overline \phi_-|_{\widetilde M_+\setminus \overline D}$
is just the natural identification
$\widetilde M_+\setminus \overline D\simeq M_+\setminus
\EE_+\simeq M_-\setminus \EE_-$.
Therefore we have morphisms
$\phi_-\times  \phi_+:\widetilde M_-\to M_-\times M_+$,
$\overline \phi_-\times  \overline \phi_+:\widetilde M_+\to M_-\times M_+$,
which by the above are injective and easily seen to be injective on tangent
vectors. Furthermore
$(\phi_-\times  \phi_+)(\widetilde M_-\setminus D)=
(\overline \phi_-\times  \overline \phi_+)
(\widetilde M_+\setminus \overline D)$.
Therefore $\widetilde M_-$ and $\widetilde M_+$ are isomorphic
and in fact both  isomorphic to the closure of the graph of the obvious
rational map $M_-\to M_+$.
\end{pf}

In the following theorem we put together the main results we
have obtained so far.

\begin{thm}\label{zblowup}
Let $S$ be a surface with either $p_g(S)=0$ or $K_S$ trivial.
Let $c_1\in NS(S)$, $c_2\in \Z$ and put
$N:=4c_2-c_1^2-3\chi(\oo_S)+q(S)$.
Let $W$ be a good wall of type $(c_1,c_2)$ and let $H_-$,  $H_+$ be
ample divisors on $S$ in neighbouring chambers separated by
$W$. Then for all $a\in [0,1]$ there exist spaces $M_a(c_1,c_2)$
and a finite set  of miniwalls dividing  $[0,1]$ into
finitely many minichambers
such that the following holds:
\begin{enumerate}
\item $M_0(c_1,c_2)=M_{H_-}(c_1,c_2)$, $M_1(c_1,c_2)=M_{H_+}(c_1,c_2)$.
\item If $a_1$, $a_2$ are in the same minichamber
 then $M_{a_1}(c_1,c_2)=M_{a_2}(c_1,c_2)$.
\item If $a_-<a<a_+$ and $a$ is the unique miniwall
between $a_-$ and $a_+$ then
$M_{a_+}(c_1,c_2)$ is obtained from $M_{a_-}(c_1,c_2)$
as follows:
We blow up $M_{a_-}(c_1,c_2)$ along  the disjoint smooth subvarieties
$\EE_{\xi}^{n,m}$, with  $(\xi,n,m)\in A^+(a)$  (see \ref{defminiwall})
which fulfill $0\le dim(\EE_{\xi}^{n,m})< N$
and blow-down the exceptional divisors   to $\EE_{-\xi}^{m,n}$
respectively.
Then we remove the  $\EE_{\xi}^{n,m}$  with $(\xi,n,m)\in\A^+(a)$ and
$dim(\EE_{\xi}^{n,m})= N$
(which are unions of connected  components of $M_{a_-}(c_1,c_2)$)
and take the disjoint  union with all $\EE_{-\xi}^{m,n}$ with
$(\xi,n,m)\in\A^+(a)$ and
$\EE_{\xi}^{n,m}= \emptyset$.
\item If $H$ is an ample divisor on $W$ which lies in the
closure of both
of the chambers containing $H_-$ and $H_+$, then, for all
$b\in [0,1]$, the space $M_b(c_1,c_2)$
is a  moduli space of  H-semistable sheaves on $S$ with a
suitable parabolic structure.
\end{enumerate}
\end{thm}

In \cite{Mu1} Mukai defines  elementary transforms of a
symplectic variety $X$
as follows. Assume $X$ contains a subvariety $P$, which has codimension
$n$ and is a $\P_n$-bundle over a variety $Y$. Let $\widetilde X$
be the blow-up of $X$ along $P$. Then the exceptional divisor $E$
is isomorphic to the incidence correspondence in $P\times_Y P'$,
where $P'$ is  the dual projective bundle to $P$.
One can then blow down $E$ to $P'$ to obtain a smooth
symplectic variety $X'$.
We will for the moment call $Y$ the center of such an
elementary transformation.

So by the above we obtain the following:

\begin{cor}
Let $S$ be a K3-surface or an abelian surface.
Let $H_-,H_+$ be polarisations which both do not lie on a wall.
Then $M_{H_+}(c_1,c_2)$ is obtained from $M_{H_-}(c_1,c_2)$
by a series of elementary transforms, whose centers are of the form
$M(1,(c_1+\xi)/2,n)\times M(1,(c_1-\xi)/2,m)$ for $\xi$
defining a wall between
$H_-$ and $H_+$ and $(n,m)$ running through the nonnegative integers with
$n+m=(4c_2-c_1^2+\xi^2)/4$.
\end{cor}

\begin{rem} If $q(S)\ne 0$ we can also, for $A\in Pic(S)$, $c_2\in \Z$
and an ample divisor $H$,
study the moduli space $\widetilde M_H(A,c_2)$ of rank $2$
torsion-free sheaves
$E$ on $S$ with $det(E)=A$ and $c_2(E)=c_2$. Then there is a morphism
$M_H(c_1,c_2)\maps Pic^0(S)$, whose fibres are the various $\widetilde
M_H(A,c_2)$ for $A$ with $c_1(A)=c_1$.  Then, by restricting our arguments
to the fibres, we get that theorem \ref{zblowup}
also holds with the obvious changes for  $\widetilde M_H(A,c_2)$.
\end{rem}

\section{The change of the Donaldson invariants in terms of Hilbert schemes}

In this section we assume that $q(S)=0$.
Let $\gamma_{c_1,c_2,g}$ be the Donaldson polynomial
with respect to a Riemannian metric $g$
associated to the principal $SO(3)$-bundle $P$ on $S$ whose
second Stiefel-Whitney class
$w_2(P)$ is the reduction of $c_1$ mod $2$ and whose first
 Pontrjagin class is
$p_1(P)=(c_1^2-4c_2)$. Then $\gamma_{c_1,c_2,g}$
is a homogeneous polynomial on $H_*(S,\Q)$ of weight
$2N=2(4c_2-c_1^2-3\chi(\oo_S))$,
 where the elements of  $H_i(S,\Q)$  have weight
$4-i$. In case $p_g(S)>0$ it is known that $\gamma_{c_1,c_2,g}$
does not depend on the metric (as long as it is generic).

In [Ko] the invariant has been introduced for $4$-manifolds $M$ with
$b_+(M)=1$.
In [K-M] it has been shown that in case $b_+(M)=1$, $b_1(M)=0$
it  depends only on the chamber of the
period point of the metric in the positive cone of $H^2(M,\R)$.

The algebro-geometric analogues of the Donaldson polynomials are
defined as follows:

\begin{defn} \label{algana} (\cite{OG1}, \cite{OG2})
Assume that $M_H(c_1,c_2)$ is a fine moduli space, i.e.
$M_H(c_1,c_2)=M_H^s(c_1,c_2)$, and there is a universal sheaf $\U$ on
$S\times M_H(c_1,c_2)$. We define a linear map
$$\nu_{c_1,c_2,H}:H_i(S,\Q)\to H^{4-i}(M_H(c_1,c_2),\Q);
\ \ \nu_{c_1,c_2,H}(\alpha):=(c_2(\U)-{1\over 4}c_1^2(\U))/\alpha,$$
where  $/$ denotes the slant product.
We assume furthermore that $M_H(c_1,c_2)$ is of  the expected dimension
$N:=4c_2-c_1^2-3\chi(\oo_S)$.
Given classes $\alpha_{s}\in H_{2j_s}(S,\Q)$, for  $s=1,\ldots,k$
with $2k-\sum_s j_s=N$,
we set
$$\Phi_{c_1,c_2,H}(\alpha_1,\ldots,\alpha_k):=\intt_{M_H(c_1,c_2)}
\nu_{c_1,c_2,H}(\alpha_1)\cdot \ldots \cdot\nu_{c_1,c_2,H}(\alpha_k).$$
As $c_1,c_2$ are fixed in our paper, we will write
$\nu_H=\nu_{c_1,c_2,H}$ and $\Phi_{H}=\Phi_{c_1,c_2,H}$.

Let $\pt\in H_0(S,\Z)$ be the class of a point in $S$. Knowing
$\Phi_{H}$ is equivalent to knowing the numbers
$$\Phi_{H,l,r}(\alpha)
:=\intt_{M_H(c_1,c_2)}\nu_H(\alpha)^l\cdot\nu_H(\pt)^r.$$
for all $l,r$ with $l+2r=N$ and all $\alpha\in
H_2(S,\Q)$.
\end{defn}

\begin{defn}
Following \cite{OG2}, we call $M_H(c_1,c_2)$
admissible if the  following holds:
\begin{enumerate}
\item $H$ does not lie on a wall of type $(c_1,c_2)$;
\item $dim(M_H(c_1,c_2))=N$,
\item if $c_1$ is divisible by $2$ in $NS(S)$, then
$N> (4c_2-c_1^2)/2$;
\item $dim(M_H(c_1,k))+2(c_2-k)<N$ for all $k<c_2$.
\end{enumerate}
\end{defn}

For admissible  $M_H(c_1,c_2)$  the results of
  \cite{Mo} and \cite{Li} give
$$\Phi_{H}|_{H^2(S,\Q)}=(-1)^{(c_1^2+\<c_1\cdot K_S\>)}
\gamma_{c_1,c_2,g(H)}|_{H^2(S,\Q)},$$ where
$g(H)$ is the Fubini-Study metric associated to $H$.
Furthermore if $c_2>>0$, then
$\Phi_{H}=(-1)^{(c_1^2+\<c_1\cdot K_S \>)}\gamma_{c_1,c_2,g(H)}$.

We now want to determine  how $\Phi_{H}$ changes, when $H$
passes through a wall.
We assume that
if $c_1$ is  divisible by $2$ in $NS(S)$ then
$(4c_2-c_1^2)$ is not divisible by $8$. Then, by the criterion of
\cite{Ma2}, $M_H(c_1,c_2)$ is a fine moduli space, unless
$H$ lies on a wall.

Now we assume that  we are in the situation of section 3, i.e.
$H_-$ and $H_+$ are ample divisors lying in neighbouring chambers separated by
$W$, and $H$ a polarization on the wall $W$ not lying  on any other
wall and lying in the closure of both the chambers containing $H_-$ and $H_+$.
We assume furthermore that $W$ is a good wall.
For $b\in [0,1]$ we have $M_b(c_1,c_2)$ as in section 3.

\begin{defn} By remark \ref{fine} we see that, for $b$ not on a miniwall,
$M_b(c_1,c_2)=M^s_b(c_1,c_2)$ and there is a universal sheaf on
$M_b(c_1,c_2)$.

Assume  that $b\in [0,1]$ does not lie on a miniwall. Then analoguosly to
the definition of
$\Phi_{H}$ and $\Phi_{H,l,r}$ in \ref{algana},
we may define $\Phi_{b}$ and $\Phi_{b,l,r}$ by always replacing
 $M_H(c_1,c_2)$ by $M_b(c_1,c_2)$.
\end{defn}

We notice that $\Phi_{H_-}=\Phi_{0}$ and $\Phi_{H_+}=\Phi_{1}$
and it is obvious that $\Phi_{b}$ only depends on the
minichamber containing $b$.
We therefore  have to determine the change of $\Phi_b$ when $b$
passes through a miniwall.

We will make  the same assumptions as in section 4, i.e. let
$a$ be a miniwall and let $(\xi,n,m)\in A^+(a)$. Let $a_-<a_+$
lie in neighbouring minichambers separated by $a$.
To  simplify the  notation we will  for the moment assume that
$A^+(a)=\{ (\xi,n,m)\}$ and that $H^2(S,\Z)$ contains no $2$-torsion.
We also assume that either $p_g(S)=0$ or $K_S$ is trivial.

\begin{nota}\label{notado}
We use the notations and definitions of sections  4 and 5.
If the change is birational, i.e.  we are not in case (1) of \ref{newcomp},
we shall write $\widetilde M$ instead of $\widetilde M_-$.
Let $d:=d_\xi=n+m$,  $e_-=rk (\A_-)$, $e_+=rk (\A_+)$,
then $N=2d+e_-+e_+-1$ if $p_g(S)=0$ and $N=2d+e_-+e_+-2$
if $K_S$ is trivial.
We put $\nu_+:=\nu_{a_+}$, $\nu_-:=\nu_{a_-}$, $\Phi_+:=\Phi_{a_+}$,
 $\Phi_-:=\Phi_{a_-}$, $\Phi_{+,l,r}:=\Phi_{a_+,l,r}$ and
$\Phi_{-,l,r}:=\Phi_{a_-,l,r}$.

Note that  the condition $q(S)=0$ implies $Pic(S)\simeq NS(S)$.
For $\beta\in NS(S)$ we may therefore denote  by $\oo_S(\beta)$
the corresponding
line bundle. Let $q_1,q_2$ be the two projections of $T=\Hilb^n(S)\times
\Hilb^m(S)$.
\end{nota}

\begin{rem}\begin{enumerate}
\item If the change is birational, then
by the projection formula
$\Phi_+,$ $\Phi_{+,l,r}$ (resp. $\Phi_-$, $\Phi_{-,l,r}$)
coincide with the numbers which are defined analogously by  replacing
$M_{a}(c_1,c_2)$ by $\widetilde M$ and the universal sheaf by
$\V_+$ (resp. $\V_-$).
\item Assume $p_g(S)=0$ and say $\EE_+=\emptyset$.
Let $\E$ be the universal
sheaf on $\EE_-$ from (\ref{globext}), then we can define
$\sigma_-:H_i(S,\Q)\maps H^{4-i}(\EE_-,\Q)$ and  $\delta_-$
and  $\delta_{-,l,r}$
in the same way as   $\nu_-$ and $\Phi_-$ and   $\Phi_{-,l,r}$ by
 replacing $M_-$ by $\EE_-$ and
the universal sheaf on $M_-$ by $\E$.
Then $\Phi_+-\Phi_-=-\delta_-$.
\end{enumerate}
\end{rem}

\begin{defn}\label{hilbdef}
Let $Z_n(S)\subset S\times \Hilb^n(S)$ be the universal subscheme.
In  $S\times \Hilb^n(S)\times \Hilb^m(S)$,
we put $\ZZ_1:=(\id_S\times q_1)^{-1}(Z_n(S))$,
$\ZZ_2:=(\id_S\times q_2)^{-1}(Z_m(S))$ and denote by
$\I_{\ZZ_1}$, $\I_{\ZZ_2}$ the corresponding idealsheaves.
Let $F_1:=\oo_S((c_1+\xi)/2)$, $F_2:=\oo_S((c_1-\xi)/2)$.
By our assumptions
$T=\Hilb^n(S)\times \Hilb^m(S)$ and $\F_1'=\I_{\ZZ_1}(F_1)$,
$ \F_2'=\I_{\ZZ_2}(F_2)$.
Let $h_n:\Hilb^n(S)\maps S^{(n)}$ be the Hilbert-Chow morphism
\cite{Fo}, where $S^{(n)}$ is the $n$-fold symmetric power of $S$
with the quotient map $\phi_n:S^n\maps S^{(n)}$.
For $i=1,\ldots,n$ we denote by $p_i:S^n\to S$
the projection to the $i^{th}$ factor.
We denote
$\Delta_{i}:=\big\{(x,x_1,\ldots,x_n)\in S\times S^n\bigm| x=x_i\big\}$
and $Y_n:=(\id_S\times \phi_n)(\Delta_{1}).$
We have linear maps
\begin{eqnarray*}
\iota_n:H_i(S,\Q)\maps H^{4-i}(\Hilb^n(S),\Q);&& \iota_{n}(\alpha)=
[Z_n(S)]/\alpha\   \    \text{ and }\\
\bar\iota_n:H_i(S,\Q)\maps H^{4-i}S^{(n)},\Q);&& \bar\iota_{n}(\alpha)=
[Y_n]/\alpha.\end{eqnarray*}
For $\alpha\in H^i(S,\Q)$ put  $\alpha_{n,m}:=
[\ZZ_1]/\alpha+[\ZZ_2]/\alpha=q_1^*(\iota_n(\alpha))+q_2^*(\iota_n(\alpha))
\in H^{4-i}(T,\Q)$.
\end{defn}

The map $\iota_{n}$ is in fact easy to describe:

\begin{lem} \label{hilbkuenn}
\begin{enumerate}
\item $[Z_n(S)]=(\id_S\times h_n)^*([Y_n])$.
\item $(\id_S\times\phi_n)^*([Y_n])=\sum_i[\Delta_i]$
\item For $\alpha\in H^i(S,\Q)$ we have
$\iota_n(\alpha)=h_n^*(\bar\iota_n(\alpha))$
and $\phi_n^*(\bar\iota_n(\alpha))=\sum_{i=1}^n p_i^*(\check
\alpha)$, where $\check
\alpha$ is the Poincar\'e dual of $\alpha$.
\end{enumerate} \end{lem}

\begin{pf}
(1).
Out of codimension $3$  on $S\times \Hilb^n(S)$ we have $\oo_{Z_n(S)}=
(\id_S\times h_n)^*(\oo_{Y_n})$. So we get
$[Z_n(S)]=(\id_S\times h_n)^*([Y_n]).$
Out of codimension $3$ we also have $(\id_S\times \phi_n)^*(\oo_{Y_n})=
\bigoplus_i\oo_{\Delta_{i}}.$
Therefore (2) follows in the same way as (1).
(3) follows immediately from (1) and (2).\end{pf}

\begin{rem}\label{restch} For the total Chern classes we have
$c((\id_S\times j)^*\V_-)=(\id_S\times j)^*c(\V_-)$ and
$c((\id_S\times j)^*\V_+)=(\id_S\times j)^*c(\V_+)$, where,
as above, $j:D\maps \widetilde M$ is the embedding of the exceptional divisor.
\end{rem}

\begin{pf}
We have to see that $\Tor_k(\V_-,\oo_{S\times D})=0$
for all $k>0$
(and similarly for $\V_+$). This follows however easily from the flatness
of $\V_-$ over $\widetilde M_-$.
\end{pf}

\begin{lem} \label{nu}\begin{enumerate}
\item Assume that we are in case (1) of \ref{newcomp}, i.e.
 the change of moduli is birational.
Then, for $\alpha\in H_2(S,\Q)$, we have
\begin{eqnarray*}\nu_+(\alpha)-\nu_-(\alpha)
&=&-{1\over 2}\<\xi,\alpha\>[D],\\
\nu_{+}(pt)-\nu_{-}(\pt)&=&{1\over 4}j_*([\tau_-]-[\tau_+]).
\end{eqnarray*}
\item If $\EE_+=\emptyset$ then
\begin{eqnarray*} \sigma_-(\alpha)&=& {1\over 2}\<\xi,\alpha\>[\tau_-],\\
\sigma_{-}(\pt)&=&-{1\over 4}[\tau_-]^2.
\end{eqnarray*}
\end{enumerate}\end{lem}

\begin{pf}By (\ref{dia1}) we have the sequence
$$0\maps \V_-\maps \V_+(D_S)\maps \F_{1D}(\lambda+\tau_-+\tau_+)\maps 0.$$
Using Riemann-Roch without denominators \cite{Jo} we get
\begin{eqnarray*} c_1(\F_{1D}(\lambda+\tau_-+\tau_+))&=&[D_S]\\
c_2(\F_{1D}(\lambda+\tau_-+\tau_+))&=&
-c_1(\F_{1D}(\lambda)),
\end{eqnarray*}
and thus
\begin{eqnarray*}
c_1(\V_+(D_S))&=&c_1(\V_-)+[D_S],\\
c_2(\V_+(D_S))&=&c_2(V_-)+[D_S]\cdot
c_1(V_-)-(\id_S\times j)_*(c_1(\F_{1D}(\lambda)))\\
4c_2(\V_+)-c_1(\V_+)^2&=&4c_2(\V_+(D_S))-c_1(\V_+(D_S))^2\\
&=&
4c_2(\V_-)-c_1(\V_-)^2+2[D_S]\cdot c_1(\V_-)
-[D_S]^2-4(\id_S\times j)_*(c_1(\F_{1D}(\lambda))).
\end{eqnarray*}
Let $\alpha\in H_2(S,\Q)$. As $[D_S]$ is the  pull-back of $[D]$ from
$\widetilde M$, we have
$$([D_S]\cdot c_1(\V_-))/\alpha=[D](c_1(\V_-)/\alpha)=\<c_1,\alpha\>[D].$$
Furthermore
$(\id_S\times j)_*c_1(\F_{1D}(\lambda))/\alpha=
j_*(c_1(\F_1(\lambda))/\alpha)$,
where the second slant product is taken on $S\times D$ and
$c_1(\F_1)=\pi_D^*(c_1(\F_1'))=p^*([F_1])$.
So we get $(\id_S\times j)_\ast (c_1(\F_1)/\alpha)=\<F_1\cdot \alpha\>[D]$.
As $\lambda$ is the pull-back of a divisor on $D$, we have
$(\id_S\times j)_*c_1(\lambda)/\alpha= 0$ and
 similarly $[D_S]^2/\alpha=0$.
So we get $\nu_+(\alpha)-\nu_-(\alpha)=-{1\over 2}\<\xi,\alpha\>[D]$

By $c_1(\F_1)=p^*([F_1])$, $c_1(\F_2)=p^*([F_2])$, we get
$c_1(\F_1')/\pt=c_1(\F_2')/\pt=0.$
Then  the sequence
$$0\maps \F_{1D}(\lambda)\maps \V_-|_D\maps \F_{2D}(\tau_-+\lambda)\maps 0,$$
and remark \ref{restch} give
\begin{eqnarray*}
(c_1(\V_-)\cdot [D_S])/\pt &=&
(\id_S\times j)_*(c_1(\V_-|_D))/\pt=j_*([\tau_-]+2[\lambda]),\\
(c_1(\F_{1D}(\lambda)))/\pt &=& j_*(c_1(\F_1(\lambda))/\pt)=
j_*([\lambda]),\\{}
[D_S]^2/\pt &=& [D]^2.
\end{eqnarray*}

So we get
\begin{eqnarray*}\nu_+(\pt)-\nu_-(\pt)&=&{1\over 4}([D]^2
+j_*([2\tau_-+4\lambda])-4j_*([\lambda])) \\
&=&
{1\over 4}j_*([\tau_-]-[\tau_+]).\end{eqnarray*}
(2) can be shown using essentially the same arguments.
\end{pf}

\begin{lem} \label{change1} Let $l+2r=N$.
\begin{enumerate}
\item If we are in case (1)
 of \ref{newcomp}, then
\begin{eqnarray*}
&&\Phi_{+,l,r}(\alpha)-\Phi_{-,l,r}(\alpha)\\&&
\quad=\sum_{b=0}^{l}\sum_{c=0}^r (-1)^{r-c+1}2^{b+2c-N}{l\choose b}
{r\choose c}\<\xi,\al\>^{l-b}\intt_D
\left(\al_{n,m}^{b}\pt_{n,m}^{c}\sum_{s+t=N-b-2c-1}
(-\tau_+)^s\tau_-^t\right)\end{eqnarray*}
\item $\EE_+=\emptyset$, then
\begin{eqnarray*}
&&\Phi_{+,l,r}(\alpha)-\Phi_{-,l,r}(\alpha)\\&&
\quad=\sum_{b=0}^{l}\sum_{c=0}^r (-1)^{r-c+1}2^{b+2c-N}{l\choose b}
{r\choose c}\<\xi,\al\>^{l-b}\intt_{\EE_-}
\left(\al_{n,m}^{b}\pt_{n,m}^{c}\tau_-^{N-b-2c}\right)\end{eqnarray*}
\end{enumerate}\end{lem}

\begin{pf} (1) By remark \ref{restch} we get for $\alpha\in H_{i}(S,\Q)$
that
$[D]\cdot \nu_+(\alpha)= j_*((4c_2(\V_+|_D)-c_1(\V_+|_D)^2)/4\alpha$
(and similar for $\nu_-$).
By the sequences
\begin{eqnarray*}&&0\maps \F_{2D}(-\tau_++\lambda)\maps \V_+|_{D_S}
\maps \F_{1D}(\lambda)\maps 0\\
&&0\maps \F_{1D}(\lambda)\maps \V_-|_{D_S}\maps \F_{2D}(\tau_-+\lambda)\maps
0\end{eqnarray*}
we get
\begin{eqnarray*}&&4c_2(\V_+|_D)-c_1(\V_+|_D)^2
=4(c_2(\F_{1D})+c_2(\F_{2D}))-(c_1(\F_{2D})-c_1(\F_{1D})-[\tau_+])^2\\
&&4c_2(\V_-|D)-c_1(\V_-|D)^2
=4(c_2(\F_{1D})+c_2(\F_{2D}))-(c_1(\F_{1D})-c_1(\F_{2D})-[\tau_-])^2
\end{eqnarray*}
By the above we have
$c_1(\F_{1D})=p^*([F_1])$, $c_1(\F_{2D})=p^*([F_2])$,
$c_2(\F_{1D})=(\id_S\times \pi_D)^*(c_2(\I_{\ZZ_1}))=(\id_S\times
\pi_D)^*([\ZZ_1])$ and
$c_2(\F_{2D})=(\id_S\times \pi_D)^*([\ZZ_2])$, where, as above,  $\pi_D:D\maps
T$ is
the projection.
So we have
\begin{eqnarray*}4c_2(\V_+|_D)-c_1(\V_+|_D)^2
&=&4(\id_S\times \pi_D)^*([\ZZ_1]+[\ZZ_2])-(p^*(\xi)+[\tau_+])^2,\\
4c_2(\V_-|_D)-c_1(\V_-|_D)^2
&=&4(\id_S\times \pi_D)^*([\ZZ_1]+[\ZZ_2])-(p^*(\xi)-[\tau_-])^2,
\end{eqnarray*}
and thus for $\alpha\in H_2(S,\Q)$:
\begin{eqnarray*}
j^*(\nu_+(\alpha))&=&\alpha_{n,m}+{1\over 2}\<\xi,\al\>[-\tau_+]\\
j^*(\nu_-(\alpha))&=&\alpha_{n,m}+{1\over 2}\<\xi,\al\>[\tau_-]\\
j^*(\nu_+(\pt))&=&\pt_{n,m}-{1\over 4}[\tau_+]^2\\
j^*(\nu_-(\pt))&=&\pt_{n,m}-{1\over 4}[\tau_-]^2
\end{eqnarray*}
We write
\begin{eqnarray*}
&&\Phi_{+,l,r}(\al)-\Phi_{-,l,r}(\alpha)=\\
&&\qquad=\intt_{\widetilde M}
\left(\nu_+(\al)^{l}(\nu_+(\pt)^r-\nu_-(\pt)^r)+
\nu_+(\pt)^r(\nu_+(\al)^{l}-\nu_-(\al)^{l})\right)\\
&&\qquad=\intt_D\left({1\over 4}([-\tau_+]+[\tau_-])
j^*\left(\sum_{s+t=r-1}\nu_+(\pt)^s\nu_-(\pt)^t
\nu_+(\alpha)^l\right)
\right.\\
&&\qquad\qquad-\left.
{1\over 2}\<\xi,\al\>\
j^*\left(\sum_{s+t=l-1}\nu_+(\alpha)^s \nu_-(\alpha)^t
\nu_-(\pt)^r\right)
\right).\end{eqnarray*}
Now the claim follows after a straightforward computation.
(2) follows easily from lemma \ref{nu}(2).\end{pf}

\begin{prop}\label{donmin}
\begin{enumerate}
\item If  $S$ is a $K3$ surface and $N>0$, then
$\Phi_+=\Phi_-$.
\item If $p_g(S)=0$, then for $\alpha\in H_2(S,\Q)$ and $l,r$
with $l+2r=N$ we have
\begin{eqnarray*}
&&\Phi_{+,l,r}(\alpha)-\Phi_{-,l,r}(\alpha)\\
&&\quad=\sum_{b=0}^{l}\sum_{c=0}^r (-1)^{r-c+e_-}2^{b+2c-N}{l\choose b}
{r\choose c}\<\xi,\al\>^{l-b}\intt_{T}
\left(\al_{n,m}^{b} \pt_{n,m}^{c}
s_{2d-b-2c}(\A_+'\oplus{\A_-'}^{\vee})\right)
\end{eqnarray*}
\end{enumerate}
\end{prop}

\begin{pf}
(1) It easy to show using Riemann-Roch, that the condition
$N>0$ implies  $e_->1$ and $e_+>1$. Therefore, as
$\al_{n,m}$ and $\pt_{n,m}$ are pull-backs from $T$,
it is enough to show that for $k\le e_-+e_+-2$
we have
$$(\pi_D)_*\left(\sum_{s+t=k}(-\tau_+)^s \tau_-^t\right)=0.$$
Now $D$ is the projectivisation $\P(Q)$ where
$Q=\A_-/\tau_-$ over $\P_{-}=\P(\A_-')$. Therefore
\begin{eqnarray*}
(\pi_D)_*\left(\sum_{s+t=k}(-\tau_+)^s\tau_-^t\right)&=&
(\pi_{-})_*\left(\sum_{s+t=k}s_{s-e_++2}(Q)\tau_-^t\right)\\
&=&(\pi_{-})_*(s_{k-e_++2}(\A_-))\\
&=&(\pi_{-})_*\pi_-^*(s_{k-e_++2}(\A_-'))=0.\end{eqnarray*}
Here $\pi_{-}:\P_-\maps T$ is the projection.
(2) We just note that
$\pi_{+})_*((-\tau_+)^k)=s_{k-e_++1}(\A_+')$ and
 $(\pi_{-})_*(\tau_-^k)=(-1)^{e_-+1}s_{k-e_-+1}({\A_-'}^{\vee})$.
Then the result follows immediately from the definitions and
lemma \ref{change1}.
\end{pf}

For the rest of the chapter we assume that  $p_g(S)=q(S)=0$.
On the other hand we allow  $NS(S)=H^2(S,\Z)$ to contain torsion.

\begin{defn}\label{defchange}
Let $\xi\in H^2(S,\Z)$ be a class defining a good wall of type
$(c_1,c_2)$.
Let
$d_\xi:=(4c_1-c_1^2+\xi^2)/4$, $e_{\xi}:=-\<\xi\cdot(\xi-K_S)\>/2+d_\xi+1$
and
$$T_\xi:=\Hilb^{d_\xi}(S\sqcup S)=
\coprod_{n+m=d_\xi}\Hilb^n(S)\times \Hilb^m(S).$$
Let $q_\xi:S\times T_\xi\maps T_\xi$ be the projection.
Let  $V_\xi$ be the sheaf $p^*(\oo_S(-\xi)\oplus\oo_S(-\xi+K_S))$
on $S\times T_\xi$.
Let $\ZZ^\xi_{1}$ (resp.$\ZZ^\xi_{2}$) be the subscheme of
$S\times T_\xi$ which restricted to each component
$S\times\Hilb^n(S)\times \Hilb^m(S)$ is the subscheme $\ZZ_1$
(resp. $\ZZ_2$)
from \ref{hilbdef}. Let $\I_{\ZZ^\xi_{1}}$, $\I_{\ZZ^\xi_{2}}$
be the corresponding ideal sheaves. For $\alpha\in H_i(S,\Q)$ let
$\widetilde\alpha\in H^{4-i}(T_\xi,\Q)$ be the class  whose
restriction to each component
$\Hilb^n(S)\times \Hilb^m(S)$ of $T_\xi$ is $\alpha_{n,m}$.
Then for all $l,r$ with $l+2r=N$ we define a map
$H_2(S,\Q)\maps \Q$ by
\begin{eqnarray*}
\delta_{\xi,l,r}(\alpha)&:=&
\sum_{b=0}^{l}\sum_{c=0}^r (-1)^{r-c+e_\xi}2^{b+2c-N}{l\choose b}
{r\choose c}\<\xi,\al\>^{l-b}\\
&&\intt_{T_\xi}
\left(\widetilde\al^{b} \widetilde\pt^{c}s_{2d_\xi-2c-b}
(\Ext^1_{q}(\I_{\ZZ^\xi_{1}},
\I_{\ZZ^\xi_{2}}\otimes V_\xi)\right)\end{eqnarray*}
\end{defn}

\begin{thm}\label{donch1}
Let $S$ be a surface with $p_g(S)=q(S)=0$.
Let $c_1\in H^2(S,\Z)$ and  $c_2\in \Z$. Assume that, if
$c_1$ is  divisible by $2$ in $H^2(S,\Z)$ then
 $(4c_2-c_1^2)$ is not divisible by $8$.
Let $W$ be a good wall of type $(c_1,c_2)$ and let $H_-$ and   $H_+$ be
ample divisors on $S$
lying in neighbouring chambers separated by
$W$. Let $n_2$ be the number of $2$-torsion points in $H^2(S,\Z)$.
Then for all  $l,r$ with $l+2r=N=(4c_2-c_1^2)-3$ we have
$$\Phi_{H_+,l,r}-\Phi_{H_-,l,r}=n_2\sum_{\xi\in A^+(W)}
\delta_{\xi,l,r}.$$
Here, as above,
$$A^+(W)=\big\{ \xi\in H^2(S,\Z)\bigm | Z\hbox{ defines the wall W and }
\<\xi\cdot H_+\> >0 \}.$$
Therefore we get for a class $\alpha\in H_2(S,\Q)$
\begin{eqnarray*}
(\gamma_{c_1,c_2,g(H_+)}- \gamma_{c_1,c_2,g(H_-)})
(\underbrace{\pt,\ldots,\pt}_r,
\underbrace{\alpha,\ldots,\alpha}_l)
=(-1)^{(c_1^2+\<c_1\cdot K_S \>)}
n_2\sum_{\xi\in A^+(W)}\delta_{\xi,l,r}(\alpha).
\end{eqnarray*}
\end{thm}

\begin{pf}
If $H^2(S,\Z)$ contains no $2$-torsion,
and $a_-<a_+$ are in neighbouring minichambers separated by a miniwall
$a$ with $A^+(a)=\{(\xi,n,m)\}$, then  proposition \ref{donmin}
computes $\Phi_{a_+,l,r}-\Phi_{a_-,l,r}$.
By Serre duality and the definitions we see that
in the notations of proposition \ref{donmin}
$\A_+'\oplus{\A_-'}^\vee=\Ext^1_{q}(\I_{\ZZ_1},\II_{\ZZ_2}\otimes V_\xi)$.
Thus, if for all miniwalls $a$ the set $A^+(a)$ consists of
only one element, the theorem follows.

If $N_2\subset H^2(S,\Z)$ is the subgroup of $2$-torsion,
then $T_\xi^{n,m}\simeq N_2\times\Hilb^n(S)\times \Hilb^m(S)$. So
the exceptional divisor in $\widetilde M$ has $n_2$
isomorphic components
(or we add $n_2$ isomorphic connected components to $\widetilde M$
or subtract them), and each component gives the same contribution
to $\Phi_{a_+,l,r}-\Phi_{a_-,l,r}$.

Assume that
$A^+(a)=\{(\xi_1,n_1,m_1),\ldots,(\xi_2,n_2,m_2)\}.$
Then, as we have seen above, the $\EE_{\xi_i}^{n_i,m_i}$
are disjoint, and, as the change
$\Phi_{a_+,l,r}-\Phi_{a_-,l,r}$ can be computed on
the exceptional divisor (or the added components),
it is just the sum of the contributions for all $(\xi_i,n_i,m_i)$.
The result now follows by adding up the contributions of
all the miniwalls.
\end{pf}

By  the results we have obtained so far, in order
to compute explicitly
the change of the Donaldson invariants, when the polarisation
passes through a good wall $W=W^\xi$, we have
first to determine the Chern classes of the
bundles
$\Ext^1_{q_\xi}(\I_{\ZZ^\xi_{1}},\II_{\ZZ^\xi_{2}}\otimes V_\xi)$
 on $T_\xi$, and then make explicit computations in the
cohomology ring of $\Hilb^d(S\sqcup S)$.

In the rest of this section we will again use the
assumptions and notations from \ref{notado}, and will
adress the first question, i.e. we express the
Chern classes of the vector bundles
$\Ext^1_{q}(\I_{\ZZ_{1}},\II_{\ZZ_{2}}\otimes V)$ on
$T=\Hilb^n(S)\times \Hilb^m(S)$, (where we have written $V:=V_\xi$)
in terms of those of ``standard bundles''.

\begin{defn} Using the projections $p:S\times T \maps S$ and
$q:S\times T\maps T$ we associate to
  a vector bundle  $U$ of rank $r$ on $S$ the vector bundles
$[U]_1:=q_*(\oo_{\ZZ_2}\otimes p^*(U))$ and
$[U]_2:=
q_*(\oo_{\ZZ_1}\otimes p^*(U))$ of ranks $rn$
(resp. $rm$) on $T$.
\end{defn}

For a Cohen-Macaulay scheme $Z$, we denote by $\omega_Z$
its dualizing sheaf.

\begin{lem}\label{extlem}
$$\Ext^2_q(\oo_{\ZZ_1},\oo_{\ZZ_2}\otimes p^*V)=
q_*(\omega_{\ZZ_1}\otimes
 \omega_T^{-1}\otimes\oo_{\ZZ_2}\otimes p^*V)$$
and $\Ext^i_q(\oo_{\ZZ_1},\oo_{\ZZ_2}\otimes p^*V)=0$ for $i\ne 2$.
\end{lem}

\begin{pf}
Let
\begin{eqnarray}\label{resol}&&0\maps B_2\maps B_1
\maps \oo_{S\times T}\maps \oo_{\ZZ_1}\maps 0\end{eqnarray}
be a locally free resolution on $S\times T$. We  apply
$\HOM(\cdot,\oo_{\ZZ_2}\otimes p^*V)$ to
obtain the complex
$$0\maps \oo_{\ZZ_2}\otimes p^*V\maps B_1^*\otimes
\oo_{\ZZ_2}\otimes p^*V\maps
B_2^*\otimes \oo_{\ZZ_2}\otimes p^*V\maps
 0,$$
whose cohomologies are the $\EXT^i(\oo_{\ZZ_1},\oo_{\ZZ_2}\otimes p^*V)$.
We can arrive at this complex  differently,
 namely by first dualizing and then tensorizing by  $\oo_{\ZZ_2}\otimes p^*V$.
By dualizing and using that $\ZZ_1$ is Cohen-Macauley we obtain
$$0\maps \oo_{S\times T}\maps B_1^*\maps B_2^*\maps
 \omega_{\ZZ_1}\otimes \omega_T^{-1}\maps 0.$$
Tensorizing by $\oo_{\ZZ_2}\otimes p^*V$ gives the sequence
$$0\maps \oo_{\ZZ_2}\otimes p^*V\maps B_1^*\otimes \oo_{\ZZ_2}\otimes p^*V\maps
B_2^*\otimes \oo_{\ZZ_2}\otimes p^*V\maps
\omega_{\ZZ_1}\otimes
 \omega_T^{-1}\otimes\oo_{\ZZ_2}\otimes p^*V
\maps 0,$$
which is exact  by
the corollaire on p. V.20 in \cite{Se} because
$\ZZ_1$ and $\ZZ_2$ are Cohen-Macaulay and intersect properly.
Hence $\EXT^2(\oo_{\ZZ_1},\oo_{\ZZ_2}\otimes p^*V)=
\omega_{\ZZ_1}\otimes
 \omega_T^{-1}\otimes\oo_{\ZZ_2}\otimes p^*V$ and
$\EXT^i(\oo_{\ZZ_1},\oo_{\ZZ_2}\otimes p^*V)=0$ for $i<2$.
As $\ZZ_2$ and $\ZZ_1$ are flat of dimension $0$ over $T$,
the result follows
by applying $q_*$.
\end{pf}

\begin{prop}\label{grot}
In the Grothendieck ring of sheaves on $T$ we have the equality
\begin{eqnarray*}\Ext^1_q(\I_{\ZZ_1},\I_{\ZZ_2}\otimes p^*V)&=&
[V]_2+([V^{\vee}(K_S)]_1)^{\vee}
+(H^1(S,\oo_S(-\xi))\oplus H^1(S,\oo_S(-\xi+K_S)))\otimes\oo_T\\
&&\qquad -
q_*(\omega_{\ZZ_1}\otimes\omega_T^{-1}\otimes\oo_{\ZZ_2}\otimes p^*V).
\end{eqnarray*}
\end{prop}
\begin{pf}{\it Case n=0:}
We will use repeatedly that $\xi$ defines a good wall, so in particular
$q_*(p^*V)=R^2q_*(p^*V)=0$.
We apply
$\Hom_q(\oo_{S\times T},\cdot)$ to the sequence
$$0\maps \I_{\ZZ_2}\otimes p^*V\maps p^*V\maps \oo_{\ZZ_2}\otimes p^*V\maps 0$$
 to obtain
\begin{eqnarray}\label{ex1}
0\maps [V]_2\maps \Ext^1_q(\oo_{S\times T},\I_{\ZZ_2}\otimes p^*V) \maps
 R^1q_*(p^*V)\maps 0.
\end{eqnarray}
The surjectivity follows as $\ZZ_2$ is flat of dimension $0$ over $T$
and the injectivity by $q_*p^*V=0$.

\noindent {\it General case:}
We apply $\Hom_q(\cdot,\I_{\ZZ_2}\otimes p^*V)$ to the sequence
$0\maps \I_{\ZZ_1}\maps\oo_{S\times T}\maps \oo_{\ZZ_1}\maps 0$
to get
\begin{eqnarray}\label{ex2}
\qquad\qquad 0\to \Ext^1_q(\oo_{S\times T},\I_{\ZZ_2}\otimes p^*V)\to
\Ext^1_q(\I_{\ZZ_1},\I_{\ZZ_2}\otimes p^*V)\to
\Ext^2_q(\oo_{\ZZ_1},\I_{\ZZ_2}\otimes p^*V)\to 0.
\end{eqnarray}
The exactness on the left follows from the fact that
$q_*(\I_{\ZZ_2}\otimes p^*V)=0$ and so
$\Ext^1_q(\oo_{\ZZ_1},\I_{\ZZ_2}\otimes p^*V)$ is
torsion-free being  a subsheaf of the locally free sheaf
$R^1q_*(\I_{\ZZ_2}\otimes p^*V)$. Its
  support is  contained in $q(\ZZ_1\cap\ZZ_2)$ and thus it is the zero sheaf.

We apply $\Hom_q(\oo_{\ZZ_1},\cdot)$ to
$0\maps \I_{\ZZ_2}\otimes p^*V\maps p^*V\maps \oo_{\ZZ_2}\otimes p^*V\maps 0$
 and use
lemma \ref{extlem} to obtain
\begin{eqnarray}\label{ex3}
&&\qquad 0\maps \Ext^2_q(\oo_{\ZZ_1},\I_{\ZZ_2}\otimes p^*V)\maps
\Ext^2_q(\oo_{\ZZ_1},p^*V)
\maps
q_*(\omega_{\ZZ_1}\otimes\omega_T^{-1}\otimes\oo_{\ZZ_2}\otimes p^*V)
\maps 0.
\end{eqnarray}
By duality
$\Ext^2_q(\oo_{\ZZ_1},p^*V)=q_*(\oo_{\ZZ_1}\otimes
p^*(V^\vee(K_S))^\vee=[V^\vee(K_S)]_1^\vee$.
Thus the result follows by putting   \ref{ex1} to \ref{ex3} together.
\end{pf}

\def\intt{\int\limits}

\def\U{{\cal U}}
\def\spec{{\hbox{\rom{Spec}\,}}}
\def\ppt{{\Bbb P}_2^*}
\def\Pt{{\Bbb P}_2}
\def\M{{\cal M}}
\def\mmm{{m}}
\def\pt{\Pt}
\def\SS{{\cal S}}
\def\eps{\varepsilon}
\def\R{{\Bbb R}}
\def\pp{{\cal P}}
\def\A{{\cal A}}
\def\G{{\cal G}}
\def\V{{\cal V}}
\def\Q{{\Bbb Q}}
\def\Z{{\Bbb Z}}
\def\N{{\Bbb N}}
\def\P{{\Bbb P}}
\def\C{{\Bbb C}}
\def\I{{\cal I}}
\def\L{{\tilde L}}
\def\X#1#2{X(#1,#2)}
\def\id{{\hbox{\rom{id}}}}
\def\Coh{{\cal Coh}}
\def\II{{\cal I}}
\def\F{{\cal F}}
\def\E{{\cal E}}
\def\OO{{\cal O}}
\def\oo{\OO}
\def\PP{{\cal P}}
\def\eps{\varepsilon}
\def\Ext{{\hbox{\rom{Ext}}}}
\def\Hom{{\hbox{\rom{Hom}}}}
\def\EXT{{\cal Ext}}
\def\HOM{{\cal Hom}}
\def\hxn{{X^{[n]}}}
\def\hsn{{S^{[n]}}}
\def\xn{{X(n)}}
\def\xtn{{\widetilde \xn}}
\def\uberschrift#1{\bigskip\bigskip {\bf #1}\medskip}
\def\Theor#1#2{\medskip \noindent{{\bf Theorem~#1.~~}}{\it #2}\medskip}
\def\Thm#1#2{\medskip \noindent{{\bf Theorem~#1.~~}}{\it #2}\medskip}
\def\Prop#1#2{\medskip \noindent{{\bf Proposition~#1.~~}}{\it #2}\medskip}
\def\Cor#1#2{\medskip \noindent{{\bf Corollary~#1.~~}}{\it #2}\medskip}
\def\Lemma#1#2{\medskip \noindent{{\bf Lemma~#1.~~}}{\it #2}\medskip}
\def\Def#1{\medskip \noindent{{\bf Definition~#1.~~}}}
\def\Not#1{\medskip \noindent{{\bf Notation~#1.~~}}}
\def\Rmk#1{\medskip \noindent{{\bf Remark~#1.~~}}}
\def\Rem#1{\Rmk{#1}}
\def\Ex#1{\medskip \noindent{{\bf Example~#1.~~}}}
\def\Parag#1#2{\bigskip \noindent {\S}{\bf #1.~#2.} \medskip}
\def\Pf{\demo{Proof}}
\def\acapo{\vskip 1pt \noindent}
\def\whsq{\vbox to 5.8pt
{\offinterlineskip\hrule
\hbox to 5.8pt{\vrule height
5.1pt\hss\vrule height 5.1pt}\hrule}}
\def\blsq{\vrule height 8pt width 6pt}
\def\Qed{{\hfill {\whsq}}\enddemo}
\def\Cvd{\Qed}
\def\ppp{{\Bbb P}}
\def\maps{\longrightarrow}
\def\oo{{\cal O}}
\def\ii{{\cal I}}
\def\diagramm#1{
\def\normalbaselines{\baselineskip18pt\lineskip7pt\lineskiplimit7pt}
 \matrix{#1}}

\def\mapr#1{\smash{
 \mathop{\longrightarrow}\limits^{#1}}}

\def\mapsr#1{{
 \mathop{\longrightarrow}\limits^{#1}}}

\def\mapinto#1{\smash{
 \mathop{\longrighthook}\limits^{#1}}}

\def\mapl#1{\smash{
 \mathop{\longleftarrow}\limits^{#1}}}

\def\mapd#1{\Big\downarrow
\rlap{$\vcenter{\hbox{$\scriptstyle#1$}}$}}
\def\mapu#1{\Big\uparrow
\rlap{$\vcenter{\hbox{$\scriptstyle#1$}}$}}

\def\mapsec#1#2{\mapdownl#1\mapup#2}

\def\mapne#1{\nearrow
\rlap{$\vcenter{\hbox{$\scriptstyle#1$}}$}}
\def\mapse#1{\searrow
\rlap{$\vcenter{\hbox{$\scriptstyle#1$}}$}}
\def\mapnw#1{\nwarrow
\rlap{$\vcenter{\hbox{$\scriptstyle#1$}}$}}
\def\mapsw#1{\swarrow
\rlap{$\vcenter{\hbox{$\scriptstyle#1$}}$}}

\def\relgl{\mathrel{\smash=}}
\def\longleich{\relgl\joinrel\relgl\joinrel\relgl}

\def\sn{\Hilb^n(S)}
\def\snn{{S^{[n-1,n]}}}
\def\sym{S^{(n)}}
\def\H{{\cal H}}
\def\C{{\Bbb C}}
\def\R{{\Bbb R}}
\def\P{{\Bbb P}}
\def\N{{\Bbb N}}
\def\PP{{\cal P}}
\def\Z{{\Bbb Z}}
\def\Q{{\Bbb Q}}
\def\T{{\cal T}}
\def\eps{\varepsilon}
\def\Oplus{\textstyle\bigoplus\limits}
\def\phi{\varphi}
\def\cinf{C^{\infty}}
\def\Symg{\frak S}

\def\({\left(}
\def\){\right)}
\def\pB{{\cal B}}
\def\sB{{\cal B}^\prime}
\def\Supp{\hbox{Supp}\,}

\def\part{P(n)}
\def\bb{{\text{\rom{\bf b}}}}
\def\hsb{S^{(\bb)}}

\def\EE{{\text{\rom{\bf E}}}}

\def\al{\alpha}
\def\Ga{\Gamma}
\def\bGa{\bar \Gamma}
\def\ga{\gamma}
\def\peta{\bar \eta}

\def\ZZ{{\cal Z}}

\def\pt{{pt}}
\def\al{\alpha}

\def\Sn{{S^{(n)}}}
\def\Sm{{S^{(m)}}}
\def\<{\langle}
\def\>{\rangle}
\def\Hilb{{\text{\rom{Hilb}}}}
\def\W{{\cal W}}
\def\Tor{{\hbox{\rom{Tor}}}}
\def\Ext{{\hbox{\rom{Ext}}}}

\def\stil#1{\tilde S^#1}
\def\pigr{\pi}
\def\st{{{\hbox{$*$}}}}
\def\gi{g}
\def\Fi{\phi}
\def\ti{\tilde}
\def\alb{{\bar\al}}

\section{Explicit computations on Hilbert schemes of points}

The aim of this section is to make theorem \ref{donch1} more explicit.
We want to compute the contributions $\delta_{\xi}$ to the change
of the Donaldson invariants for
a class $\xi$  defining a good wall, in terms of cohomology classes
 and intersection
numbers  on $S$.
We do not succeed in determining $\delta_{\xi}$ completely. It turns
however out that $\delta_{\xi}$ can be developed in terms of powers of $\xi$
and we will compute the six lowest order terms (as predicted by the
conjecture of Kotschick and Morgan half of them are zero).

\begin{nota} \label{xino}
In this section we fix a class $\xi\in H^2(S,\Z)$ which defines a good wall
of type $(c_1,c_2)$ and will therefore drop $\xi$ in our notation.
In particular we write
$d:=d_\xi$, $e:=e_\xi$ and  $T:=\Hilb^{d}(S\sqcup S)$. As usual
let $p$ and $q$
be the projections of $S\times T$ to $S$ and $T$ respectively. We write
$V:=\oo_S(-\xi)\oplus\oo_S(-\xi+K_S)$, $\ZZ_1:=\ZZ_1^\xi$,
$\ZZ_2:=\ZZ_2^\xi$ and $\delta_{l,r}:=\delta_{\xi,l,r}$.
We put $\Gamma:=q_*(\omega_{\ZZ_1}\otimes\omega_{S\times T}^{-1}\otimes
\oo_{\ZZ_2}\otimes p^*V)$.
\end{nota}

We see by theorem \ref{donch1} that, in order to compute the change
$\delta_{l,r}$, it is enough to compute $\int_T
s(\Ext^1_q(\I_{\ZZ_1},\I_{\ZZ_2}\otimes p^*V))\cdot \gamma$
for all classes $\gamma\in H^*(T,\Q)$
which are pull-backs from $S^{(d)}$ via the natural map
$\Hilb^d(S\sqcup S)\maps (S\sqcup S)^{(d)}\maps S^{(d)}$.
By proposition \ref{grot}
we have
\begin{eqnarray*}
\int_T s(\Ext^1_q(\I_{\ZZ_1},\I_{\ZZ_2}\otimes p^*V))\cdot\gamma&=&
\int_T s([V^\vee(K_S)]_1^\vee \oplus [V]_2)\cdot
\gamma\\
&&\quad +\int_T (c(\Gamma)-1)s([V^\vee(K_S)]_1^\vee
\oplus [V]_2)\cdot\gamma.\\
\end{eqnarray*}

In the first part of this section we compute the first of these two
integrals.
As said in the beginning of this section, we only want to compute
the terms of lowest order of the change of the Donaldson invariants.
This corresponds to restricting our attention to a big open
subset of the Hilbert scheme of points.

\begin{nota} A point $\sigma\in \Sn$ is a formal linear
combination  $\sum_i m_i x_i$
of points on $S$ with positive integer coefficients and $\sum_i m_i=n$.
The support $supp(\sigma)$ is the set of points $x_i$. For all $i\le n$
let
$$S^{(n)}_i:=\big\{\sigma\in S^{(n)}\bigm|
\# supp(\sigma)\ge n-i+1\big\}.$$
Furthermore,  for any variety $X$ with a canonical morphism
$f:X\maps S^{(n)}$, we denote $f^{-1}S^{(n)}_i$ by $X_i$.
For the universal family $Z_n(S)\subset S\times \Hilb^n(S)$ we denote by
$Z_n(S)_i$ the preimage of $\Hilb^n(S)_i$.
\end{nota}

In order to compute the first integral we will use an inductive approach,
which is based on results of \cite{E1},\cite{F-G} and which is
similar to computations in
\cite{Go2} on the Hilbert scheme of $3$ points.

\begin{defn}(\cite{E1},\cite{F-G})
Let $S^{[n-1,n]}\maps S\times \Hilb^{n-1}(S)$ be the blow-up
along the universal
family $Z_{n-1}(S)$, and let $F_n$ the exceptional divisor.
Contrary to our conventions in the previous section for any
vector bundle $E$ on $S$ we will denote
by $E[n]$ the vector bundle $q_*(\oo_{Z_n(S)}\otimes p^* E)$ on
 $\Hilb^n(S)$.
\end{defn}

\begin{thm} \label{snn}(\cite{E1})
$S^{[n-1,n]}$ is smooth. There is a natural morphism
$S^{[n-1,n]}\maps \Hilb^n(S)$,
and on $S^{[n-1,n]}$ we have an exact sequence
\begin{eqnarray}\label{hilbseq}
0\maps V(-F_n)\maps V[n]\maps V[{n-1}]\maps 0,
\end{eqnarray}
where we have used convention \ref{convent}.
\end{thm}

It is easy to see that the induced map $S^{[n-1,n]}\maps S\times \Hilb^n(S)$
factors through $Z_n(S)\subset  S\times \Hilb^n(S)$, and that
the map $S^{[n-1,n]}\maps Z_n(S)$ is an isomorphism over
the open set
$Z_n(S)_1$. We denote by $S^{[n-1,n]}_i$ the preimage of $Z_n(S)_i$.

\begin{lem}\label{hilbnorm}
Let $N^\vee_n$ be the conormal sheaf of $Z_n(S)$ in $S\times \Hilb^n(S)$.
Then we have an exact sequence on  $S^{[n-1,n]}_2$
\begin{eqnarray*}\label{hilbnseq}
0 \maps N^\vee_n\maps T_S^\vee\maps \oo_{F_n}(-F_n)\maps 0.
\end{eqnarray*}
Here we have used the convention \ref{convent}.
In particular on  $S^{[n-1,n]}_2$ we get
$$s(N^\vee_n)=s(T^\vee_S){1-F_n\over 1-2F_n}.$$
\end{lem}
\begin{pf}
It is easy to see that
$S^{[n-1,n]}_2\maps \Hilb^n(S)_2$ is a branched $n$-fold cover,
\'etale out of $F_n$ and with ramification of order $1$ along
$F_n$.
So the result follows in the same way as in the proof of (\cite{F-G},
 lemma  2.10).
\end{pf}

\begin{lem} \label{discard} Let $i$ be a positive integer
 and assume that
$\alpha_1,\alpha_2\in A^{*}(\Hilb^n(S))$
have the same pull-back to $ \Hilb^n(S)_i$. Then
$$\intt_{\Hilb^n(S)}\alpha_1\cdot \beta=
\intt_{\Hilb^n(S)}\alpha_2\cdot \beta$$
for all $\beta\in H^{4n-4i-2}(S^{(n)},\Q)$.
The same result holds if we replace
$\Hilb^n(S)_i$ by $S^{[n-1,n]}_i$.
\end{lem}
\begin{pf}
Let $j:\Hilb^n(S)\setminus \Hilb^n(S)_i\maps Hilb^n(S)$ be the inclusion.
We get $\alpha_1=\alpha_2+j_{*}(\mu)$ for a class
$\mu\in A^*(\Hilb^n(S)\setminus \Hilb^n(S)_i).$ As the codimension of the
complement of $S^{(n)}_i$ in
$S^{(n)}$ is $2i$, the result follows by the projection formula.
\end{pf}

\begin{nota} \label{tplusminus}
For all $l\ge 1$ we denote by $\Delta_l$ the "small" diagonal
$\{(x,\ldots,x)\ | \   x\in S\}$
and by $[\Delta_l]$ its cohomology class.

We define classes $t_{1-},t_{2-},t_{3-}\in H^*(S,\Q)$ by
\begin{eqnarray*}
t_{1-}&:=&1+(2\xi-K_S)+(3\xi^2-3\xi K_S+K_S^2),\\
t_{2-}&:=&3+(18\xi-13K_S)+(63\xi^2-91\xi K_S+33K_S^2+5s_2(S)),\\
t_{3-}&:=&27+(270\xi-237K_S).
\end{eqnarray*}
Here $s_i(S):=s_i(T_S)$ is the $i^{th}$ Segre class of $S$.
We define $t_{1+}$, $t_{2+}$, $t_{3+}$ by replacing $K_S$ by $(-K_S)$
 in the definition
of $t_{1-}$, $t_{2-}$, $t_{3-}$ respectively
and put $t_i:=t_{i-}+t_{i+},$ i.e.
\begin{eqnarray*}
t_1&=&2+4\xi+6\xi^2+2K_S^2,\\
t_2&=&6+36\xi+126\xi^2+66K_S^2+10s_2(S),\\
t_3&=&54+540\xi.
\end{eqnarray*}
\end{nota}

\begin{lem}\label{hilbind}
Let $\gamma\in H^{4n-2k}(S^{(n)},\Q)$ with $k\le 5$.
Then
\begin{eqnarray*}
n\intt_{\Hilb^n(S)}s(V[n])\cdot\gamma
=\sum_{l=1}^3\intt_{S^l\times
\Hilb^{n-l}(S)} (-1)^{l-1} [\Delta_l]p_1^*t_{l-}\cdot s(V[n-l])\cdot
\gamma,
\end{eqnarray*}
where $p_1:S^l\maps S$ is the projection to the first factor.
\end{lem}
\begin{pf}
By theorem  \ref{snn} we have the identity
$s(V[n])=s(V(-F_n))s(V[n-1])$  on $S^{[n-1,n]}$ and furthermore
$$s(V(-F_n))=\sum_{i,j\ge 0}{i+j+1\choose i+1} s_i(V) F_n^j.$$
So we get
\begin{eqnarray}
n\int_{\Hilb^n(S)}s(V[n])\cdot\gamma
&=&\intt_{S\times\Hilb^{n-1}(S)}s(V)s(V[{n-1}])\cdot\gamma\label{ha}\\
\label{secint} +\sum_{i,j\ge 0}&&\intt_{\snn}F_n {i+j+2\choose i+1}
s_i(V)F_n^j s(V[{n-1}])\cdot\gamma.\label{hb}
\end{eqnarray}
By using  $V=\oo_S(-\xi)\oplus\oo_S(-\xi+K_S)$, we
see immediately that $s(V)=t_{1-}$.
We denote for all $i$ by $f_i$ the composition
$$S^{[n-i,n-i+1]}\to S\times \Hilb^{n-i}(S)\to S\times S^{(n-i)}
\to S^i\times S^{(n-i)} \to S^{(n)},$$ where  the second map is induced
by the diagonal map $S\maps S^i$ and put we $\gamma_i:=f_i^*(\gamma)$.
The integral (\ref{secint})
can be expressed as an integral over $F_n$.
We push it forward to $Z_{n-1}(S)\subset S\times \Hilb^{n-1}(S)$
and pull
back to $S^{[n-2,n-1]}$. Note that $f_2$ maps $S^{[n-2,n-1]}_i$ to
$S^{(n)}_{i+1}$. So we get, in view of lemma \ref{hilbnorm}
and lemma  \ref{discard},
\begin{eqnarray}(\ref{hb})=-\intt_{S^{[n-2,n-1]}}\sum_{i,j>0}
{i+j+3\choose i+1} s_i(V) s_j(N^\vee_{n-1})s(V[{n-1}])\cdot \gamma_2.
\label{hc}
\end{eqnarray}
We now  again use the identity $s(V[{n-1}])=s(V(-F_{n-1}))s(V[{n-2}])$
on  $S^{[n-2,n-1]}$ and obtain
\begin{eqnarray}
(\ref{hc})=\!\!\!\!\!
\intt_{S\times \Hilb^{n-2}(S)} \sum_{i+j+l\le 2}{i+l+3\choose i+1}
s_i(V) s_j(V)s_l(T^\vee_S) s(V[{n-2}])\cdot \gamma_2\label{hd}\\
\ \ \ \ \ \ \ \
+\!\!\!\!\!\!\!\intt_{S^{[n-2,n-1]}}\!\!\sum_{i,j,l}\!\!\!\!{i+l+3\choose i+1}
s_i(V) s(V[{n-2}])
\!\left(s_j(V(F_{n-1}))s_l(N^\vee_{n-1})-s_j(V)s_l(T^\vee_S)
\right)\!\cdot\! \gamma_2\label{he}
\end{eqnarray}
By explicit calculation and  the definition of $V$, we get
for the first integral
\begin{eqnarray*}&&\sum_{i+j+l\le 2}{i+l+3\choose i+1} s_i(V)
s_j(V)s_l(T^\vee_S)\\
&&\quad=
3+9s_1(V)-4K_S+13s_2(V)+6s_1(V)^2-14s_1(V)K_S+5s_2(S)\\
&&\quad=t_{2-}.
\end{eqnarray*}
Now we compute the integral (\ref{he}).
We use the formula
$$s(N^\vee_{n-1})=s(T^\vee_S){1-F_{n-1}\over 1-2F_{n-1}}$$
and the notation
$$2^{[l]}=\begin{cases} 1&l<0;\\
2^l&l\ge 0.\end{cases}$$
to obtain
\begin{eqnarray*}\qquad
(\ref{he})&=&- \intt_{S^{[n-2,n-1]}}\sum_{i,j_1,j_2,k_1,k_2}
F_{n-1}{i+k_1+k_2+3\choose i+1}{j_1+j_2+1\choose j_1+1}
2^{[k_2-1]} \cdot\\
&&\qquad\qquad\qquad\qquad\qquad\qquad
\cdot s_i(V)s_{j_1}(V)s_{k_1}(T^\vee_S)F_{n-1}^{j_2+k_2-1}
s(V[{n-2}])\cdot\gamma_2.
\end{eqnarray*}
This can  again be expressed as an integral over $F_{n-1}$.
We push forward to $Z_{n-2}(S)\subset S\times\Hilb^{n-2}(S)$
and then pull back to
$S^{[n-3,n-2]}$. Note that $f_3$ maps $S^{[n-3,n-2]}_i$ to
$S^{(n)}_{i+2}$. Therefore using lemma  \ref{hilbnorm} and lemma
\ref{discard}
to see that we can replace the push-forward of $F_{n-1}^l$ by the pull-back
of $s_{l-2}(T^\vee_S)$ via the projection $S^{[n-3,n-2]}\maps
S\times\Hilb^{n-3}(S)$.  We then push forward to
$S\times
\Hilb^{n-3}(S)$ and notice that by  theorem \ref{snn} and lemma
\ref{discard} we can replace the push-forward of $s(V[{n-2}])$ by
$s(V)s(V[{n-3}])$.
Putting all this together we obtain
\begin{eqnarray*}
(\ref{he})&=&
\intt_{S\times Hilb^{n-3}(S)}\sum_{i+j_1+j_2+k_1+k_2+l\le 3}
{i+k_1+k_2+3\choose i+1}{j_1+j_2+1\choose j_1+1}
2^{[k_2-1]} \cdot\\
&&\qquad\qquad\cdot
s_i(V)s_{j_1}(V)s_l(V)s_{k_1}(T^\vee_S)s_{j_2+k_2-2}(T^\vee_S)
s(V[{n-3}])\cdot\gamma_3.
\end{eqnarray*}
We obtain, again by direct calculation,
\begin{eqnarray*}
\sum_{i+j_1+j_2+k_1+k_2+l\le 3}
&&{i+k_1+k_2+3\choose i+1}{j_1+j_2+1\choose j_1+1}
2^{[k_2-1]} s_i(V)s_{j_1}(V)s_l(V)s_{k_1}(T^\vee_S)s_{j_2+k_2-2}(T^\vee_S)\\&&
=
27+135s_1(V)-102K_S=27+270\xi-237 K_S.
\end{eqnarray*}
This completes the proof.
\end{pf}

\begin{rem}\label{hilbind1}
Let $\gamma\in H^{4n-2k}(S^{(n)},\Q)$ with $k\le 5$.
Then the same proof shows
\begin{eqnarray*}
n\intt_{\Hilb^n(S)}s((V^\vee(K_S)[n])^\vee)\cdot\gamma
=\sum_{l=1}^3\intt_{S^l\times
\Hilb^{n-l}(S)} (-1)^{l-1} [\Delta_l]p_1^*t_{l+}\cdot
 s((V^\vee(K_S)[n-l])^\vee)\cdot
\gamma.
\end{eqnarray*}
\end{rem}

We will now introduce a compact notation for some  symmetric cohomology classes
on $S^n$
that will also help us in organizing our
combinatorical calculations.

\begin{defn}\label{polnota}
We denote by $\Symg_{n}$ the symmetric
group on $n$ letters,
which acts on $S^n$ by permuting the factors.
For $\alpha\in H^{2*}(S^k,\Q)$ and  $\beta\in H^{2*}(S^l,\Q)$
we define $\alpha\st \beta\in
H^{2*}(S^{k+l},\Q)^{\Symg_{k+l}}$ by putting
$$\alpha\st \beta:={1\over (k+l)!} \sum_{\sigma\in \Symg_{k+l}}
(p_{\sigma(1)}\times\ldots\times p_{\sigma(k)})^\st \alpha
\cdot
(p_{\sigma(k+1)}\times\ldots\times p_{\sigma(k+l)})^*\beta.$$
It is easy to see that $*$ is a commutative and associative operation.
We will  denote
$$\alpha^{\st k}:=\underbrace{\alpha\st \alpha\st \ldots\st \alpha}_k.$$
\end{defn}

\begin{rem} \label{starrem}
The following elementary properties of $*$ will be very
important for our further computations:
\begin{enumerate}
\item
For  $\alpha\in H^{2*}(S^k,\Q)$, $\beta\in H^{2*}(S^l,\Q)$
and $w\in  H^*(S^{k+l},\Q)^{\Symg_{k+l}}$ it follows immediately from the
symmetry of $w$ that
\begin{eqnarray*}\int_{S^{k+l}}(\alpha\st \beta)\cdot w
&=&\int_{S^{k+l}}(p_{1}\times\ldots \times p_{k})^*\alpha
\cdot
(p_{k+1}\times\ldots \times p_{k+l})^*\beta\cdot w\\
&=&\sum_{(w_1,w_2)}\int_{S^{k}}\alpha w_1
\cdot \int_{S^{l}}\beta w_2.
\end{eqnarray*}
Here
$w=\sum_{(w_1,w_2)} w_1\cdot w_2$ is the K\"unneth decomposition.
Analogous results hold if more then  two factors are multiplied via $*$.
\item Let $1$ denote the neutral element of the ring $H^*(S,\Q)$.
Then  $1^{\st k}$ is the neutral element of $H^*(S^k,\Q)$.
\item
It is also easy to see from the definitions
that $*$ fulfills the distributive law
$\alpha \st  (\beta_1+\beta_2)=\alpha \st  \beta_1+\alpha \st  \beta_1$.
In fact $+$ and $\st $  make
$\bigoplus_{n\ge 0}H^{2*}(S^n,\Q)^{\Symg_n}$ a commutative ring.
\item
In particular the binomial formula holds:
$$\sum_{k+l=n} {n\choose k}\alpha^{\st k}\st\beta^{\st l}
=(\alpha+\beta)^{\st n}.$$
\end{enumerate}
\end{rem}

\begin{nota}\label{starnota}
For a class $\alpha\in H^*(S,\Q)$ and a postive integer $i$
we denote by $(\alpha)_i:=[\Delta_i]p_1^*\alpha\in H^*(S^i,\Q)^{\Symg_i}$,
where $[\Delta_i]$ is the (small) diagonal $\{(x,\ldots,x) \ | \ x\in S\}$
in $S^i$. In particular $(\alpha)_1=\alpha$. We will in the future
write $(\alpha)_i(\beta)_j$ instead of $(\alpha)_i\st (\beta)_j$ and
$(\alpha)_i^{m}$ instead of $(\alpha)_i^{\st m}$. Furthermore we write
$\alpha^{\st m}\beta^{\st l}$ and $\alpha^{\st m}(\beta)_i$ instead of
$\alpha^{\st m}\st \beta^{\st l}$ and $\alpha^{\st m}\st (\beta)_i$.
\end{nota}

\begin{prop}\label{erstint}
Let $\gamma\in H^{4d-2k}S^{(d)},\Q)$ with $k\le 5$
 and $w\in H^{4d-2k}(S^d,\Q)$ its pull-back to $S^d$.
Then
\begin{eqnarray*}
&&d!\intt_{\Hilb^d(S\sqcup S)}
s([V^\vee(K_S)]^\vee_1\oplus [V]_2)\cdot \gamma\\
&&\qquad\qquad
 = \intt_{S^d}\left(t_1^{\st d}-{d\choose 2}(t_2)_2 t_1^{\st (d-2)}
+2{d\choose 3}(t_3)_3t_1^{\st (d-3)}
+3{d\choose 4}(t_2)_2^{2}t_1^{\st (d-4)}
\right)\cdot w.
\end{eqnarray*}
\end{prop}

\begin{pf}
Let $n,m$ be nonnegative integers with $n+m=d$.
Let
$\gamma_1\cdot\gamma_2\in H^*(S^{(n)})\times S^{(m)},\Q)\setminus 0$ be a
K\"unneth component of the pull-back of $\gamma$ via
$S^{(n)}\times S^{(m)} \maps S^{(d)}$. Let $w_1\cdot w_2\in
H^{4n-2l}(S^n,\Q)^{\Symg_n}\times H^{4n-2r}(S^m,\Q)^{\Symg_m}$
be the pull-back of
$\gamma_1\cdot\gamma_2$. Then $0\le l,r\le 5$. By an easy induction
using
lemma \ref{hilbind},  remark \ref{hilbind1}
and  remark  \ref{starrem} and ignoring all terms of codimension
$\ge 6$ we get
\begin{eqnarray*}
n!\intt_{\Hilb^n(S)}
s((V^\vee(K_S)[n])^\vee)\cdot \gamma_1&=&
\int_{S^n} P_{n}\cdot w_1 \\ {}
m!\intt_{\Hilb^m(S)}
s(V[m])\cdot \gamma_2&=&
\int_{S^m}Q_m\cdot {w_2},
\end{eqnarray*}
where
\begin{eqnarray*}
P_{n}&=&
t_{1+}^{\st n}-\sum_{i=2}^n(i-1)(t_{2+})_2t_{1+}^{\st (n-2)}+
\sum_{i=2}^{n-2}\sum_{j=i+2}^n(i-1)(j-1)(t_{2+})_2^{2}t_{1+}^{\st (n-4)}\\
&&\qquad +\sum_{i=3}^n(i-1)(i-2)(t_{3+})_3 t_{1+}^{\st (n-3)}.
\end{eqnarray*} and $Q_m$ is defined analogously to $P_n$ replacing
$n,$ $ t_{1+},$  $t_{2+}$ and $t_{3+}$ by
$m,$ $ t_{1-},$  $t_{2-}$ and $t_{3-}$ respectively.

Applying again remark  \ref{starrem} we obtain
$$n!m!\intt_{\Hilb^n(S)\times \Hilb^m(S)}
s((V^\vee(K_S)[n])^\vee)s(V[m])\cdot \gamma=
\intt_{S^d}(P_n\st Q_m)\cdot w,
$$
and thus
$$
d!\intt_{\Hilb^d(S\sqcup S)}
s([V^\vee(K_S)]^\vee_1)s([V]_2)\cdot \gamma\\
=
\sum_{n+m=d} {d\choose n} \intt_{S^d} (P_n\st
Q_m)\cdot w$$
Finally we have
\begin{eqnarray*}
&&\sum_{n+m=d} {d\choose n}  P_n\st Q_m\\
&&\quad =\sum_{n+m=d} {d\choose n}
\Bigg(t_{1+}^{\st n} t_{1-}^{\st m}-{n\choose 2}
(t_{2+})_2t_{1+}^{\st(n-2)}\st t_{1-}^{\st m}
-{m\choose 2}(t_{2-})_2t_{1+}^{*n} t_{1-}^{\st (m-2)}\\
&&\quad\quad + 3{n\choose 4}(t_{2+})_2^{ 2}
t_{1+}^{\st(n-4)}t_{1-}^{\st m}+
{n\choose 2}{m\choose 2}(t_{2+})_2 (t_{2-})_2
t_{1+}^{\st(n-2)}t_{1-}^{\st(m-2)}
+ 3{m\choose 4}(t_{2-})_2^{ 2}   t_{1+}^{\st n}  t_{1-}^{\st(m-4)}\\
&&\quad\quad +2{n\choose 3}(t_{3+})_3  t_{1+}^{\st (n-3)}
t_{1-}^{\st m}
 +2{m\choose 3}(t_{3-})_3  t_{1+}^{\st n}  t_{1-}^{\st(m-3)}\Bigg)\\
&&=\quad t_1^{\st d}-{d\choose 2}(t_2)_2  t_1^{\st (d-2)}
+3{d\choose 4}(t_2)_2^{  2}  t_1^{\st (d-4)}
+2{d\choose 3}(t_3)_3  t_1^{\st (d-3)}.
\end{eqnarray*}
\end{pf}

\def\pt{{pt}}
\def\al{\alpha}

\def\Se{{S^{(n)}}}
\def\Sz{{S^{(m)}}}

\def\Sn{{S^{(n)}}}
\def\Sm{{S^{(m)}}}
\def\<{\langle}
\def\>{\rangle}
\def\Hilb{{\text{\rom{Hilb}}}}
\def\W{{\Cal W}}
\def\Tor{{\hbox{\rom{Tor}}}}

\def\HSn{\Hilb^{n} (S)}
\def\HSm{\Hilb^{m} (S)}
\def\stiln{\tilde S^{n}}
\def\HSm{\Hilb^{m} (S)}
\def\stilm{\tilde S^{m}}
\def\pigr{\pi}
\def\gi{g}
\def\Fi{\phi}
\def\ti{\tilde}
\def\alb{{\bar\al}}
\def\Star#1{{(#1)}_*}
\def\ztt{{\ti{\cal Z_2}}}
\def\Deltaz{\Delta^0}
\def\E#1#2{E_{#1#2}}
\def\Ez#1#2{E^0_{#1#2}}
\def\bE#1#2{\bar E_{#1#2}}
\def\bEz#1#2{\bar E^0_{#1#2}}
\def\bF{\bar F}

Now we want to compute the second integral
$$\int_T(c(\Gamma)-1)s([V^\vee(K_S)]_1^\vee)s([V]_2)\cdot\gamma$$
for $\gamma\in H^{4d-2k}(S^{(d)},\Q)$ with $k\le 5$,
The conventions of \ref{convent} stay in effect.

\begin{defn}\label{diazwei}
Let $n,m$ be nonnegative integers with $n+m=d$.
We consider the following diagram
$$\CD
S\times \HSn  \times\HSm & @>q
>> &\HSn \times \HSm &@>\gi>> &\Se\times \Sz\cr
@AA{\ti\Fi}A&&@A{\Fi}AA&&@A\Fi'AA\cr
S\times \stiln \times\stilm & @>\tilde q>> &\stiln \times \stilm&@>\ti\gi>>
&S^{n}\times S^{m}\cr
\endCD$$
Here, as above, $q$ and $\tilde q$  are  the projections, $g:\HSn \times \HSm
\maps
\Se\times \Sz$ is the product of the Hilbert-Chow morphisms  and
$\phi':S^{n}\times S^{m}
\maps \Se\times \Sz$ is the product of the quotient maps,
and all the other varieties and maps are defined via pull-back.
For $i=1,2$ we put  $\widetilde Z_i:=\ti\phi^{-1}(\ZZ_i)$ and $\widetilde
Z_{1,2}:=\widetilde Z_1\cap \widetilde Z_2$,
(i.e. the scheme-theoretic intersection).
\end{defn}

\begin{nota} We denote by
\begin{eqnarray*}&&(\Se\times \Sz)_{*}
 :=\Big\{ (\sigma_+,\sigma_-)\in
\Se\times\Sz
\Bigm|\\
 &&\qquad\qquad\qquad\qquad\qquad\# supp(\sigma_+)\ge n-1,\,
 \# supp(\sigma_-)\ge m-1,\,
\# supp(\sigma_++\sigma_-)\ge
d-2\Big\}.
\end{eqnarray*}
Furthermore for all $X$ with a natural morphism
$f:X\maps \Se\times \Sz$  we denote
$X_{*}:=f^{-1}(\Se\times \Sz)_{*}$.
We put  $$\Gamma_{n,m}:=\phi^*(\Gamma|_{(\HSn\times \HSm)_{*}})$$
(see 7.1 for the definition of $\Gamma$).
For $1\le i\le n$ (resp. $1\le j\le m$) we denote by $p_{i+}$ (resp.
$p_{j-}$) the projection from
$S^{n}\times S^{m}$ onto the $i^{th}$ factor of $S^{n}$
(resp. the $j^{th}$ factor of $S^{m}$).

For $\epsilon =+,-$, $\eta =+,-$
we put
\begin{eqnarray*}\Delta_{0,i}^{\epsilon }&:=&\Bigm\{(x,x_1^+,\ldots x_{n}^+,
x_1^-,\ldots x_{m}^-)\in S\times S^{n}\times S^{m}\Bigm|
 x=x_i^\epsilon \Big\}\\
 \Delta_{i,j}^{\epsilon \eta }&:=&\Big\{(x_1^+,\ldots x_{n}^+,
x_1^-,\ldots x_{m}^-)\in  S^{n}\times S^{m}\Bigm| x_i^\epsilon
=x_j^\eta \Big\}\\
 \Delta_{0,i,j}^{\epsilon \eta }&:=&\Big\{(x,x_1^+,\ldots x_{n}^+,
x_1^-,\ldots x_{m}^-)\in S\times S^{n}\times S^{m}\Bigm|
x=x_i^\epsilon =x_j^\eta \Big\}\\
\end{eqnarray*}
We will also denote by $\Delta_{i,j}^{+-}$,
$\Delta_{0,i,j}^{+-}$,
$\Delta^{\epsilon }_{0,i}$
the pull-backs $\tilde g^{-1}(\Delta_{i,j}^{+- })$,
$(\id_S\times \tilde g)^{-1}(\Delta_{0,i,j}^{+-})$,
$(\id_S\times \tilde g)^{-1}(\Delta^{\epsilon }_{0,i})$.
We denote   $D_{i,j}:=
\tilde g^{-1}(\Delta^{++}_{i,j})$ and $E_{i,j}:=
\tilde g^{-1}(\Delta^{--}_{i,j})$.
$D_{i,j}$ and $E_{i,j}$ are divisors (see below), we denote
$F_{i}:=\sum_{j<i}D_{i,j}$ and $G_{i}:=\sum_{j<i}E_{i,j}$.
\end{nota}

\begin{rem}\label{hilbfacts}
The following easy facts will be used throughout the computation.
\begin{enumerate}
\item It is well known that $(\Hilb^{n}(S)\times \Hilb^{m}(S))_{*}$ is
obtained from
$(S^{n}\times S^{m})_{*}$ by blowing up all the $\Delta^{++}_{i,j}$
and $\Delta^{--}_{i,j}$ and taking the quotient by the action of the
product of the symmetric groups
$\Symg_{n}\times \Symg_{m}$. It follows that in fact
$(\stiln\times \stilm)_{*}$
is just the blow up of $(S^{n}\times S^{m})_{*}$ along the
(disjoint) smooth subvarieties  $(\Delta^{++}_{i,j})_*$
and $(\Delta^{--}_{i,j})_*$ and the $(D_{i,j})_*$ and $(E_{i,j})_*$ are the
exceptional
divisors.
\item It is also easy to see that
$(\widetilde Z_1)_{*}=\bigcup_{i=1}^{n} (\Delta_{0,i}^{+})_*$,
$(\widetilde Z_2)_{*}=\bigcup_{j=1}^{m} (\Delta_{0,j}^{-})_*$
and therefore
$$(\widetilde Z_{1,2})_{*}=\bigcup_{i=1}^{n}\bigcup_{j=1}^{m}
(\Delta_{0,i,j}^{+-})_{*}.$$
(We mean here the scheme theoretic union, i.e. the scheme defined by
the intersection of the ideals).
\item
For  $i\ne j$ we have (scheme-theoretically)
\begin{eqnarray*}
(\Delta_{0,i}^{+})_{*}\cap \Delta_{0,j}^{+}&=&
(\Delta_{0,i}^{+})_{*}\cap D_{i,j},\\
(\Delta_{0,i}^{-})_{*}\cap \Delta_{0,j}^{-}&=&
(\Delta_{0,i}^{-})_{*}\cap E_{i,j}.
\end{eqnarray*}
\end{enumerate}
\end{rem}

\begin{lem}\label{comalg}
\begin{enumerate}
\item Let $X$ be a smooth variety, and let $Y$ and $Z$ be
 Cohen-Macauley subschemes
of $X$ such that the ideal $\I_{Z/(Y\cup Z)}$ of $Z$ in
$Y\cup Z$ is $\oo_Y(-D)$ for a divisor $D$ on $Y$.
Then in the Grothendieck ring of $X$ we have
\begin{eqnarray*}\oo_{Y\cup Z}&=&\oo_Y(-D)+\oo_Z\ \hbox{ and }\\
\omega_X^{-1}\otimes \omega_{Y\cup Z}&=&\omega_X^{-1}\otimes\omega_Y(D)
+\omega_X^{-1} \otimes\omega_Z.\end{eqnarray*}
\item Let $f:X\maps Y$ be a morphism between smooth varieties. Let
$Z\subset Y$ be a Cohen-Macauley subscheme of codimension $2$ and
assume $W:=f^{-1}(Z)$ has pure codimension $2$ in $X$. Then
$$f^{*}(\omega_Y^{-1}\otimes \omega_Z)=\omega_X^{-1}\otimes \omega_W.$$
\item Let $X$ be a smooth variety and $Y$ and $Z$ Cohen-Macauley
 subschemes of codimension
$2$ intersecting properly. Then in the Grothendieck ring of $X$
we have
$$\oo_{Y}\otimes\oo_{Z}=\oo_{Y\cap Z}.$$
\end{enumerate}
\end{lem}

\begin{pf}
(1) The first identity follows from the standard exact sequence
$$0\maps \oo_Y(-D)\maps \oo_{Y\cup Z}\maps \oo_Z\maps 0.\eqno (*)$$
Now we dualize $(*)$ and use that for a two codimensional
Cohen-Macauley subscheme
$W\subset X$ we have
$$\EXT^i(\oo_W,\oo_X)=\begin{cases} 0& i<2,\\
\omega_X^{-1}\otimes \omega_X& i=2\end{cases}$$
to obtain the sequence
$$0\maps \omega_X^{-1}\otimes
\omega_Z\maps  \omega_X^{-1}\otimes \omega_{Y\cup Z}\maps
\omega_X^{-1}\otimes\omega_Y(F)\maps 0$$ and thus the second identity.

(2) We take a locally free resolution
$$0\maps B\maps A\maps \oo_Y\maps \oo_Z\maps 0.$$
Pulling it back we obtain the sequence
$$0\maps f^*B\maps f^*A\maps \oo_X\maps \oo_W\maps 0,$$
which stays exact by the Hilbert-Birch theorem (see e.g. \cite{P-S}
lemma 3.1).
Dualizing we obtain the exact sequence
$$0\maps \oo_X\maps f^*A\maps f^*B\maps
\omega_X^{-1} \otimes\omega_Z\maps 0.$$
We can also arrive at this sequence differently, by first dualizing and
then pulling back. This way we obtain the sequence
$$0\maps \oo_X\maps f^*A\maps f^*B\maps
f^*(\omega_X^{-1} \otimes\omega_Z)\maps 0,$$
and (2) follows.

(3) By the corollaire on p. 20 in \cite{Se} we have
$\Tor_i(\oo_Y,\oo_Z)=0$ for $i>0$, and (3) follows.
\end{pf}

\begin{lem}\label{grot2}
In the Grothendieck ring of
$(\widetilde S^{n}\times \widetilde S^{m})_{*}$
we have the equality
$$\phi^*(\Gamma_{n,m})=\sum_{i=1}^{n}\sum_{j=1}^{m}\Big(
\oo_{\Delta_{i,j}^{+-}}(F_i-G_j-p_{i+}^*\xi)+
\oo_{\Delta_{i,j}^{+-}}(F_i-G_j-p_{i+}^*(\xi+K_S))\Big).$$
\end{lem}

\begin{pf}
Using remark \ref{hilbfacts}(2) and remark \ref{hilbfacts}(3) and applying
lemma \ref{comalg}(1) inductively we obtain in the Grothendieck ring of
$S\times(\widetilde S^n\times \widetilde S^m)_*$ the equalities
\begin{eqnarray}
\label{gr1}\oo_{\widetilde Z_1}&=&\sum_{i=1}^n
\oo_{\Delta_{0,i}^+}(-F_i)\\
\label{gr2}\oo_{\widetilde Z_2}&=&\sum_{j=1}^m
\oo_{\Delta_{0,j}^-}(-G_j)\ \hbox{ and}
\end{eqnarray}
\begin{eqnarray*}
\widetilde\phi^*(\omega_{S\times\Hilb^n(S)\times\Hilb^m(S)}^{-1}\otimes
\omega_{\ZZ_1})&=&
\omega_{S\times\widetilde S^n\times \widetilde S^m}^{-1}
\otimes\omega_{\widetilde Z_1}\\
&=&\sum_{i=1}^n
 \omega_{S\times\widetilde S^n\times \widetilde S^m}^{-1}
\omega_{\Delta_{0,i}^+}(F_i)\\
&=&\sum_{i=1}^n
\oo_{\Delta_{0,i}^+}(-p_{i+}^*K_S+F_i),
\end{eqnarray*}
where in the third and the last line we have used
lemma \ref{comalg}(2).
Now using lemma \ref{comalg}(3) and tenzorizing by
$p^*V$ we obtain in the Grothendieck ring of
$S\times(\widetilde S^n\times \widetilde S^m)_*$ the equality
$$(\id_S\times \phi)^*(\omega_T^{-1}\otimes \omega_{\ZZ_1}\otimes
\oo_{\ZZ_2}\otimes p^*V)
=\sum_{i=1}^{n}\sum_{j=1}^{m}\Big(\oo_{\Delta_{0,i,j}^{+-}}(F_i-G_j
-p_{i+}^*\xi)+ \oo_{\Delta_{0,i,j}^{+-}}(F_i-G_j
-p_{i+}^*(\xi+K_S)\Big)
.$$
The morphism $\phi:(\widetilde S^{n}\times \widetilde S^{m})_{*}
\maps \Hilb^{n}(S)\times \Hilb^{m}(S)$
is  flat.
Therefore we get by (\cite{Ha} prop.III.9.3)
\begin{eqnarray*}
&&\phi^*q_*(\omega_T^{-1}\otimes \omega_{Z_1}\otimes \oo_{Z_2}\otimes p^*V)\\
&&\qquad =
\tilde q_*(\id_S\times \phi)^*(\omega_T^{-1}\otimes
\omega_{Z_1}\otimes \oo_{Z_2}\otimes p^*V)\\
&&\qquad=\sum_{i=1}^{n}\sum_{j=1}^{m}\tilde q_*(
(\oo_{\Delta_{0,i,j}^{+-}}(F_i-G_j
-p_{i+}^*\xi)+ \oo_{\Delta_{0,i,j}^{+-}}(F_i-G_j
-p_{i+}^*(\xi+K_S))\\
&&\qquad
=\sum_{i=1}^{n}\sum_{j=1}^{m}(\oo_{\Delta_{i,j}^{+-}}(F_i-G_j
-p_{i+}^*\xi)+ \oo_{\Delta_{i,j}^{+-}}(F_i-G_j -p_{i+}^*(\xi+K_S)),
\end{eqnarray*}
in the Grothendieck ring of
$(\widetilde S^{n}\times \widetilde S^{m})_{*}$. The last identity
follows from the fact that the projection
$\tilde q|_{\Delta_{0,i,j}^{+-}}:{\Delta_{0,i,j}^{+-}}\maps \Delta_{i,j}^{+-} $
 is an isomorphism.
\end{pf}

\begin{lem}\label{rrsd} Let $X$ be a smooth variety and let
$i:Y\maps X$ be the closed
embedding of a smooth subvariety of codimension $2$ with conormal bundle
$N^\vee$. Let
$D$ be a  divisor on $Y$. Then $$c(i_*(\oo_Y(-D)))=1-i_*\Big(\sum_{k,l\ge 0}
{k+l+1\choose l} D^l s_k(N^\vee)\Big).$$
\end{lem}

\begin{pf}
This is a straightforward application of Riemann-Roch without denominators
\cite{Jo}.
\end{pf}

\begin{lem}\label{cgamma}
Let $1\le k\le 5$. Then
$c_k(\Gamma_{n,m}|_{(\widetilde S^n\times \widetilde S^m)_*})$
is the part of degree $k$ of
$$\matrix
\displaystyle -\sum_{(i,j)}
[\Delta_{i,j}^{+-}]\Big(2+4p_{i+}^*\xi+4(G_j-F_i)
+p_{i+}^*(6\xi^2+3s_2(S)-K_S^2)
+12p_{i+}^*\xi(G_j-F_i)\cr
\displaystyle+6(F_i^2+G_j^2)+p_{i+}^*(24\xi^2+12s_2(S)-4K_S^2)
(G_j-F_i)
+24p_{i+}^*\xi(G_j^2+F_i^2)+8(G_j^3-F_i^3)\Big)\cr
+\displaystyle \sum_{(i,j)\ne (i_1,j_1)}
[\Delta_{i,j}^{+-}][\Delta_{i_1,j_1}^{+-}](4+8p_{i+}^*\xi+8p_{i_1+}^*\xi)
\endmatrix\eqno (7.20.1)$$
Here  $(i,j)$ and $(i_1,j_1)$ run through
$\{1,\ldots, n\}\times \{1,\ldots, m\}$.
\end{lem}

\begin{pf}
We compute on $(\widetilde S^n\times \widetilde S^m)_*$.
We notice that $[\Delta_{i,j}^{+-}]$ is just the pull-back of
the corresponding class  in $S^{n}\times S^{m}$ via $\tilde g$
and the conormal bundle of $\Delta_{i,j}^{+-}$ is just the pull-back
of the conormal bundle, i.e. $p_{i+}^*(T^\vee_S)$.
Furthermore we note that on
$(\widetilde S^{n}\times \widetilde S^{m})_{*}$ we have
$[\Delta_{i,j}^{+-}]\cdot F_i\cdot G_j=0$.
Therefore we obtain by  lemma  \ref{rrsd}
after some calculation
that for $1\le k\le 5$ the Chern class
$c_k(\oo_{\Delta_{i,j}^{+-}}(-p_{i+}^*\xi +F_i-G_j))$
is the part of degree $k$ of
\begin{eqnarray*}&&-[\Delta_{i,j}^{+-}]\Big(1+p_{i+}^*(2\xi-K_S)+2(G_j-F_i)
+p_{i+}^*(3\xi^2-3\xi K_S+s_2(S))
\\
&&\qquad+p_{i+}^*(6\xi-3K_S)(G_j-F_i)+3(G_j^2+F_i^2)\\
&&\qquad
+p_{i+}^*(12\xi^2-12K_S\xi+4s_2(S))(G_j-F_i)
+p_{i+}^*(12\xi-6K_S)(G_j^2+F_i^2)+4(G_j^3-F_i^3)\Big).\\
\end{eqnarray*}
Analogously we obtain that $c_k(\oo_{\Delta_{i,j}^{+-}}
(-p_{i+}^*(\xi+K_S) +F_i-G_j))$ is the part of degree $k$ of
\begin{eqnarray*}
&&1-[\Delta_{i,j}^{+-}]\Big(1+p_{i+}^*(2\xi+K_S)+2(G_j-F_i)
+p_{i+}^*(3\xi^2+3\xi K_S+s_2(S))\\
&&\qquad
+p_{i+}^*(6\xi+3K_S)(G_j-F_i)+3(G_j^2+F_i^2)\\
&&\qquad
+p_{i+}^*(12\xi^2+12K_S\xi+4s_2(S))(G_j-F_i)
+p_{i+}^*(12\xi+6K_S)(G_j^2+F_i^2)+4(G_j^3-F_i^3)\Big).\\
\end{eqnarray*}
We notice that
$[\widetilde\Delta_{i,j}^{+-}]^2=[\Delta_{i,j}^{+-}]p_{i+}^*(c_2(S))$.
Thus, by
multiplying out, we get that $c_k(\oo_{\Delta_{i,j}^{+-}}(-p_{i+}^*\xi
-F_i+G_j)\oplus \oo_{\Delta_{i,j}^{+-}} (-p_{i+}^*(\xi+K_S) -F_i+G_j))$
is the
part of degree $k$ of
\begin{eqnarray*}
&&1-[\Delta_{i,j}^{+-}]\Big(2+4p_{i+}^*\xi+4(G_j-F_i)\\
&&\qquad +p_{i+}^*(6\xi^2+3s_2(S)-K_S^2)+
12p_{i+}^*\xi(G_j-F_i)+6(G_j^2+F_i^2)\\
&&\qquad +p_{i+}^*(24\xi^2+12s_2(S)-4K_S^2)(G_j-F_i)
+24p_{i+}^*\xi(G_j^2+F_i^2)+8(G_j^3-F_i^3)\Big).\\
\end{eqnarray*}
Now we take the product over all $i,j$. We use that on
$(\tilde S^{n}\times\tilde S^{m})_{*}$ we have
$[\Delta_{i_1,j_1}^{+-}]\cdot [\Delta_{i_2,j_2}^{+-}]\cdot
F_i=[\Delta_{i_1,j_1}^{+-}]\cdot [\Delta_{i_2,j_2}^{+-}]\cdot G_j=0$ unless
$\{i_1,j_1\}=\{i_2,j_2\}$, and obtain the result.
\end{pf}

\begin{rem}
\label{cstand}
\begin{enumerate}
\item In the Grothendieck ring of
$(\widetilde S^{n}\times \widetilde S^{m})_{*}$ we
have
\begin{eqnarray*}
\phi^*([V]_2|_{\Hilb^n(S)\times \Hilb^m(S)})
&=& \sum_{j=1}^{m} p_{j-}^*V(-G_j),\\
\phi^*([V^\vee(K_S)]^\vee|_{\Hilb^n(S)\times \Hilb^m(S)})
&=& \sum_{i=1}^{n}  p_{i+}^*(V(-K_S))(F_i).\\
\end{eqnarray*}
\item Therefore, for $l\le 3$,
$s_l(\phi^*([V]_2|_{\Hilb^n(S)\times \Hilb^m(S)}))$
 is the part of degree $l$ of
$$\matrix\displaystyle\prod_{j=1}^{m}p_{j-}^*t_{1-}
+\sum_{1\le j\le j_1\le m}\Big(
+2E_{j,j_1}
+p_{j-}^*(10\xi-5K_S)E_{j,j_1}+3E_{j,j_1}^2)\cr
+p_{j-}^*(30\xi^2-30\xi K_S+9K_S^2)E_{j,j_1}
+p_{j-}^*(18\xi-9K_S)E_{j,j_1}^2+4E_{j,j_1}^3 \Big)
\displaystyle\prod_{j_2\not\in \{j,j_1\}} p_{j_2-}^*t_{1-},
\endmatrix\eqno (7.21.1)$$
and $s_l(\phi^*([V^\vee(K_S)])^\vee|_{\Hilb^n(S)\times \Hilb^m(S)}
)$ is the part of degree $l$ of
$$\matrix
\displaystyle\prod_{i=1}^{n}p_{i+}^*t_{1+} +
\sum_{1\le i\le i_1\le n}
\Big( -2D_{i,i_1}
-p_{i+}^*(10\xi+5K_S)D_{i,i_1}+3D_{i,i_1}^2 \cr
-p_{i+}^*(30\xi^2+30\xi K_S+9K_S^2)D_{i,i_1}
+p_{i+}^*(18\xi+9K_S)D_{i,i_1}^2-4D_{i,i_1}^3\Big)
\displaystyle\prod_{i_2\not\in \{i,i_1\}}
p_{i_2+}^*t_{1+}.
\endmatrix
\eqno (7.21.2)$$
\end{enumerate}
\end{rem}

\begin{pf} (1) follows from the formulas $\ref{gr1}$, $\ref{gr2}$
by tensorizing with $p^*V$ (resp. $p^*(V^\vee(K_S))$) and pushing
down via
$\tilde q_*$.
(2) is just a straightforward computation
using that $E_{i,j}\cdot E_{k,l}=D_{i,j}\cdot D_{k,l}=0$  for $\{i,j\}\ne
\{k,l\}$.
\end{pf}

\begin{rem}\label{discard2}
Let $k\le 5$ and $\gamma\in H^{4d-2k}(\tilde S^{n}\times \tilde
S^{m},\Q)$ and assume that $\alpha_1,\alpha_2\in A^*(\tilde S^{n}\times \tilde
S^{m})$ have the same pull-back to $(\tilde S^{n}\times \tilde
S^{m})_{*}$. Then, for all $i\le n$, $j\le m$, we get
analogously to lemma  \ref{discard}
$$\intt_{\tilde S^{n}\times \tilde
S^{m}}\Delta_{i,j}^{+-}\cdot(\alpha_1-\alpha_2)\cdot\gamma=0.$$
\end{rem}

\begin{prop}\label{intzwei}
Let $\gamma\in H^{4d-2k}(S^{(d)},\Q)$ with  $k\le 5$, and
let $w\in H^{4d-2k}(S^d,\Q)$ be the  pull-back of $\gamma$ to $S^d$.
Then
\begin{eqnarray*}
&&d!\intt_{\Hilb^d(S\sqcup S)}
(c(\Gamma)-1)s([V^\vee(K_S)]_1^\vee)s([V]_2)\cdot\gamma\\
&&\qquad =
\int_{S^d}\Big( -d(d-1)(2+12\xi+42\xi^2+3s_2(S)+K_S^2)_2
  t_1^{\st (d-2)}\\
&&\qquad\quad\qquad\quad +d(d-1)(d-2)(30+260\xi)_3  t_1^{\st (d-3)}\\
&&\qquad\quad\qquad\quad
+2d(d-1)(d-2)(d-3)(2+12\xi)_2^{  2}  t_1^{\st (d-4)}\Big)\cdot w,
\end{eqnarray*}
and, with
\begin{eqnarray*}
R_d&:=&t_1^{\st d}-d(d-1)(5+30\xi+105\xi^2+8s_2(S)+34K_S^2)_2
  t_1^{\st (d-2)}\\
&&\quad+d(d-1)(d-2)(48+440\xi)_3  t_1^{\st (d-3)}\\
&&\quad +{d(d-1)(d-2)(d-3)\over 2}
(5+30\xi)_2^{  2}  t_1^{\st (d-4)},\end{eqnarray*}
we get
\begin{eqnarray*}
d!\intt_{\Hilb^d(S\sqcup S)}
c(\Gamma)s([V^\vee(K_S)]_1^\vee)s([V]_2)\cdot\gamma
 =\intt_{S^d}R_d\cdot w.
\end{eqnarray*}
\end{prop}
\begin{pf}
We fix $n$ and $m$ with $n+m=d$ and start by computing on
$(\widetilde S^{n}\times \widetilde S^{m})$. Using remark
\ref{discard2} we can restrict our
attention to $(\widetilde S^{n}\times \widetilde S^{m})_{*}$.
We multiply out the formulas  (7.21.1),(7.21.2) and (7.20.1)
and push down to $S^{n}\times S^{m}$.
We shall use the following facts:
On $(\tilde S^{n}\times \tilde S^{m})_{*}$ any of
$D_i$ and $E_j$,
 gives zero when multiplied by
$[\Delta_{i,j}^{+-}][\Delta_{i_1,j_1}^{+-}]$.
Furthermore
$\tilde g_*(D_{i,j})=\tilde g_*E_{i,j}=0$,
$\tilde g_*(D_{i,i_1}^2)=-[\Delta_{i,i_1}^{++}]$,
$\tilde g_*(E_{j,j_1})^2=-[\Delta_{j,j_1}^{--}]$,
$\tilde g_*(D_{i,i_1}^3)=[\Delta_{i,i_1}^{++}]p_{i+}^*(K_S)$,
$\tilde g_*(E_{j,j_1}^3)=[\Delta_{j,j_1}^{--}]]p_{j-}^*(K_S)$.

Below we  collect the result of the push-down in ten terms according to
the factors that they contain {\it before} the push-down.
All the summands contain at least one diagonal factor
$[\Delta_{i,j}^{+-}]$ and at most two diagonal factors
 $[\Delta_{i,j}^{+-}]$,
$[\Delta_{i_1,j_1}^{+-}]$. The first seven terms come from summands
containing precisely one factor $[\Delta_{i,j}^{+-}]$. So to define these
summandss we can fix $i$ and $j$.   The first
term corresponds to summands not containing any exceptional divisor $D_{i,i_1}$
or  $E_{j,j_1}$.
The second to seventh summands correspond in that order
to the push-downs of the terms containing only powers of $D_{i,i_1}$ with
$i_1<i$,   $E_{j,j_1}$ with $j_1<j$, $D_{i,i_1}$ with
$i_1>i$,   $E_{j,j_1}$ with $j_1>j$, $D_{i_1,i_2}$ with
$i\not\in \{i_1,i_2\}$ and    $E_{j_1,j_2}$ with $j\not\in \{j_1,j_2\}$.
Notice that on $(\widetilde S^{n}\times \widetilde S^{m})_{*}$
the class
$[\Delta_{i,j}^{+-}]D_{i_1,i_2}E_{j_1,j_2}$ is zero for all
 $i_1,i_2,j_1,j_2$ and
$[\Delta_{i,j}^{+-}]D_{i_1,i_2}D_{i_3,i_4}
=[\Delta_{i,j}^{+-}]E_{j_1,j_2}E_{j_3,j_4}=0$
unless $\{i_1,i_2\}=\{i_3,i_4\}$ (resp. $\{j_1,j_2\}=\{j_3,j_4\}$).
The last three summands correspond to terms containing
two diagonal factors $[\Delta_{i,j}^{+-}][\Delta_{i_1,j_1}^{+-}]$.
In that order they correspond to the possibilties that
$j=j_1$, that $i=i_1$ and finally that $i\ne i_1$ and $j\ne j_1$.
After a long but elementary computation we get that, if
$k\le 5$,
$\tilde g_*(\phi^*(c_k(\Gamma_{n,m})-1)s([V^\vee(K_S)[n]^\vee)s(V[m]))$
is the part of degree $k$ of
\begin{eqnarray*}
&&\sum_{(i,j)}\Bigg(
- [\Delta_{i,j}^{+-}] p_{i+}^*(2+12\xi+42\xi^2+3s_2(S)+K_S^2)
\prod_{i_1\ne i}p_{i_1+}^*(t_{1+})
\prod_{j_1\ne j} p_{j_1-}^*(t_{1-})\\
&&\qquad\quad
+\sum_{i_1<i}
[\Delta_{i,j}^{+-}][\Delta_{i,i_1}^{++}]
p_{i+}^*(20+160\xi+90K_S)
\prod_{i_2\not \in \{i,i_1\}}p_{i_2+}^*t_{1+}
\prod_{j_1\ne j} p_{j_1-}^*t_{1-}\\
&&\qquad\quad
+\sum_{j_1<j}[\Delta_{i,j}^{+-}][\Delta_{j,j_1}^{--}]
p_{i+}^*(20+160 \xi-90K_S)
\prod_{i_1\ne i }p_{i_1+}^*t_{1+}
\prod_{j_2\not \in \{j,j_1\}}
p_{j_2-}^*t_{1-}
\\
&&\qquad\quad
+\sum_{i_1>i}
[\Delta_{i,j}^{+-}][\Delta_{i,i_1}^{++}]
p_{i+}^*(6+60\xi+20K_S)
\prod_{i_2\not \in \{i,i_1\}}p_{i_2+}^*t_{1+}
\prod_{j_1\ne j} p_{j_1-}^*t_{1-}\\
&&\qquad\quad
+\sum_{j_1>j}[\Delta_{i,j}^{+-}][\Delta_{j,j_1}^{--}]
p_{j-}^*(6+60 \xi-20K_S)
\prod_{i_1\ne i }p_{i_1+}^*t_{1+}
\prod_{j_2\not \in \{j,j_1\}}
p_{j_2-}^*t_{1-}\\
&&\qquad\quad
+\sum_{i_1\ne i}\sum_{i_2\ne i,i_2< i_1}\Bigg(
[\Delta_{i,j}^{+-}][\Delta_{i_1,i_2}^{++}]
p_{i+}^*(2+12 \xi)
p_{i_1+}^*(3+18\xi+13K_S)\\
&&\qquad\quad\cdot
\prod_{i_3\not\in\{ i,i_1,i_2\} }p_{i_3+}^*t_{1+}
\prod_{j_1\ne j}
p_{j_1-}^*t_{1-}\Bigg)
+\sum_{j_1\ne j}\sum_{j_2\ne j,j_2< j_1}\Bigg(
[\Delta_{i,j}^{+-}][\Delta_{j_1,j_2}^{--}]
p_{j-}^*(2+12 \xi)\\
&&\qquad\quad\cdot
p_{j_1-}^*(3+18\xi-13K_S)
\prod_{i_1\ne j}p_{i_1+}^*t_{1+}
\prod_{j_3\not\in\{ j,j_1,j_2\} }p_{j_3-}^*t_{1-}\Bigg)
\\
&&\qquad\quad
+\sum_{i_1<i}
[\Delta_{i,j}^{+-}][\Delta_{i_1,j}^{+-}]
p_{i+}^*(4+40\xi+4K_S)
\prod_{i_2\not \in \{i,i_1\}}p_{i_2+}^*t_{1+}
\prod_{j_1\ne j} p_{j_1-}^*t_{1-}\\
&&\qquad\quad
+\sum_{j_1<j}[\Delta_{i,j}^{+-}][\Delta_{i,j_1}^{+-}]
p_{j-}^*(4+40 \xi-4K_S)
\prod_{i_1\ne i }p_{i_1+}^*t_{1+}
\prod_{j_2\not \in \{j,j_1\}}
p_{j_2-}^*t_{1-}\\
&&\qquad\quad
+\sum_{i_1<i}\sum_{j_1\ne j}
[\Delta_{i,j}^{+-}][\Delta_{i_1,j_1}^{+-}]
p_{i+}^*(2+12 \xi)
p_{i_1+}^*(2+12\xi)
\prod_{i_2\not\in\{ i,i_1\} }p_{i_2+}^*t_{1+}
\prod_{j_2\not \in \{j,j_1\}}
p_{j_2-}^*t_{1-}\Bigg).
\end{eqnarray*}
Now we want to translate this result into the  notation \ref{starnota}.
Using remark \ref{starrem} and notation \ref{starnota} we see that for
$w\in H^{4d-2k}(S^d,\Q)^{\Symg_{d}}$ and $a\in H^*(S,\Q)$
we have
\begin{eqnarray*}
\int_{S^d}[\Delta^{+-}_{i,j}] p_{i+}^*a\prod_{i_1\ne i}p_{i_1+}^*t_{1+}
\prod_{j_1\ne i}p_{j_1-}^*t_{1-}\cdot w
&=&\int_{S^d}(a)_2  t_{1+}^{\st  (n-1)} t_{1-}^{\st (m-1)}  \cdot w.
\end{eqnarray*}
Now assume $j\ne j_1$. Then
\begin{eqnarray*}
\int_{S^d}[\Delta^{+-}_{i,j}][\Delta^{+-}_{i,j_1}]p_{i+}^*a
\prod_{i_1\ne i}
p_{i_1+}^*t_{1+}
\prod_{j_2\not\in \{j,j_1\}}p_{j_2-}^*t_{1-}\cdot w &=&
\int_{S^d}(a)_3  t_{1+}^{\st  (n-1)}  t_{1-}^{\st  (m-2)}\cdot w.
\end{eqnarray*}
We also see that
$[\Delta^{+-}_{i,j}][\Delta^{+-}_{i,j_1}]
=[\Delta^{+-}_{i,j}][\Delta^{--}_{j,j_1}]$ and
$[\Delta^{+-}_{i,j}][\Delta^{+-}_{i_1,j}]
=[\Delta^{+-}_{i,j}][\Delta^{++}_{i,i_1}]$.
If $i\ne i_1$ and $j\ne j_1$ we get similarly
\begin{eqnarray*}
\int_{S^d}[\Delta^{+-}_{i,j}][\Delta^{+-}_{i_1,j_1}]
p_{i+}^*a_1 p_{i_1+}^*a_2\!
\!\!\!\!\!\prod_{i_2\not\in \{i,i_1\}}\!\!\!\!\!
p_{i_2+}^*t_{1+}
\!\!\!\!\!\prod_{j_2\not\in \{j,j_1\}}\!\!\!\!\!p_{j_2-}^*t_{1-}\cdot w
\!\!&=&\! \int_{S^d}\!\!(a_1)_2(a_2)_2t_{1+}^{\st  (n-2)}
 t_{1-}^{\st  (m-2)}\cdot w.
\end{eqnarray*}
We can translate our result into this notation and simplify it by
collecting the terms number $2,4,8$ and  the terms
$3,5,9$ respectively. So we get for
$w\in H^{4d-2k}(S^d,\Q)^{\Symg_{d}}$ with $k\le 5$:
\begin{eqnarray*}
&&\int_{S^d}\tilde g_*\big(\phi^*\big((c(\Gamma)-1)
s([V^\vee(K_S)[n]^\vee)s(V[m])\big)\big)\cdot w\\
&&\quad =\int_{S^d}\Bigg(-nm
(2+12\xi+42\xi^2+3s_2(S)+K_S^2)_2
  t_{1+}^{\st  (n-1)}
  t_{1-}^{\st  (m-1)}\\
&&\qquad\quad
+{n\choose 2}m
(30+260\xi+114K_S)_3
  t_{1+}^{\st  (n-2)}
  t_{1-}^{\st  (m-1)}\\
&&\qquad\quad
+{m\choose 2}n(30+260\xi-114K_S)_3
  t_{1+}^{\st  (n-1)}
  t_{1-}^{\st  (m-2)}\\
&&\qquad\quad
+2{n\choose 2}{m\choose 2}
(2+12 \xi)_2^{  2}
  t_{1+}^{\st  (n-2)}  t_{1-}^{\st  (m-2)}\\
&&\qquad\quad
+mn{n-1\choose 2}
(2+12 \xi)_2
(3+18\xi+13K_S)_2  t_{1+}^{\st  (n-3)}
  t_{1-}^{\st  (m-1)}\\
&&\qquad\quad
+nm{m-1\choose 2}
(2+12 \xi)_2
(3+18\xi-13K_S)_2  t_{1-}^{\st  (m-3)}
  t_{1+}^{\st  (n-1)}\Bigg)\cdot w
\end{eqnarray*}
Now we  sum over all $m,n$
and keep in mind that the map
$\phi:(\widetilde S^{n}\times \widetilde S^{m})_{*}\maps
\Hilb^{n}(S)\times \Hilb^{m}(S)$ has degree $m! n!$. So  we obtain
\begin{eqnarray*}
&&d!\int_{\Hilb^d(S\sqcup S)}
(c(\Gamma)-1)s([V^\vee(K_S)]_1^\vee)s([V]_2)\cdot\gamma\\
&&\qquad
=\sum_{n+m=d} {d\choose n}
\int_{S^d}\tilde
g_*\big(\phi^*\big((c(\Gamma)-1)s([V^\vee(K_S)[n]^\vee)s(V[m])\big)\big)
\cdot w\\
&&\qquad
=\sum_{n+m=d} \int_{S^d}\Bigg(- d(d-1)  {d-2\choose n-1}
(2+12\xi+42\xi^2+3s_2(S)+K_S^2)_2
  t_{1+}^{\st  (n-1)}
  t_{1-}^{\st  (m-1)}\\
&&\qquad\quad +d(d-1)(d-2)
{d-3\choose n-2}(30+260\xi)_3
  t_{1+}^{\st  (n-2)}
  t_{1-}^{\st  (m-1)}\\
&&\qquad\quad
+2{d(d-1)(d-2)(d-3)}{d-4\choose n-2}
(2+12 \xi)_2^{  2}
  t_{1+}^{\st  (n-2)}
  t_{1-}^{\st  (m-2)}\Bigg)\cdot w\\
&&\qquad
=\int_{S^d}\Big(-d(d-1)(2+12\xi+42\xi^2+3s_2(S)+K_S^2)_2  t_1^{\st (d-2)}\\
&&\qquad\quad +d(d-1)(d-2)(30+260\xi)_3  t_1^{\st (d-3)}\\
&&\qquad\quad +2{d(d-1)(d-2)(d-3)}(2+12 \xi)_2^{  2}  t_1^{\st (d-4)}\Big)\cdot
w.
\end{eqnarray*}
This shows the first formula. The second follows by combining this formula with
proposition  \ref{erstint}.
\end{pf}

Now we have described the intersection numbers
$\int_T s(\Ext^1_q(\I_{Z_1},\I_{Z_2}\otimes p^*V)\cdot\gamma$,
and are in a position to
finish our computation of the leading terms of the change of the Donaldson
invariants $\delta_{l,r}(\alpha)$.
We first want to compute a formula for the change of $\delta_{N,0}(\alpha)$ and
then compute how one has to modify this formula to get  $\delta_{l,r}(\alpha)$.
The reason that the computation of $\delta_{N,0}(\alpha)$ is easier,
is the following fact:

\begin{rem} \label{pteinszwei}
Let $l,j,k$ be  positive integers, $\alpha\in H^2(S,\Q)$,
$\beta\in H^{2i}(S,\Q)$ and
$\gamma\in H^*(S^k,\Q)^{\Symg_{k}}$. Then we get
\begin{eqnarray}\label{pteinszw1}
\qquad\qquad\int_{S^{k+j}} (\beta)_j\st \gamma\cdot
(p_1^*\alpha+\ldots +p_{k+j}^*\alpha)^l= j^{2-i}\int_{S^{k+j}}
\beta\st pt^{\st (j-1)}\st \gamma\cdot
(p_1^*\alpha+\ldots +p_{k+j}^*\alpha)^l
\end{eqnarray}
\end{rem}

\begin{pf} For the diagonal $\Delta_{j}\subset S^j$
and a class $\alpha\in H^2(S,\Q)$, we have
$(p_1^*\alpha+\ldots +p_j^*\alpha)\cdot
[\Delta_{j}]=jp_1^*(\alpha)[\Delta_{j}]$.
By remark \ref{starrem} the left hand side of (\ref{pteinszw1}) is equal to
\begin{eqnarray*}
\left(\int_{S^{k}} \Delta_{j}p_1^*\beta\cdot
(p_1^*\alpha+\ldots +p_{k}^*\alpha)^{2-i}\right)\left(\int_{S^{j}}\gamma\cdot
(p_1^*\alpha+\ldots +p_{j}^*\alpha)^{l+i-2}\right).
\end{eqnarray*}
So the result follows.
\end{pf}

\begin{nota}
We denote by $q_S$ the quadratic form  on $H_2(S,\Z)$
and, for $\gamma\in H^2(S,\Q)$, we let  $L_{\gamma}$ be the linear form on
$H_2(S,\Q)$ given by
$\alpha\mapsto \<\gamma,\alpha\>$.
For a class $\beta\in H_i(S,\Q)$ we denote
$\bar \beta:=p_1^*\check \beta+ \ldots +p_d^*\check \beta\in
H^{4-i}(S^d,\Q)$, where
 as above, $\check\beta$ is the
Poincar\'e dual of $\beta$.
Note that by lemma \ref{hilbkuenn} and definition \ref{defchange}
$\phi^*(\widetilde\beta|_{\Hilb^{n}(S)\times \Hilb^{m}(S)})$ is the pullback of
$\bar \beta$.
Let $N=4c_2-c_1^2-3$ again  be the expected dimension of $M_H(c_1,c_2)$.
\end{nota}

\begin{lem}\label{formel1}
For all $x,y\ge 0$ and all $\alpha\in H_2(S,\Q)$ we have
\begin{eqnarray*}\int_{S^d}\xi^{\st x}  \pt^{\st y}  1^{\st (d-x-y)}
\cdot \bar\alpha^{2d-x-2y}={(2d-x-2y)!\over 2^{d-x-y}}
q_S(\alpha)^{d-x-y}\<\xi,\alpha\>^x.\end{eqnarray*}
\end{lem}
\begin{pf}
By remark  \ref{starrem} we have
\begin{eqnarray*}
\int_{S^d}\!\!\!\! \xi^{\st x}  \pt^{\st y}1^{\st(d-x-y)}
\cdot{\bar\alpha^{N-m}}=
\int_{S^d} \!\!\! p_1^{*}\xi\cdot \ldots \cdot p_x^{*}\xi
\cdot p_{x+1}^{*}\pt\cdot
\ldots\cdot p_{x+y}^{*}\pt
\cdot (p_1^*\check\alpha +\ldots +p_d^*\check\alpha)^{N-m},
\end{eqnarray*}
and it is easy to see that this is just
${(2d-x-2y)!\over 2^{d-x-y}}
q_S(\alpha)^{d-x-y}\<\xi,\alpha\>^x.$
\end{pf}

\begin{thm} \label{chthm1} In the polynomial ring on
$H^*(S,\Q)$ we have
$$\delta_{\xi,N,0}\equiv(-1)^{e_\xi}\sum_{k=0}^2
{N!\over (N-2d+2k)!(d-k)!}Q_{k}(N,d,K_S^2) L_{\xi/2}^{N-2d+2k}q_S^{d-k}\
\hbox{ modulo } L_\xi^{N-2d+6},$$ where, by convention
${1\over m!}=0$ for $m<0$ and
\begin{eqnarray*}
Q_{0}(N,d,K_S^2)&=& 1\\
Q_{1}(N,d,K_S^2)&=& 2N+2K_S^2-2d+8\\
Q_{2}(N,d,K_S^2)&=& 2N^2-4dN+4NK_S^2+21N+2d^2
-4dK_S^2 -18d+2(K_S^2)^2+18K_S^2+49.
\end{eqnarray*}
\end{thm}

\begin{pf}
Let $R_d\in H^*(S^d,\Q)^{\Symg_{d}}$ be the class from
proposition \ref{intzwei}
 with
$$d!\intt_{\Hilb^d(S\sqcup S)}
c(\Gamma)s([V^\vee(K_S)]_1^\vee)s([V]_2)\cdot \gamma=\intt_{S^d}
R_d\cdot {w}.$$
 By remark \ref{pteinszwei}
there is a class $U'_d$ which is a linear combination of classes of the form
$\xi^{\st x}\pt^{\st y} 1^{\st (d-x-y)}$ with
$\int_{S^d}
R_d\cdot\bar \alpha^b=\int_{S^d}
U'_d\cdot {\bar \alpha^b}$
for all $\alpha\in H_2(S,\Q)$.
We write $U'_d:=\sum_{x,y\ge 0}u_{x,y}\xi^{\st x}\pt^{\st y} 1^{\st (d-x-y)}$
and $U_d:=\sum_{x+y\le 2}u_{x,y}\xi^{\st x}\pt^{\st y} 1^{\st (d-x-y)}$
By definition  \ref{defchange} and theorem  \ref{donch1}
we see that
\begin{eqnarray*}
\delta_{\xi,N,0}(\alpha)=\sum_{i=0}^{2d}A_i\cdot\<\xi,\alpha\>\int_{S^d}
\{U'_d\}_i\bar\alpha^{2d-i},
\end{eqnarray*} where $\{\ \}_i$ denotes the part of degree $i$, and the
$A_i$ are suitable rational numbers.
Thus
$\delta_{\xi,N,0}$ modulo $L_\xi^{n-2d+6}$ is already determined by $U_d$.
As $S$ is a surface
with $p_g(S)=q(S)=0$, we have $12=12\chi(\oo_S)=K_S^2+c_2(S)$
and thus we can replace
$s_2(S)$ by $2K_S^2-12$.
So, using
 proposition  \ref{intzwei}, we obtain after a short calculation
that
\begin{eqnarray*}
U_d&=&
2^{d}1^{\st d}+2^{d+1}d \cdot1^{\st (d-1)}\st\xi
+2^{d+1}d(d-1)1^{\st (d-2)}  \xi^{\st 2}\\
&&\quad+2^d d(3\xi^2+K_S^2-5d+5)1^{\st (d-1)} \st  \pt\\
&&\quad+2^d d(d-1)(6\xi^2+2K_S^2-10d+5)1^{\st (d-2)}\st   \xi\st   \pt\\
&&\quad+2^{d-2} d(d-1)(18(\xi^2)^{2}+12\xi^2K_S^2+2(K_S^2)^{ 2}
-60d\xi^2-20dK_S^2+50d^2\\
&&\quad +15\xi^2-10K_S^2-34d-36)\pt^{\st 2}  1^{\st (d-2)},
\end{eqnarray*}
where we  view $\xi^2$ and $K_S^2$ as integers and not as cohomology classes.

Now we apply definition \ref{defchange} and lemma  \ref{formel1}.
Then, after some computation, we get the result with
$Q_0(N,d,K_S^2)$, $Q_1(N,d,K_S^2)$, $Q_2(N,d,K_S^2)$
replaced by
\begin{eqnarray*}
P_0(N,d,K_S^2,\xi^2)&=&1\\
P_1(N,d,K_S^2,\xi^2)&=&8N-26d+6\xi^2+2K_S^2+26\\
P_2(N,d,K_S^2,\xi^2)&=&18(\xi^2)^2+12(\xi^2)(K_S^2)+2(K_S^2)^2+48N\xi^2
-156d\xi^2\\
&&\quad-52dK_S^2+338d^2 +16K_S^2N+32N^2-208d
N+207\xi^2\\
&&\quad+54K_S^2+264N
-882d+508.
\end{eqnarray*}
We notice that by definition $d=(4c_2-c_1^2+\xi^2)/4$
and $N=4c_2-c_1^2-3$ and thus
$\xi^2=4d-N-3$.
Substituting  this into the $P_i(N,d,K_S^2,\xi^2)$ we obtain the result.
\end{pf}

We see that the result is compatible with the conjecture
of Kotschick and Morgan. In fact it suggests a slightly sharper statement.

\begin{conj}  In the polynomial ring on $H^2(S,\Q)$ we have
$$\delta_{\xi,N,0}=(-1)^{e_\xi}\sum_{k=0}^d
{N!\over (N-2d+2k)!(d-k)!}Q_{k}(N,d,K_S^2) L_{\xi/2}^{N-2d+2k}q_S^{d-k},$$
where $Q_{k}(N,d,K_S^2)$ is a  polynomial of degree
$k$ in $N,d,K_S^2$, which is independent of $S$ and $\xi$.
\end{conj}

Now we want to compute $\delta_{l,r}$ in general.
We shall see that there is reasonably simple relationship between the
formula for $\delta_{N,0}$ and that for $\delta_{l,r}$
(with $l+2r=N$),
which is however obscured  by the existence of
a correction term  coming from  the failure of
remark  \ref{pteinszwei}  for classes of the form
$\bar\alpha^{k-2}\bar \pt$ (instead of $\bar\alpha^{k}$).

\begin{lem}\label{formel2}
\begin{enumerate}
\item
For all $x,y\ge 0$, all $c\le r$   and all $\alpha\in H_2(S,\Q)$ we
have with $m:=2d-2c-x-2y$:
\begin{eqnarray*}
\int_{S^d}
  \xi^{\st x}  \pt^{\st y}1^{\st (d-x-y)}\cdot
 {\bar\pt^c\bar\alpha^{m}}={(d-x-y)!\over
(d-x-y-c)!}{m!\over 2^{d-x-y-c}}q_S(\alpha)^{d-x-y-c}
\<\xi,\alpha\>^{x}
\end{eqnarray*}
\item
\begin{eqnarray*}
\int_{S^d}(1)_2  1^{\st (d-2)}\cdot {\bar\pt\, \bar\alpha^{2d-4}}
={4d-6\over d-1}\int_{S^d}\pt\st   1^{\st (d-1)}\cdot
{\bar\pt\, \bar\alpha^{2d-4}}
\end{eqnarray*}
\end{enumerate}
\end{lem}
\begin{pf}
(1) By remark  \ref{starrem}  we have
\begin{eqnarray*}
&&\int_{S^d}\xi^{\st x}  \pt^{\st y}\cdot{\bar\pt^c\bar\alpha^{m}}
=
\int_{S^d} p_1^{*}\xi\ldots p_x^{*}\xi \cdot p_{x+1}^{*}\pt\ldots
p_{x+y}^{*}\pt
\cdot (p_1^*\check\alpha +\ldots +p_d^*\check\alpha)^{m}
\cdot (p_1^*\check\pt +\ldots +p_d^*\check\pt)^{c}
\end{eqnarray*}
and it is elementary to show that this is just
$$
{(d-x-y)!\over
(d-x-y-c)!}{m!\over 2^{d-x-y-c}}q_S(\alpha)^{d-x-y-c}
\<\xi,\alpha\>^{x}.
$$
(2) By remark \ref{starrem} and remark  \ref{pteinszwei} we have
\begin{eqnarray*}&&\int_{S^d} (1)_2  1^{\st (d-2)}\cdot{\bar\pt\,
\bar\alpha^{2d-4}}\\
&&\qquad
=\int_{S^d} [\Delta_{1,2}]
\cdot (p_1^* \check\alpha +\ldots +p_d^*\check\alpha)^{2d-4}
\cdot (p_1^*\check\pt +\ldots +p_d^*\check\pt)\\
&&\qquad=2\int_{S^{d-2}} (p_1^*\check\alpha  +\ldots
+p_{d-2}^*\check\alpha)^{2d-4}+(d-2)\int_{S^{d-1}} [\Delta_{1,2}]
\cdot (p_1^*\check\alpha +\ldots +p_d^*\check\alpha)^{2d-4}\\
&&\qquad=(4d-6)\int_{S^{d-2}}\bar\alpha^{2d-4}\\
&&\qquad=(4d-6)\int_{S^{d-1}}\pt*1^{\st (d-2)}\bar\alpha^{2d-4}\\
&&\qquad={4d-6\over
d-1}\int_{S^d}\pt\st   1^{ \st  (d-1)}\cdot{\bar\pt\,\bar\alpha^{2d-4}}
\end{eqnarray*}
\end{pf}

\begin{thm} \label{chthm2}
Let $l,r$ be nonnegative integers with $l+2r=N$.
Then in the polynomial ring on $H^*(S,\Q)$ we get
$$
\delta_{\xi,l,r}\equiv\sum_{c=0}^2{(-1)^{r-c+{e_\xi}}\over 2^{-3c+2r}}
{r\choose c}\sum_{k=c}^{2}
{l!\over (l-2d+2k)!(d-k)!} Q_{k-c,c}(l,d,K_S^2,\xi^2)
L_{\xi/2}^{l-2d+2k}q_S^{d-k}$$  modulo $ \xi^{N-2d+6},$
where
\begin{eqnarray*}
Q_{m,c}(l,d,K_S^2,\xi^2)&=& P_{m}(l,d,K_S^2,\xi^2)+21mc\ \hbox{ for }
m+c\le 2
\end{eqnarray*}
Here the  $P_{i}(N,d,K_S^2,\xi^2)$ are the polynomials from the
proof of theorem \ref{chthm1}.
\end{thm}

\begin{pf}
For  $i\le r$ and a class $\gamma\in H^*(S^d,\Q)^{\Symg_{d}}$
we denote by  $W_{l,r,c}(\gamma)$ the map that associates to
$\alpha\in H_2(S,\Q)$ the number
\begin{eqnarray*}
\sum_{b=0}^{l}(-1)^{r-c+e_\xi}2^{b+2c-N}{l\choose
b}{r\choose c}\<\xi,\al\>^{l-b}
\int_{S^d}\gamma\cdot{\bar\al^b\bar\pt^c}.
\end{eqnarray*}
Let $R_d$, $U_d\in H^*(S^d,\Q)^{\Symg_d}$ be the classes from the proof of
\ref{chthm1}.
By thm \ref{donch1} and proposition  \ref{intzwei} we get
$$\delta_{l,r}\equiv\sum_{c=0}^r
W_{l,r,c}(R_d).$$
By lemma  \ref{pteinszwei} we see that $W_{l,r,0}(R_d)\equiv W_{l,r,0}(U_d)$
 modulo $L_\xi^{N-2d+6}$.
Furthermore  we get  modulo $L_\xi^{N-2d+6}$
\begin{eqnarray*}
 W_{l,r,k}(R_d)&\equiv& 0 \ \hbox{ for } k>2\\
W_{l,r,2}(R_d)&\equiv & W_{l,r,2}(2^d 1^{\st d})\\
W_{l,r,1}(R_d)&\equiv& W_{l,r,1}(\bar R_d),
\end{eqnarray*}
where $\bar R_d=t_1^{\st d}-5\cdot 2^{d-2}d(d-1)(1)_21^{\st (d-2)}$.
By lemma  \ref{formel2}(2)
$W_{l,r,1}(\bar R_d)\equiv W_{l,r,1}(\bar U_d)$ where
$$\bar U_d=
t_1^{\st d}-5\cdot 2^{d-2}d(4d-6)\pt\st 1^{\st(d-1)}.$$
Now the result follows by applying lemma  \ref{formel2}(1)
and some computation.
\end{pf}
\begin{rem}
Using $\xi^2=4d-2r-l-3$ we get equivalently
\begin{eqnarray*}
Q_{0,0}(l,d,K_S^2,\xi^2)&=&Q_{0,1}(l,d,K_S^2,\xi^2)=
Q_{0,2}(l,d,K_S^2,\xi^2)=1\\
Q_{1,0}(l,d,K_S^2,\xi^2)&=&2l-2d-12r+2K_S^2+8\\
Q_{2,0}(l,d,K_S^2,\xi^2)&=&72r^2-24rl+24dr-24K_S^2r+2l^2-4dl+4K_S^2l\\
&&\quad +2d^2-4dK_S^2+2(K_S^2)^2-198r+21l-18d+18K_S^2+49\\
Q_{1,1}(l,d,K_S^2,\xi^2)&=&2 l - 2 d - 12 r +2K_S^2+ 29.\\
\end{eqnarray*}
\end{rem}

\begin{rem} We see that our results contain as a special case
the formulas for the change for $d\le 2$.
In the case that $d=3$ we notice that the spaces $X_*$ and
$Y_{*}$ (for $X$ and $Y$  schemes with a natural morphism
to $S^{(d)}$ and $S^{(n)}\times
S^{(m)}$ respectively)
just coincide with $X$ respectively $Y$. Therefore our computations are valid
on the whole of $\Hilb^d(S\sqcup S)$ and  our methods will also give complete
formulas for the change of the Donaldson invariants in case $d=3$. We
however do not carry out the elementary but long computations here.
\end{rem}


\begin{thebibliography}{ABCD}

\bibitem[A-K]{A-K} Altman, A., Kleiman, S.,
{\em Compactifying the Picard Scheme,}
Adv. Math. {\bf 35} (1980), 50--112.


\bibitem[BPS]{BPS} Banica, C., Putinar, M., Schumacher, G.,
{\em Variation der globalen Ext in Deformationen kompakter komplexer
R\"aumen}, Math.\ Ann.\ {\bf 250} (1980), 135--155.


\bibitem[Do]{Do} Donaldson, S.K., {\em Polynomial invariants for smooth
four-manifolds}, Topology {\bf 29} (1990), 257--315.

\bibitem[E1]{E1} Ellingsrud, G., Another proof of the
irreducibility of the punctual Hilbert scheme of a
smooth surface, Bergen Univ. preprint 1991.


\bibitem[D-H]{D-H} Dolgachev, I., Hu, Y., {\em Variation of geometric
invariant theory quotients}, preprint 1993.

\bibitem[Fo]{Fo} Fogarty, J.,{\em Algebraic families on an algebraic
surface}, Amer. J. Math. {\bf 90} (1968), 511-521


\bibitem[F-G]{F-G} Fantechi, B., G\"ottsche, L., {\em
The cohomology ring of the Hilbert scheme of 3 points on a smooth
projective variety}, J. reine angew. Math. {\bf 439} (1993), 147-158.



\bibitem[F-M]{F-M} Friedman, R., Morgan, J., {\em On the diffeomorphism
types of certain algebraic surfaces, I}, J. Diff. Geom. {\bf 27} (1988),
297--369.

\bibitem[F-Q]{F-Q} Friedman, R., Qin, Z.,
{\em Flips of moduli spaces and transition formulas for Donaldson
    polynomial invariants of rational surfaces}, preprint.

\bibitem[Fu]{Fu} Fulton, W., {\em Intersection Theory}, Ergebnisse der
Mathematik und ihrer Grenzgebiete, Springer 1984.


\bibitem[Gi]{Gi}
Gieseker, D.,{\em On the moduli of vector bundles on an algebraic surface},
 Ann. Math. {\bf 106} (1977),  45--60.

\bibitem[G\"o1]{Go1} G\"ottsche, L.,
{\em Change of polarization and Hodge
numbers of moduli spaces of torsion free sheaves on surfaces},
 MPI-preprint 94-9 (1994).

\bibitem[G\"o2]{Go2} \bysame {\em Trisecant formulas for
smooth projective varieties},
in  Projective Geometry with Applications, Marcel Decker Inc. 1994.



\bibitem[Ha]{Ha} Hartshorne, R., {\em Algebraic Geometry}, Graduate Texts
in Math., Springer Verlag, 1977.


\bibitem[Jo]{Jo} Jouanolou, J.P., {\em Riemann-Roch sans d\'enominateurs},
Invent. Math. {\bf 11} (1970), 15-26.



\bibitem[Ko]{Ko} Kotschick, D., {\em $SO(3)$-invariants for $4$-manifolds
with $b^+=1$}, Proc. London Math. Soc. {\bf 63} (1991), 426--448.

\bibitem[K-M]{K-M} Kotschick, D., Morgan, J.,
{\em $SO(3)$-invariants for $4$-manifolds with $b^+=1$ II},
 J. Diff. Geom. {\bf 39} (1994), 433--456.

\bibitem[La]{La} Lange, H.,\ \em Universal Families of Extensions.

\em J.\ of Alg.\ 83 (1983), 101--112.

\bibitem [L-Q]{L-Q} Li, W.P., Qin, Z.,
{\em Lower-degree Donaldson polynomials of rational surfaces},
J. Alg. Geom. {\bf 2} (1993), 413--442.

\bibitem[Li]{Li} Li, J., {\em Algebraic geometric interpretation of
Donaldson's polynomial invariants}, J. Diff. Geom.
 {\bf 37} (1993), 417--465.

\bibitem [Ma1]{Ma1} Maruyama, M.,
{\em Moduli of stable sheaves I}, J. Math. Kyoto Univ.,
{\bf 17} (1977), 91--126.

\bibitem [Ma2]{Ma2} \bysame {\em Moduli of stable sheaves II},
J. Math. Kyoto Univ. {\bf 18} (1978), 557--614.

\bibitem[Ma-Yo]{Ma-Yo}
Maruyama, M., Yokogawa, K.,
{\em Moduli of  parabolic stable sheaves},
Math. Ann. {\bf 293} (1992), 77--99.

\bibitem[M-W]{M-W} Matsuki, K., Wentworth, R.,
{\em Mumford-Thaddeus principle on the moduli space
of vector bundles on an algebraic surface}, preprint 1994.

\bibitem[Mo]{Mo} Morgan, J., {\em Comparison of the
Donaldson polynomial invariants with their algebro-geomtric analogues},
Topology {\bf 32} (1993),
449--488.

\bibitem[Mu1]{Mu1}
 Mukai, S., {\em Symplectic structure of the moduli space
 of sheaves on an abelian or K3 surface},
Invent. Math. {\bf 77} (1984), 101--116.

\bibitem[Mu2]{Mu2}
\bysame {\em On the moduli space of bundles on K3 surfaces I},
Vector bundles on  Algebraic varieties, Proc.~Bombay 1984,
 Oxford U.~Press (1987), 241--413.

\bibitem[Mu3]{Mu3} \bysame,\ \em Moduli of vector bundles on K3 surfaces,
and symplectic manifolds. \em Sugaku Expositions 1(2) (1988), 138--174.

\bibitem[OG1]{OG1} O'Grady, K.G., {\em Algebro-geometric analogues of
Donaldson's polynomials}, Invent. Math. {\bf 107} (1992), 351--395.

\bibitem[OG2]{OG2} \bysame,
{\em Relations among Donaldson polynomials of certain algebraic surfaces},
preprint (1993).

\bibitem[P-S]{P-S} Peskine, C., Szpiro, L.,
{\em Dimension projective finie et cohomologie locale}, Publ. I.H.I.S.
{\bf 42} (1972), 47--119.



\bibitem[Q1]{Q1} Qin, Z.,
{\em Equivalence classes of polarizations and moduli spaces of sheaves},
J. Diff. Geom. {\bf 37} (1993), 397--413.

\bibitem[Q2]{Q2} \bysame,
{\em Moduli of stable sheaves on ruled surfaces and their Picard groups},
 J. reine angew. Math. {\bf 433} (1992), 201--219.

\bibitem[Q3]{Q3}\bysame,
{\em Complex structures on certain differentiable $4$-manifolds},
Topology {\bf 32} (1993), 551--566.

\bibitem[Se]{Se}Serre, J.P., {\em Alg\`ebre locale. Multiplicit\'es},
Lecture notes in Math. {\bf 11}, Springer Verlag, 1965.


\bibitem[Th1]{Th1} Thaddeus, M.,
{\em Stable pairs, linear systems and the Verlinde formula},
Invent. Math. {\bf 117} (1994), 317--353.

\bibitem[Th2]{Th2} \bysame,
{\em Geometric invariant theory and flips},
preprint (1994).

\bibitem[Yo]{Yo}
 Yokogawa, K., {\em  Compactification of moduli of parabolic
sheaves and moduli for parabolic Higgs sheaves}, J. Math.
Kyoto University
{\bf 33} (1993), 451--504 .
\end{thebibliography}
\end{document}